\documentclass[aap,MSNbibl,dvips]{arximspdf}
\usepackage{newsym2e}

\input xypic
\xyoption{curve}

\doi{10.1214/13-AAP987} 
\volume{24}
\issue{6}
\pubyear{2014}
\firstpage{2595}
\lastpage{2643}

\makeatletter
\newcommand{\widebar}{\overline}
\newcommand{\bigtimes}{\mathop{\!\mbox{\parbox[c][9pt][b]{18pt}{\fontsize{18}{18}\selectfont{$\times$}}}\!\!\!\!}}
\newcommand{\xrightarrow}[1]{\stackrel{#1}{\longrightarrow}} 
\newtheorem{teo}{Theorem}[section]
\newtheorem{prop}[teo]{Proposition}
\newtheorem{lem}[teo]{Lemma}
\newtheorem{cor}[teo]{Corollary}
\newproclaim{defn}[teo]{Definition}
\newproclaim{rem}[teo]{Remark}
\makeatother

\begin{document}
\begin{frontmatter}

\title{A class of nonergodic interacting particle systems with unique
invariant measure\thanksref{T1}}
\runtitle{Nonergodic interacting particle systems}

\begin{aug}
\author{\fnms{Benedikt} \snm{Jahnel}\ead[label=e1]{Benedikt.Jahnel@ruhr-uni-bochum.de}}
\and
\author{\fnms{Christof} \snm{K\"ulske}\corref{}\ead[label=e2]{Christof.Kuelske@ruhr-uni-bochum.de}}
\runauthor{B. Jahnel and C. K\"ulske}
\affiliation{Ruhr-Universit\"at Bochum}
\dedicated{Dedicated to A. van Enter on the occasion of his 60th birthday}
\address{Fakult\"at f\"ur Mathematik\\
Ruhr-Universit\"at Bochum\\
D44801 Bochum\\
Germany\\
\printead{e1}\\
\phantom{E-mail:\ }\printead*{e2}}
\end{aug}
\thankstext{T1}{Supported by the Sonderforschungsbereich
SFB|TR12-Symmetries and Universality in Mesoscopic Systems.}

\received{\smonth{12} \syear{2012}}
\revised{\smonth{10} \syear{2013}}

%
\begin{abstract}
We consider a class of discrete $q$-state spin models defined in terms
of a translation-invariant quasilocal specification with discrete clock-rotation invariance which
have extremal Gibbs measures $\mu'_{\varphi }$ labeled by the uncountably
many values of
$\varphi $ in the one-dimensional sphere (introduced by van Enter, Opoku,
K\"ulske [\textit{J. Phys. A} \textbf{44} (2011) 475002, 11]).
In the present paper we construct an associated Markov jump process
with quasilocal rates whose semigroup $(S_t)_{t\geq0}$
acts by a continuous rotation $S_t(\mu'_{\varphi })=\mu'_{\varphi +t}$.

As a consequence our construction provides examples of interacting
particle systems with unique translation-invariant invariant measure, which
is not long-time limit of all starting measures, answering an old
question (compare Liggett
[\textit{Interacting Particle Systems} (1985) Springer], question
four, Chapter one).
The construction of this particle system is inspired by recent
conjectures of Maes and Shlosman about
the intermediate temperature regime of the nearest-neighbor clock model.
We define our generator of the interacting particle system as a
(noncommuting) sum
of the rotation part and a Glauber part.

Technically the paper rests on the control of the spread of weak nonlocalities
and relative entropy-methods, both in equilibrium and dynamically,
based on Dobrushin-uniqueness bounds for conditional measures.
\end{abstract}

%

\end{frontmatter}

\section{Introduction}\label{sec1}

Consider an interacting particle system (IPS) on the infinite
$d$-dimensional integer lattice with finite local state space and
quasilocal rates. Such an IPS is a Markov process in continuous time
where particles (or spins) which sit on the lattice sites taking one of
finitely many spin values are updated
after exponential waiting times to take new states with probabilities
which depend
in an (essentially) local way on the states of the neighboring
particles. Assume that these updating
rules are lattice translation-invariant. Such infinite-volume processes
may possess multiple equilibria (time-invariant measures).
Indeed, any Gibbsian potential (Hamiltonian) for a discrete-spin model
allows one to prescribe rates defining
a Glauber dynamics for which the corresponding Gibbs measures are
time-invariant and moreover reversible. Consequently, if there is a
phase-transition (meaning that there is nonuniqueness of the Gibbs measures
for this Hamiltonian),
the set of time-invariant measures has more than one point; see \cite
{Li85}. To prove on the other hand that for a Glauber dynamics
there are no time-invariant measures
other than Gibbs measures is more intricate, and in general dimensions
this statement is only known to be true if one assumes all measures to
be lattice-translation invariant; see \cite{Ho71,Li85} and compare
Proposition~\ref{TheGlauberDynamics}.

To pose our problem let us start now from any lattice
translation-invariant IPS without assumptions on reversibility.
Consider a lattice translation-invariant measure which is
invariant under the IPS dynamics. Suppose there is only one such measure.
Is it true that the dynamics is necessarily ergodic? The notion of
ergodicity for an IPS
means that for any starting measure the time-evolved measures converge to
the unique invariant measure.

This is an old question which was picked up again in a recent very
interesting paper by Maes and Shlosman \cite{MaSh11} about dynamics of
clock models; see \mbox{\cite{FrSp81,FrSp82,BrRePl10}} and \cite{NeSc82}.
In their paper the authors conjecture that this may not be the case and
suggest a mechanism producing
time-periodic behavior
of rotating infinite-volume states.
The concrete model they suggest to analyze is the discrete rotator model
with standard scalarproduct nearest-neighbor interactions at
intermediate temperatures,
and a nonreversible time-evolution. Nonergodicity could appear
because if one uses one of the Gibbs measures as the initial measure,
the discrete rotators would keep rotating coherently,
and so the starting distribution would be repeated periodically under
the dynamics.
While these conjectures seemed plausible, at the same time no simple
proof based on their heuristics in their model
seemed possible.

To see naively how periodicity can create nonergodicity think of the
example of a two-state
discrete time Markov chain with
transition matrix $
{0\ \ 1 \choose 1\ \ 0}$. This chain has the unique invariant distribution
$(\frac{1}{2},\frac{1}{2})$, but obviously never forgets its initial
condition.
The same phenomenon of a unique stationary measure which does not
attract all starting measures
occurs for a Markov chain if the state space is finite and transition
graph is bipartite.

Can such a periodic behavior with unique invariant measure
persist for Markov processes with time-simultaneous updating of all
spins with local rules
on the infinite-lattice? Yes, and an example for nonergodicity of
discrete-time,
parallel updating PCA (probabilistic cellular automaton),
was only recently given in \cite{ChMa11}.
However, the issue of existence of a nonergodic IPS which interests us here
can not immediately be reduced to that of a nonergodic PCA.
Indeed, continuous Markovian time-evolutions in comparison
with discrete time-evolution have a tendency to wash out synchronization
and forget initial conditions. (Clearly the continuous time version of the
simple two-state Markov chain example mentioned above \textit{is} ergodic.)

In the present paper we construct a dynamics for a $q$-state particle
system ($q$~possibly large
but finite) which does the job: it has a unique translation-invariant
invariant measure for which the dynamics is not ergodic.
Our construction is inspired by the conjectures of Maes and Shlosman
(and different from \cite{ChMa11}) which we put to a situation where
they can be proved.

In order to do this we will relate an IPS to a hidden system of
continuous \mbox{$S^1$-}valued spins via a discretization transformation
which acts on each local state space. This will allow us to carry over
knowledge about phase transitions in the continuous system to the discrete
system we want to analyze. Technically it builds on earlier works \cite
{EnKuOp11,KuOp08} about the preservation
of the Gibbs property under such discretization transformations. While
these results concern properties
of equilibrium measures the main new idea of the present paper is the
definition of an associated
nonreversible dynamics.
This dynamics is chosen in such a way that it preserves the set of
equilibrium measures. It does not (unlike a Glauber
dynamics) preserve the individual equilibrium measures but rotates the
lattice translation-invariant equilibrium measures
into each other periodically. In this way a periodic orbit of measures
is constructed.
That such a dynamics can be realized by means of a generator with
quasilocal jump rates is one
main result of this paper; that this dynamics has a unique
time-invariant translation-invariant measure is
another main result.

The interest in the study of rotation dynamics also has an independent
source which comes
from biological applications like interacting neurons or collective
behavior of animals.
Usually the models studied in this context
are of mean-field type like the famous Kuramoto model. This is natural
from the
perspective of many applications and also has the technical advantage
of reducing all relevant questions
to questions about (paths of) empirical distributions which makes them
more tractable than lattice systems.
In these models one usually studies $S^1$-valued spins under diffusive
time-evolutions which
contain a mean-field coupling that tends to synchronize the rotators.
Often these models contain additional sources of quenched randomness
(modeling individual rotation frequencies) which lead to
a nonreversible character and a periodic orbit which is deformed in a
way which depends
on fluctuations of the realizations of the rotation frequencies.
The relevant questions starting with
existence of synchronized rotating states and their finer
properties have been very successfully studied in particular in the
Kuramoto model \cite{AcBoPeRiSp05,CoDa12,GiPaPePo}.

Viewed in this light our construction of a lattice dynamics
hints at the existence of synchronization phenomena also on the
lattice, even for discrete local spaces.
It would be interesting to know more about the domain of attraction of
the periodic orbit whose existence we prove, but
we do not tackle this ongoing issue in this paper where we only analyze
properties \textit{on} the cycle.
Let us mention in this context that our construction of a rotation
dynamics to implement the Maes--Shlosman
mechanism of nonergodic behavior can be carried over to a mean-field setup.
We perform the construction of such a dynamics and the analysis of its
properties
in the related paper \cite{JaKu13}.
In that paper synchronization for discrete rotators is actually proved,
and a Lyapunov function
is constructed to prove attractivity of the cycle of rotating Gibbs measures.

\subsection{Main result}\label{sec1.1}
To construct our IPS we have to introduce
a continuous-spin model first which will be given in terms of a
Gibbsian specification for
an absolutely summable Hamiltonian (energy function) acting on
continuous spins.
The particle dynamics will be related to this model in a further step.
To define this continuous-spin model we consider an $S^1$-rotation
invariant and translation-invariant Gibbsian specification $\gamma
^\Phi$ on
the lattice $G={\mathbb Z}^d$, with
local state space $S^1=[0,2\pi)$. Let this specification $\gamma
^\Phi=(\gamma
^\Phi_\Lambda )_{\Lambda \subset G}$
be given in the standard way by an absolutely summable, $S^1$-invariant
and translation-invariant potential
$\Phi=(\Phi_A)_{A{\subset}G, A\ \mathrm{finite}}$, w.r.t. to the Lebesgue
measure $\lambda $ on the spheres.
This means that the Gibbsian specification is given by the family of
probability kernels
%
\begin{equation}
\label{FirstLayerSpecification*} \gamma ^\Phi_\Lambda (B|\eta)=\frac{\int1_B(\omega _\Lambda
\eta_{\Lambda ^c}) \exp(-H_\Lambda (\omega _\Lambda \eta
_{\Lambda ^c}))\lambda ^{\otimes\Lambda }(d\omega _\Lambda
)}{\int\exp(-H_\Lambda (\omega _\Lambda \eta_{\Lambda
^c}))\lambda
^{\otimes\Lambda }(d\omega _\Lambda )}
\end{equation}
for finite $\Lambda {\subset}G$ and Hamiltonian $H_\Lambda =\sum_{A\cap\Lambda \neq\varnothing
}\Phi_A$ applied to a measurable set $B{\subset}(S^1)^G$ and a boundary
condition $\eta\in(S^1)^G$; for details on Gibbsian specifications, see
\cite{Ge11}. We use notation $\Lambda ^c:=G\setminus\Lambda $.
$H_\Lambda $ also has to
be differentiable under variation at a single site and these partial
derivatives have to be uniformly bounded. A standard example of such a
model is provided by the nearest-neighbor scalarproduct interaction
rotator model with Hamiltonian
%
\begin{eqnarray}\label{MetricFamily}
&& H_{\Lambda}(\omega _\Lambda \eta_{\Lambda ^c})
\nonumber\\[-8pt]\\[-8pt]
&&\qquad  =
- \beta\sum_{i,j\in\Lambda\dvtx \langle  i,j \rangle
}\cos(\omega _i-\omega
_j) - \beta\sum_{i\in\Lambda,j \in\Lambda^c\dvtx  \langle i,j \rangle} \cos(\omega
_i-\eta_j).\nonumber
\end{eqnarray}

Denote by $\mathcal{G}(\gamma ^\Phi)$ the simplex of the Gibbs measures
corresponding to this specification,
which are the probability
measures $\mu$ on $(S^1)^G$ which satisfy the \mbox{DLR-}equation $\int\mu
(d\eta)
\gamma ^\Phi_\Lambda (B|\eta)=\mu(B)$ for all finite $\Lambda $.
Denote by $\mathcal{G}_{\theta
}(\gamma ^\Phi)$ the lattice
translation-invariant Gibbs measures.

We will make as an assumption on the class of potentials (Hamiltonians)
we discuss moreover that it has a continuous symmetry breaking in the
following sense. Assume that the extremal translation-invariant Gibbs
measures can be obtained as weak limits with
homogeneous boundary conditions, that is, with $\eta_\varphi \in(S^1)^G$
defined as $(\eta_\varphi )_i=\varphi $ for all $i\in G$ and
$\varphi \in S^1$ we
have
\[
\operatorname{ex} \mathcal{G}_{\theta}\bigl(\gamma ^\Phi
\bigr)= \Bigl\{ \mu _\varphi \big|\mu_\varphi = \lim
_{\Lambda
\nearrow G} \gamma ^{\Phi}_{\Lambda }(\cdot|
\eta_\varphi ), \varphi \in S^1 \Bigr\}.
\]
We further assume that
different boundary conditions $\eta_\varphi $ yield different
measures so
that there is a unique labeling of states $ \mu_\varphi $
by the angles $\varphi $ in the sphere $S^1$. It is a nontrivial proven
fact that
this assumption is true in the case of the standard rotator model (\ref
{MetricFamily}) in $d=3$ for $\lambda $-a.a. temperatures in the
low-temperature region as discussed in \cite{FrPf83,MaSh11,FrSiSp76,Pf82}.

We will now describe the discretization transformation which maps the
continuous-spin model
to a discrete-spin (or particle) model on which then the dynamics will
be constructed in the following step.

Denote by $T$ the local coarse-graining with equal arcs, that is,
$T\dvtx [0,2\pi)\mapsto\{1,\ldots, q\}$ where $T(\varphi ):=k$ if and
only if
$2\pi(k-1)/q\leq\varphi <2\pi k/q$.
Extend this map to infinite-volume configurations
by performing it sitewise. We will refer to the image space $\{1,\ldots,q\}^G$ as the
coarse-grained layer.
In particular we will consider images of infinite-volume measures under $T$.

We will need to choose the parameter of this discretization $q\geq
q_0(\Phi)$ large enough
so that the image measures are again Gibbs measures for a discrete
specification on the coarse-grained layer.
That this is always possible (even for large interactions)
follows from our earlier investigations \cite{KuOp08,EnKuOp11}.
More precisely, let us assume that the condition from Theorem~2.1 of
\cite{EnKuOp11} is fulfilled
(ensuring a regime where the Dobrushin uniqueness condition holds
for the so-called constrained first-layer models---the Dobrushin
condition is a weak dependence condition implying uniqueness and
locality properties).
Note, as in our notation, the usual temperature parameter $\beta $ is
incorporated into $\Phi$, for $\beta $ tending to infinity so does
$q_0(\Phi)$.

We are now ready to describe our definition of a dynamics on the
coarse-grained layer
in terms of a generator
which plays well together with the discretization transformation $T$
just introduced. This dynamics has two parts, a reversible part and a
nonreversible part.
We begin with the more interesting nonreversible part and
define a Markov process with state space $\{1,\ldots,q\}^G$ in terms of
the generator
%
\begin{equation}
\label{rotationgenerator} (L\psi) \bigl(\omega ' \bigr):=\sum
_{i\in G} c_L \bigl(\omega ', \bigl(
\omega ' \bigr)^i \bigr) \bigl(\psi \bigl( \bigl(\omega
' \bigr)^i \bigr)-\psi \bigl(\omega '
\bigr) \bigr)
\end{equation}
acting on sufficiently smooth observables $\psi$. The jump rates are
given in terms of
certain expectations of conditional infinite-volume measures which
naturally arise
in the course of the discretization transformation.

The choice of these rates may not seem intuitive at this stage, but they
can be obtained heuristically from a straightforward computation, as we will
explain at a later stage, namely (\ref{Intuitive}). Let us at this
stage just describe their definition which is
%
%
%
%
%
%
%
%
\begin{eqnarray}
\label{rotationjumprates} \qquad c_L \bigl(\omega ', \bigl(\omega
' \bigr)^i \bigr)&:=&\frac{\int\mu_{G\setminus i}[\omega
'_{G\setminus
i}](d\omega _{G\setminus i})e^{-H_{i}(2\pi\omega '_i/q,\omega
_{G\setminus i})}
}{\int\mu_{G\setminus i}[\omega '_{G\setminus i}](d\omega
_{G\setminus i})\int\lambda (d\omega
_{i})e^{-H_{i}(\omega _i,\omega _{G\setminus i})}
1_{T(\omega _i)=\omega '_{i}}
}
\nonumber\\[-8pt]\\[-8pt]
&=&\frac{\mu_{G\setminus i}[\omega '_{G\setminus
i}](e^{-H_{i}(\omega '_i|^r,\cdot_{i^c})})
}{\mu_{G\setminus i}[\omega '_{G\setminus i}](\lambda ^i ( e^{-H_{i}}
1_{\omega '_{i}}) )
},\nonumber
\end{eqnarray}
where we have written the expression on the first line for clarity, and
the second line is the short notation we will continue to use.
Further we used the following notation: $(\omega ')^i$ is the discrete
configuration which coincides with $\omega '$ except at the site $i$ where
it is increased by the amount
of one unit (modulo $q$). The continuous
spin value $\omega '_i|^r:=2\pi\omega '_i/q\in S^1$ is the right
endpoint of the
interval in continuous single-spin space at the site
$i$ prescribed by $\omega '_i$.
(In other words, in the definition of the rate to jump up at site $i$
from $\omega '_i$ to $\omega '_i+1$, the Hamiltonian appearing
under the integral gets evaluated right at the continuous-spin boundary
$\omega '_i|^r$ between the segments of $S^1$ labeled by $\omega '_i$
and $\omega
'_i+1$.) Finally, the measure $\mu_{G\setminus i}[\omega
'_{G\setminus
i}]$ is the unique continuous-spin Gibbs measure for a system on the
smaller volume $G\setminus i$ with conditional specification obtained by
deleting all interactions with $i$ and constrained to take values
$\omega
_{G\setminus i}$ with discretization
images $T(\omega _{G\setminus i})=\omega '_{G\setminus i}$.
For more details and precise definition of $\mu_{G\setminus i}[\omega
'_{G\setminus i}]$ in terms of formulas,
see Section~\ref{2}, namely (\ref{UniqueConditionalGibbsMeasure}).
Note that these constrained Gibbs measures are well defined and well
behaved for sufficiently fine discretization $q\geq q_0(\Phi)$, see
\cite{EnKuOp11,KuOp08}
and Section~\ref{2}. For general background on constrained Gibbs
measures in the context of preservation of Gibbsianness, see \cite
{EnFeSo93,Fe05} and \cite{KuLeRe04}.

From the definition of the rates (\ref{rotationjumprates}) it is clear
that the corresponding dynamics will be irreversible since jumps are
only possible in one direction. Note that these rates depend on the original
continuous-spin Hamiltonian in two places, namely in the $H_i$ and in
the $\mu_{G\setminus i}$.

Having defined the nonreversible part of our dynamics,
we next consider a Glauber-type generator $K$ on the same space
$\{1,\ldots, q\}^G$ by putting
%
\begin{eqnarray}
\label{GlauberGenerator} (K\psi) \bigl(\omega ' \bigr)&:=&\sum
_{i\in G}  \bigl[c_K \bigl(\omega ',
\bigl(\omega ' \bigr)^i \bigr) \bigl(\psi \bigl( \bigl(
\omega ' \bigr)^i \bigr)-\psi \bigl(\omega
' \bigr) \bigr)
\nonumber\\[-8pt]\\[-8pt]
&&\hspace*{17pt}{} +c_K \bigl(\omega ', \bigl(\omega
' \bigr)^{i-} \bigr) \bigl(\psi \bigl( \bigl(\omega
' \bigr)^{i-} \bigr)-\psi \bigl(\omega '
\bigr) \bigr) \bigr]\nonumber
\end{eqnarray}
with $(\omega ')^{i-}$ being the discrete configuration which
coincides with
$\omega '$ except at the site $i$ where it is decreased by the amount
of one
unit. We choose the rates to go up and down, respectively, such that they
satisfy
%
\begin{equation}
\label{GlauberRates} %
\frac{c_K(\omega ',(\omega ')^i)}{c_K((\omega ')^i,\omega ')} 
=\frac{\mu_{G\setminus i}[\omega '_{G\setminus i}](\lambda ^{i}
(e^{-H_{i}}
1_{(\omega '_i)^i}))}{\mu_{G\setminus i}[\omega '_{G\setminus
i}](\lambda ^{i} (e^{-H_{i}}
1_{\omega '_{i}}))}. %
\end{equation}
(For clarity of notation we note that, e.g., the denominator on the RHS means
$\mu_{G\setminus i}[\omega '_{G\setminus i}](\lambda ^{i} (e^{-H_{i}}
1_{\omega '_{i}}))=
\int\mu_{G\setminus i}[\omega '_{G\setminus i}](d\omega
_{G\setminus i})\int\lambda (d\omega
_{i})e^{-H_{i}(\omega _i,\omega _{G\setminus i})}\*1_{T(\omega
_i)=\omega '_{i}}$.)
A possible choice of $K$ is obtained by identifying numerators (resp.,
denominators) on the RHS and LHS of (\ref{GlauberRates}).

Having defined the two generators $L$ and $K$, we are finally in the
position to formulate
our main result. We have the following theorem.
%
\begin{teo} \label{MainTheorem} Consider a translation-invariant,
rotation-invariant and continuously differentiable potential $\Phi$
which satisfies the decay assumption
%
\begin{equation}
\label{ExpDecayCondition} \sum_{A\ni0}\sum
_{k\in G}e^{\varepsilon |k|}\delta_k(\Phi
_A)<\infty
\end{equation}
for some $\varepsilon >0$ where $\delta _k(\Phi_A)=\sup_{\omega,\bar\omega\dvtx \omega _{k^c}=\bar\omega
_{k^c}}|\Phi_A(\omega )-\Phi_A(\bar\omega )|$
denotes the variation at the site $k$.
Assume fine enough discretization $q\geq q_0(\Phi)$, and let $\alpha
>0$ be
arbitrary.
\begin{longlist}[(3)]
\item[(1)] Then the generator $L+\alpha K$ gives rise to a
welldefined IPS
with quasilocal rates.
\item[(2)] The class of translation-invariant measures which are
invariant under the associated
Markov semigroup $(S^{L+\alpha K}_t)_{t\geq0}$ consists of a single element.
\item[(3)] There are translation-invariant measures which do not
converge under the dynamics to the unique invariant measure.
\end{longlist}
\end{teo}
Note that any finite-range potential or exponentially decaying
pair-potential satisfies (\ref{ExpDecayCondition}). We further note
that the requirements on the potential can be relaxed. For example, one
could replace exponential decay by polynomial decay of sufficiently
high order as will become clear from the proof. The conditions will be
presented whenever they get used for the first time.

\subsection{Idea of proof}\label{sec1.2}

The proof relies on the fact that the discretization transformation $T$
preserves
the Gibbsian structure of the continuous and discrete-spin system
if we assume fine enough discretization $q\geq q_0(\Phi)$, in the
following sense.

First, to talk about the correspondence between the continuous and the
discrete system we need
to make explicit the relevant Gibbsian specification for the latter.
To do so define a family of kernels $\gamma '=(\gamma '_\Lambda
)_{\Lambda \subset G, \Lambda \
\mathrm{finite}}$ for the discretized model by
%
\begin{eqnarray}
\label{CoarseSpecification} \gamma '_{\Lambda } \bigl(\omega
'_{\Lambda } | \omega ' _{G\setminus
\Lambda } \bigr)
&=& \frac{\mu_{G\setminus\Lambda }[\omega '_{G\setminus\Lambda
}](\lambda ^{\Lambda } (e^{-H_{\Lambda }}
1_{\omega '_{\Lambda }}) )
}{\mu_{G\setminus\Lambda }[\omega '_{G\setminus\Lambda
}](\lambda ^{\Lambda } (e^{-H_{\Lambda }}))},
\end{eqnarray}
where in analogy to the explanation for $\mu_{G\setminus i}[\omega
'_{G\setminus i}]$ given before, $\mu_{G\setminus\Lambda }[\omega
'_{G\setminus\Lambda }]$ is the unique continuous-spin Gibbs measure for
the continuous specification on the volume $G\setminus \Lambda $, not
interacting
with $\Lambda $ and conditioned to a discrete configuration $\omega
'_{G\setminus
\Lambda }\in\{1,\ldots,q\}^{G\setminus \Lambda }$. This $\gamma
'$ indeed is a quasilocal
specification, and the discretized Gibbs measures will be Gibbs for
this $\gamma '$. For details see Section~\ref{2}.

Further, the infinite-volume discretization map $T$ is injective when
applied to the set of translation-invariant extremal Gibbs states in
the continuum model ($\operatorname{ex} \mathcal{G}_{\theta}(\gamma
^\Phi)$).
More precisely we have the following theorem.
%
\begin{teo}\label{Bijection} $T$ is a bijection from $\operatorname{ex}
\mathcal{G}
_{\theta}(\gamma ^\Phi)$ to $\operatorname{ex} \mathcal{G}_{\theta
}(\gamma ')$
with inverse given by the kernel 
$\mu_G[\omega '](d\omega )$.
\end{teo}
%
Here $\mu_G[\omega '](d\omega )$ is the unique conditional
continuous-spin Gibbs
measure on the whole volume $G$; see (\ref{UniqueConditionalGibbsMeasure}).
It is important to understand that this kernel gets us back from a
discrete-spin Gibbs measure to a continuous-spin Gibbs measure in a way
which \textit{does not depend} on the choice of the initial measure. This
is crucial for the possibility
to construct a rotation generator $L$ with the desired properties, as
we will see.

The fact that $T \mu:=\mu\circ T^{-1}$ is Gibbs for $\gamma '$ when
$\mu$
is Gibbs for $\gamma ^\Phi$, is already proved in \cite
{EnKuOp11,KuOp08} and
based on the uniform Dobrushin condition
on the coarse-graining. The part that each translation-invariant
discrete Gibbs measure has a discretization preimage in the continuous
Gibbs measures
is new and uses the Gibbs variational principle which involves
considerations of
relative entropy densities; see \cite{Ge11}.

The following step of the proof presents the main new structure of our
paper. We show that rotation on the level of discrete extremal Gibbs
states $\mu'_{\varphi }=T\mu_{\varphi }$ can be realized
by the rotation dynamics with generator $L$ with quasilocal jump rates
as defined above. This can be formulated as follows.
%
\begin{teo}\label{TheRotationDynamics}
\textup{(1)} The\vspace*{1pt} semigroup $(S^{L}_t)_{t\geq0}$ associated to $L$ is
well defined. 

\textup{(2)} $S^{L}_t(T \mu_{\varphi })= T \mu_{\varphi +t}$ for all
$\mu_\varphi \in
\operatorname{ex} \mathcal{G}_{\theta}(\gamma ^\Phi)$ and $t\geq0$.
\end{teo}
%

The theorem expresses that a discretization of a deterministic rotation
of the continuous-spin model can be represented as a stochastic time
evolution after discretization.
The heuristic reason why this works and the heuristic route to the
identification of such a suitable
$L$ is explained in formula (\ref{Intuitive}):
the idea is to compute the time derivative $ \frac
{d}{dt}_{|_{t=0}}(T\mu
_{\varphi +t})(f)$ for indicator\vspace*{1.5pt} functions $f$,
and to identify the appearing terms as $(T\mu_\varphi )(Lf)$. During this
computation
one makes explicit the kernel from discrete
to continuous variables of Theorem~\ref{Bijection}, uses its properties
and the rates defining $L$ pop out.
If we already knew that the trajectory $t \mapsto T \mu_{\varphi +t}$ can
be realized in terms of a semigroup,
this would identify its generator.
A difficulty in the actual proof is that we do not know this a priori,
and more arguments are needed. This involves the definition of
weighted triple-norms (weighted sums of variations of observables)
to control the weak nonlocalities which are present in the rates and
the spreading of these
under the action of the dynamics.

Rephrasing the result in a group theoretical language, we can say
$(t,\mu_\varphi )\mapsto\mu_{\varphi +t}$ is an $S^1$-action on the
extremal translation-invariant Gibbs measures $\operatorname{ex}
\mathcal{G}
_{\theta
}(\gamma ^\Phi)$ and
$(t,\mu'_\varphi )\mapsto\mu'_{\varphi +t}$ is an $S^1$-action on
$\operatorname{ex} \mathcal{G}_{\theta}(\gamma ')$. The second
statement of the theorem then says that $T$ is an equivariant map
(i.e., a group-action preserving map).


Let us now turn to the discussion of the reversible generator $K$.
Having defined the discretized local specification $\gamma '=(\gamma
'_\Lambda )_{\Lambda
\subset G}$ we note that
the generator $K$ defined above plays the role of a corresponding
Glauber dynamics.
%
%
To understand the final arguments providing us with a unique
translation-invariant measure for the joint dynamics
and understand better this Glauber part of the dynamics we prove the
following intermediate result.
%
\begin{prop}\label{TheGlauberDynamics}
\textup{(1)} The semigroup $(S^{K}_t)_{t\geq0}$ associated to $K$ is
well defined.

\textup{(2)} The translation-invariant measures which are invariant under
the dynamics
$(S^{K}_t)_{t\geq0}$ are precisely the discrete Gibbs measures
$\mathcal{G}
_{\theta}(\gamma ')$.
\end{prop}

To see that invariance under this dynamics implies Gibbsianness we use
an adaptation of the relative entropy arguments exposed in Liggett
(``Holley's argument'') \cite{Ho71,Li85} from the Ising lattice gas context
to our situation. The standard idea here is to exploit the form of the
time derivatives
of relative entropy densities of the time-evolved measure relative to a
suitable finite-volume version
of a Gibbs measure.
Putting these to zero, along with translation-invariance and estimation
of boundary terms, produces a single-site DLR equation implying that
the invariant measures are Gibbs for $\gamma '$.

The technical treatment of this
beautiful argument will have to be substantially modified in view of
the new terms arising from the joint dynamics
corresponding to $L+ \alpha K$, which we want to consider finally. The
result is the following proposition which is essential for the proof of
the main theorem.
%
\begin{prop}\label{TheJointDynamics} Let $\alpha >0$.
\begin{longlist}[(3)]
\item[(1)] The semigroup $(S^{L+\alpha K}_t)_{t\geq0}$ associated to
$L+\alpha
K$ is well defined.\vadjust{\goodbreak}
\item[(2)] $S^{L+ \alpha K}_t(T \mu_{\varphi })= T \mu_{\varphi
+t}$ for all $\mu
_\varphi \in\operatorname{ex} \mathcal{G}_{\theta}(\gamma ^\Phi
)$ and $t\geq0$.
\item[(3)] The translation-invariant measures which are invariant under
the joint dynamics
$(S^{L+\alpha K}_t)_{t\geq0}$ must necessarily be elements of the discrete
Gibbs measures $\mathcal{G}_{\theta}(\gamma ')$.
\end{longlist}
\end{prop}
For the proof we use that the Glauber part leaves the discrete Gibbs
measures invariant.
Let us point out some of the issues which come into play. A bit of care
needs to be taken for the second statement since the rotation part $L$
and the Glauber part $K$ do not commute. However, one can follow the
same line of arguments as for the proof of Theorem~\ref{TheRotationDynamics}, part~(2), using weighted triple-norms, to control
the weak nonlocalities of $L$ and $K$. The idea of the third part is
this: to see that invariance under joint dynamics implies Gibbsianness
we would like to use again relative entropy arguments as in the proof
of Proposition~\ref{TheGlauberDynamics}, part (2), but note that we
now have to deal with a sum of two terms each corresponding to $L$ and
$K$. For the new part corresponding to $L$ we apply the arguments to a
finite-volume open boundary version of $L$ as well as of the measure in
the second slot of the relative entropy.
The correction term is only of boundary-order.
The bulk terms have a good sign by a finite-volume argument since the
modified $L$ is attractive to the modified measure. Together we arrive
at the desired single-site DLR equation.

Combining the second and the third part of Proposition~\ref
{TheJointDynamics} we conclude:

\begin{cor}\label{TheUniqueInvariantMeasureForJoint} Let $\alpha >0$.
Then the only translation-invariant measure which is invariant under
the joint dynamics $(S^{L+\alpha K}_t)_{t\geq0}$
is the measure $\frac{1}{2\pi}\int d\varphi T \mu_{\varphi }$.
\end{cor}
Finally, together with part (2) of Proposition~\ref{TheJointDynamics}
which shows that there is no relaxation of the pure measure $\mu
'_\varphi $ under
$(S^{L+\alpha K}_t)_{t\geq0}$, we arrive at the proof of Theorem~\ref
{MainTheorem}.

\subsection{Extensions}\label{sec1.3}
Theorem~\ref{Bijection} stays true also for models where for every
angle there are
more than one Gibbs measures. This could occur for potentials with
highly nonconvex shapes \cite{EnSh05}.
The well-definedness of the rotation semigroup is untouched and one has:
%
\begin{teo} The map $T\dvtx \operatorname{ex} \mathcal{G}_{\theta
}(\gamma ^\Phi)\mapsto
\operatorname{ex} \mathcal{G}_{\theta}(\gamma ')$
is an equivariant bijection for the $S^1$-actions on continuous and
discrete-spin Gibbs measures.
\end{teo}
The equivariance property says $S^{L+\alpha K}_t(T \mu)= TR_t\mu$
for all
$\alpha \geq0$, where $R_t\mu$ is the measure obtained by joint
rotation of
the realizations of the measure $\mu$
by an angle $t$. The conclusions of Theorems~\ref{Bijection}~and~\ref{MainTheorem}, parts (1)~and~(3) apply. Theorem~\ref{MainTheorem}, part~(2)
(the uniqueness of the invariant measure) does not
apply because Corollary~\ref{TheUniqueInvariantMeasureForJoint}
does not apply since the symmetrization over the angles will produce
more then one invariant measure.

%
\begin{figure}[b]\vspace*{-12pt}
\[
  \xymatrix{
      \operatorname{ex} \mathcal{G}_{\theta}(\gamma^\Phi)
      \ar[rrr]^{\mu\mapsto R_t\mu} \ar@/_0,5cm/[dd]_T &  &  &
      \operatorname{ex} \mathcal{G}_{\theta}(\gamma^\Phi) \ar[dd]^T  \\ \\
      \operatorname{ex} \mathcal{G}_{\theta}(\gamma')
      \ar[rrr]_{\mu'\mapsto S_t^{L+\alpha K}(\mu')}
      \ar@/_0,5cm/[uu]_{\mu'\mapsto\int\mu'(d\omega')\mu_G[\omega'](d\omega)}        &   & &
      \operatorname{ex} \mathcal{G}_{\theta}(\gamma')
  }
\]
\caption{Equivariance property of the bijective
discretization map $T$ for the deterministic rotation action
$(R_t)_{t\geq0}$ and the action of the IPS $(S_t^{L+\alpha K})_{t\geq
0}$.}\label{Diagram}
\end{figure}

The remainder of the paper contains the following:
in Section~\ref{2} we prove Theorem~\ref{Bijection} using the
variational principle. For this we need
to present generalities and facts on discretizations and recall criteria
on the preservation of Gibbsianness.
In Section~\ref{3} we consider the rotation dynamics and prove Theorem
\ref{TheRotationDynamics}.
In Section~\ref{4} we consider the Glauber dynamics and prove
Proposition~\ref{TheGlauberDynamics}.
In Section~\ref{5} we consider the joint dynamics and prove the main
Proposition~\ref{TheJointDynamics}.

\section{Discretizations}\label{sec2}\label{2}

In the present section we will give a self-contained presentation
of properties of the discretization map $T$ which maps
continuous-spin Gibbs measures to discrete-spin Gibbs measures.
We will already obtain in this section the ``vertical'' parts of the
commutating diagram of
Figure~\ref{Diagram}, that is, those parts not involving dynamics.
These generalities about local discretizations we are going to present
are easily explained in a setup
which is broader than that of \mbox{$S^1$-}valued spins on an integer lattice.

Take an underlying site space $G$, a local spin-space $S$ equipped with
a $\sigma $-algebra and the configuration space ${\Omega}=S^G$
carrying the
product-$\sigma $-algebra.
$S^1$ will often serve as an example for the local state space, but one
can also consider subsets
of the Euclidean space of any finite dimension or finite-dimensional
manifolds. We will refer to this space
as the continuous spin-space.
Consider a Gibbsian potential $\Phi=(\Phi_A)_{A{\subset}G, A\
\mathrm{finite}}$
which is absolutely summable.
Write for the Hamiltonian in the finite volume $\Lambda \subset G$,
$H_{\Lambda }=\sum_{A\cap\Lambda \neq\varnothing}\Phi_A$ and let
$\gamma ^\Phi=(\gamma ^\Phi_{\Lambda
})_{\Lambda \subset G,\Lambda \ \mathrm{finite}}$ be the associated
Gibbsian specification with a priori measure $\lambda $.
We denote by $\mathcal{G}(\gamma ^\Phi)$ the corresponding Gibbs
measures, defined
by the DLR equation
and by $\mathcal{G}_{\theta}(\gamma ^\Phi)$ the
translation-invariant Gibbs
measures. Together we call this the first-layer system.\looseness=1

Let $S=\bigcup_{s'=1}^qS_{s'}$ be a disjoint
decomposition of the local state space into sets of positive $\lambda
$-measure.
As in \cite{HaKu04,EnKuOp11,KuOp08} the map $T(s):=s'$ for $S_{s'}\ni s$
defines a deterministic
transformation on $S$, called the discretization map.
The space ${\Omega}':=\{1,\ldots, q\}^G$ will be referred to as the discrete
or coarse-grained configuration space. It is convenient to use a
notation which
identifies the label $s'\in\{1, \ldots, q\}$ with the measurable
subset of
$S$ described by it and write $1_{s'}(s)=1$ if and only if $T(s)=s'$.


\begin{lem}\label{le2.1} For each fixed discrete-spin variable $\omega '\in{\Omega
}'$ define
a family of kernels on the continuous spin-space by constraining the
continuous spins to $\omega '$ and
putting, for each finite $\Lambda {\subset}G$, and bounded measurable
observable~$\varphi $,
%
\begin{equation}
\label{woe} \gamma ^{\omega '}_{\Lambda }(\varphi |\omega
_{\Lambda ^c} ):=\frac{ \gamma ^\Phi_{\Lambda }(\varphi 1_{\omega '_{\Lambda
}}|\omega _{\Lambda ^c} )}{\gamma ^\Phi_{\Lambda }(1_{\omega
'_{\Lambda }}|\omega _{\Lambda ^c} )}.
\end{equation}

Then this family defines
a Gibbsian specification $\gamma ^{\omega '}$ on ${\Omega}^{\omega
'}=\bigtimes_{i\in G} S_{\omega
'_i}$ in the sense of \cite{Ge11,EnFeSo93}.
\end{lem}

It will be useful to sometimes indicate measurability of functions
w.r.t. sub-$\sigma $-algebras in the following way: we write $f(\omega
'_\Lambda )$
equivalently to $f(\omega ')$ if $f$ evaluates $\omega '$ only inside
the volume
$\Lambda $. For example, in the case of (\ref{woe}) we write $\gamma
^{\omega '_{\Lambda
}}_{\Lambda }$ for $\gamma ^{\omega '}_{\Lambda }$.
%

\begin{pf*}{Proof of Lemma~\ref{le2.1}}
We verify the defining properties of a specification which need to be
fulfilled to be a useful candidate
system of conditional probabilities of an infinite-volume measure.
To begin with, from the compatibility property of $\gamma ^\Phi$
follows the
compatibility property of $\gamma ^{\omega '}$ for each fixed
$\omega '$.
%
The quasilocality of $\gamma ^\Phi$ implies that of $\gamma
^{\omega '}$ for all $\omega '$.
Since $\gamma ^\Phi$ is proper it is easy to see that $\gamma
^{\omega '}$ is proper,
where properness means for all finite $\Lambda \subset G$ and
$A\subset{\Omega}$
measurable\vadjust{\goodbreak} and dependent only on sites in $\Lambda ^c$ we have
$\gamma ^\Phi_\Lambda
(A|\cdot)=1_A$.
Finally the property of nonnullness on the constrained first-layer
local spin-spaces
(uniform boundedness of local probabilities from below) follows from
the positive measure of the sets in the decomposition
and the absolute summability of $\Phi$.
\end{pf*}

We will need to choose the discretization fine enough such that there
is only one Gibbs measure
which is compatible with this specification, for any $\omega '$. One
way to
see that this is always possible
and implement this requirement is to use
Dobrushin's uniqueness theory. From general results of the theory,
further information about the unique Gibbs measure follows, and we will
make use of this later.
We define a uniform Dobrushin matrix $\widebar C=(\widebar C_{ij} )_{i,j\in G}$
which is associated
to the family of specifications\vadjust{\goodbreak} $\gamma ^{\omega '}$, indexed by
$\omega '$, by
letting their entries be
%
\begin{equation}
\label{powerpower10} \widebar C_{ij}:=\sup_{\omega '}\mathop{\sup_{\omega,\tilde\omega\dvtx }}_{\omega _{j^c}=\tilde\omega _{j^c},
T(\omega )=T(\tilde\omega )=\omega '
} \big\| \gamma ^{\omega '}\bigl( \cdot|_i |
\omega \bigr)-\gamma ^{\omega '}\bigl( \cdot|_i | \tilde\omega \bigr)\big\|_i,
\end{equation}
where $\Vert\cdot\Vert_i$ is the total variational distance at site $i$
between the marginal distributions at site $i$; for details, see \cite
{EnKuOp11,Ge11}. Notice that we used another supremum over the discrete
configurations and hence the corresponding Dobrushin constant
$\bar c:=\sup_i \sum_j \widebar C_{ij}$ is uniform in $\omega '$.

We will always suppose that the discretization is fine enough such that
$\bar c<1$.
(Later we will even suppose a slightly stronger exponential decay property
that will appear in Lemma~\ref{ExponentialRegularityRotationRates}.)

Then it follows from the theory of Dobrushin uniqueness (see Theorem
8.23 in \cite{Ge11})
that, for any fixed $\omega '$ the specification $\gamma ^{\omega '}$
has a unique Gibbs measure. Moreover,
for each finite or infinite $V\subset G$ there is a kernel from
coarse-grained configurations $\omega '$ (inside $V$) and
boundary conditions of first-layer configurations $\omega $ outside
$V$, namely
$\gamma ^{\omega '}_V(\cdot| \omega )$, which has
the infinite-volume compatibility property $\gamma ^{\omega
'}_V\gamma ^{\omega '}_W=\gamma ^{\omega '}_V$,
between all (and not only finite) subsets
of $G$.

For the unique first-layer Gibbs measure for given discretized variable
$\omega '$, we use the notation
%
\begin{equation}
\label{UniqueConditionalGibbsMeasure} \mu_G \bigl[\omega ' \bigr](d\omega ):=
\gamma ^{\omega '}_G(d\omega ).
\end{equation}
We note that $\mu[\cdot](d\omega )$ is a probability kernel from
${\Omega}'$ to ${\Omega}
$, since
it is also measurable as a function of the coarse-grained configuration.


We report the result of \cite{EnKuOp11} which gives a criterion for the
fineness of the discretization in our
main example,
the standard nearest-neighbor model (the planar rotor or $XY$-model),
with Hamiltonian as given in (\ref{MetricFamily}):
For $q\geq q(\beta )$ large enough such that
$ 2 d \beta (\sin\frac{\pi}{q})^2 < 1$ we have $\bar c<1$. Notice
similar criteria are immediate for high-dimensional rotators; for
details, see \cite{EnKuOp11}.\vadjust{\goodbreak}

There is no obstacle to use this theory also for even more general models
to which the hypothesis of Theorem~\ref{MainTheorem} apply.
We report the bound on the matrix elements of the Dobrushin matrix
as given in \cite{EnKuOp11} which takes the form
%
\begin{equation}
\label{DiamS} \widebar C_{ij}\leq\sup_{s'}
\operatorname{diam}_{ i j} S_{s'}/4
\end{equation}
with a family of
metrics $(d_{ij})_{j\in G\setminus i}$ on
the local spin-space at site $i\in G$
which are generated by variations of the energy as follows:
\[
\label{lippi} d_{ij}(\sigma _i,
\tau_i):=\mathop{\sup_{\zeta,\bar\zeta\dvtx }}_{\zeta
_{j^c}=\bar\zeta_{j^c}; T(\zeta_{j})=T(\bar\zeta_{j})} \bigl| H_i(\sigma
_i\zeta _{i^c})-H_i(\sigma _i\bar\zeta _{i^c})- \bigl(H_i(\tau_i\zeta
_{i^c})- H_i(\tau_i\bar\zeta _{i^c})
\bigr) \bigr| %
\]
%
and $\operatorname{diam}_{ij}(S_{s'}):=\sup_{s,\tilde s\in
S_{s'}}d_{ij}(s,\tilde s)$.

Using the above criterion we suppose from now on that potential and
discretization
are chosen such that we are conditionally uniformly in the Dobrushin
regime $\bar c<1$.
We note that to each quasilocal continuous-spin observable $f$
there is naturally associated a discrete-spin observable $f'(\omega
'):=\mu
_G[\omega '](f)$
which is easily seen to be quasilocal as well (but on ${\Omega}'$) using
Dobrushin uniqueness \mbox{techniques}.
Denoting by $\mathcal{F}'$ the $\sigma $-algebra over ${\Omega}$
generated by the
infinite-volume coarse-graining map $T$, we have
that $f'$ is a regular version of the conditional expectation $\mu
(f|\mathcal{F}
')(\omega ')$ for every Gibbs measure
$\mu\in\mathcal{G}(\gamma ^\Phi)$, independently of its choice.

\begin{lem}
For a continuous-spin Gibbs measure $\mu$ denote its discretization
image by $\mu'=T\mu$.
Then the measures $\mu$ and $\mu'$ are close in the sense that $\mu
(f)=\mu'(f')$ for all continuous-spin observables
$f$, and moreover differences between corresponding correlations obey
the estimate
%
\begin{eqnarray}
\label{powerpoweri}
&& \bigl| \bigl( \mu(f g)- \mu(f)\mu(g) \bigr) - \bigl(
\mu' \bigl(f'g' \bigr)-
\mu' \bigl(f' \bigr) \mu'
\bigl(g' \bigr) \bigr)\bigr|
\nonumber\\[-8pt]\\[-8pt]
&&\qquad \leq\frac{1}{4}\sum
_{i,j\in G} \delta _i(f)\delta _j(g)\widebar
D_{ij}\nonumber
\end{eqnarray}
with the matrix $(\widebar D_{ij})_{i,j\in G}:=\sum_{n\geq0} \widebar C^n$ and
$g' =\mu(g|\mathcal{F}')$.
\end{lem}

\begin{pf}
To see that (\ref{powerpoweri}) holds, write
%
\begin{eqnarray*}\label{powerpowera}
&& \mu(f g)- \mu(f)\mu(g)
\\
&&\qquad =\mu \bigl(\mu \bigl(f g|
\mathcal{F}' \bigr) \bigr) - \mu(f)\mu(g)
\\
&&\qquad =
\mu' \bigl(\mu \bigl(f g|\mathcal{F}' \bigr) -\mu
\bigl(f |\mathcal{F}' \bigr) \mu \bigl(g |\mathcal{F}'
\bigr) \bigr) + \mu' \bigl(\mu \bigl(f |\mathcal{F}'
\bigr) \mu \bigl(g |\mathcal{F}' \bigr) \bigr)
\\
&&\quad\qquad{} - \mu'
\bigl(\mu \bigl(f |\mathcal{F}' \bigr) \bigr) \mu'
\bigl( \mu \bigl(g|\mathcal{F} ' \bigr) \bigr).
\end{eqnarray*}
%
Further the standard estimate (see Proposition 8.34 in \cite{Ge11}) in
the Dobrushin uniqueness regime yields
%
\begin{equation}
\label{powerpower14} \quad\sup_{\omega '}\bigl|\mu_G \bigl[\omega
' \bigr](fg) -\mu_G \bigl[\omega '
\bigr](f) \mu _G \bigl[\omega ' \bigr](g)\bigr| \leq
\frac
{1}{4} \sum_{i,j\in G} \delta
_i(f)\delta _j(g)\widebar D_{ij},
\end{equation}
which proves (\ref{powerpoweri}).
\end{pf}

On the lattice this statement can be used to see that power law decay
of correlations for a continuous-spin
observable $f$ (as it can occur in the standard rotor model in space
dimension $2$)
carries over to power law decay between correlations in the associated
observable $f'$ when the discretization
is fine enough, since in that case the matrix elements of $D$ are
decaying exponentially fast.

%
%


It is clear that the map from $\mu$ to $\mu':=T\mu$ is injective when
viewed on the (not necessarily translation-invariant)
Gibbs measures
of the continuous-spin system:
indeed, we can restore an initial Gibbs measure $\mu$ from
its coarse-grained image via $\mu(\varphi )=\int\mu'(d\omega ')\mu
_G[\omega '](\varphi )$
where $\mu_G[\omega '](\varphi )$ does not depend on $\mu$. Hence
different $\mu$'s
must have different images $\mu'$.

Next recall the definition of the specification $\gamma '$ for the
coarse-grained system (see also \cite{KuOp08}) given in (\ref
{CoarseSpecification}), that we will sometimes also call the
second-layer system.
We have the following lemma.

\begin{lem} In the uniform Dobrushin regime, the discretization image
of any continuous-spin Gibbs measure is Gibbs for the specification
$\gamma '$.
\end{lem}

\begin{pf}
This is shown by standard arguments which we include
for convenience of the reader.
Any conditional probability with finite-volume conditioning can be
written as
%
\begin{eqnarray}
\label{powerpower15}
&& \mu' \bigl(\omega '_{\Lambda '} |
\omega ' _{\Lambda \setminus
\Lambda '} \bigr)\nonumber
\\
&&\qquad =\frac{\int\mu(d\omega _{\Lambda ^c})\gamma
_{\Lambda }(1_{\omega '_{\Lambda '}}1_{\omega '_{\Lambda
\setminus\Lambda '}}|\omega _{\Lambda ^c})}{\int\mu(d\omega
_{\Lambda ^c})\gamma _{\Lambda }(1_{\omega '_{\Lambda \setminus
\Lambda '}}|\omega _{\Lambda ^c})}
\\
&&\qquad =\frac{\int\mu(d\omega _{\Lambda ^c})(\gamma ^{\omega
'_{\Lambda \setminus\Lambda '}}|_{G\setminus\Lambda
'})_{\Lambda \setminus\Lambda '}(\lambda ^{\Lambda '}
(e^{-H_{\Lambda '}}
1_{\omega '_{\Lambda '}}) | \omega _{\Lambda ^c})
}{\int\mu(d\omega _{\Lambda ^c})(\gamma ^{\omega '_{\Lambda
\setminus\Lambda '}}|_{G\setminus\Lambda
'})_{\Lambda \setminus\Lambda '}(\lambda ^{\Lambda '}
(e^{-H_{\Lambda '}}) | \omega _{\Lambda ^c})},\nonumber
\end{eqnarray}
where $\mu'(\omega '_{\Lambda '})=\mu(1_{\omega '_{\Lambda '}})$
and $\gamma ^{\omega '}|_{G\setminus
\Lambda '}$ denote the specification on
${\Omega}^{\omega '_{G\setminus\Lambda '}}_{G\setminus\Lambda
'}=\bigtimes_{i\in G\setminus
\Lambda '} S_{\omega '_i}$
obtained by putting all potentials $\Phi_A$ with $A\cap\Lambda '\neq
\varnothing
$ equal to zero.

Then, by martingale convergence, $\mu' (\omega '_{\Lambda '} |
\omega ' _{\Lambda \setminus
\Lambda '})$ converges as $\Lambda $ tends to $G$ in the a.s.- and $L^1$-sense
to $\mu(1_{\omega '_{\Lambda '}} | \mathcal{F}' _{G\setminus
\Lambda '})(\omega _{G\setminus\Lambda
'})$ where $1_{\omega '_{\Lambda \setminus\Lambda '}}(\omega
_{G\setminus\Lambda '})=1$ for all
$\Lambda \supset\Lambda '$ and $\mathcal{F}' _{G\setminus\Lambda
'}$ is the $\sigma $-algebra over
${\Omega}$ generated by the coarse-graining map $T$ applied only in the
infinite-volume $G\setminus\Lambda '$.

%
%
%

On the other hand, for any finite $\Lambda '$, there is convergence
uniformly in the integration variable $\omega $ under the $\mu$-integrals
since the conditional
specification is in the uniform Dobrushin regime, and we have
%
\begin{eqnarray}
\label{CoarseSpecification1}
&& \gamma '_{\Lambda '} \bigl(\omega
'_{\Lambda '} | \omega ' _{G\setminus\Lambda '} \bigr)\nonumber
\\
&&\qquad = \frac{\lim_{\Lambda \uparrow G }(\gamma ^{\omega '_{G\setminus
\Lambda '}}|_{G\setminus
\Lambda '})_{\Lambda \setminus\Lambda '}(\lambda ^{\Lambda '}
(e^{-H_{\Lambda '}}
1_{\omega '_{\Lambda '}}) | \omega _{G\setminus\Lambda })
}{\lim_{\Lambda \uparrow G }(\gamma ^{\omega '_{G\setminus
\Lambda '}}|_{G\setminus\Lambda
'})_{\Lambda \setminus\Lambda '}(\lambda ^{\Lambda '}
(e^{-H_{\Lambda '}}) | \omega _{G\setminus\Lambda })}
\\
&&\qquad =\frac{\mu_{G\setminus\Lambda '}[\omega '_{G\setminus\Lambda
'}](\lambda ^{\Lambda '} (e^{-H_{\Lambda '}}
1_{\omega '_{\Lambda '}}) )
}{\mu_{G\setminus\Lambda '}[\omega '_{G\setminus\Lambda
'}](\lambda ^{\Lambda '} (e^{-H_{\Lambda
'}}))}.\nonumber
\end{eqnarray}
The limiting measure in the last line is the unique Gibbs measure
of the specification restricted to $G\setminus \Lambda '$ with open boundary
conditions, and this proves (\ref{CoarseSpecification}).
\end{pf}

It is easy to see using the standard Dobrushin estimates that the
specification $\gamma '$ built with these kernels is quasilocal.

Now we are in the position to discuss new results
which are related to the proof of the bijectivity of the map $T$.
To start, note that we also have that the influence of variations of
the boundary condition
outside $\Lambda '$ on probabilities inside $\Lambda '$ has the estimate,
uniformly in the configuration $ \omega '_{\Lambda '}$,
%
\begin{eqnarray}
\label{boundaryforgammaprime} 
\log\frac{\gamma '_{\Lambda '} (\omega '_{\Lambda '}| \omega '
_{G\setminus\Lambda '})} {
\gamma '_{\Lambda '} (\omega '_{\Lambda '} | \bar\omega '
_{G\setminus\Lambda '})} 
&\leq&4\sum
_{A\cap\Lambda '\neq\varnothing,A\cap\Lambda '^c\neq
\varnothing}\Vert\Phi _A\Vert. 
\end{eqnarray}
%

Further note that for summable potentials and $\Lambda '$ being cubes
on the
lattice, the RHS is bounded by a constant times the length
of the boundary of $\Lambda '$, in other words $\log\frac{d\gamma
'_{\Lambda '} (\cdot|
\omega ' _{G\setminus\Lambda '})}{d\gamma '_{\Lambda '}(\cdot|
\bar\omega ' _{G\setminus\Lambda
'})}=O(|\partial\Lambda '|)$, where $|\cdot|$ denotes the cardinality.


Let us now restrict to the lattice case, that is, $G={\mathbb Z}^d$ 
and discuss the relative
entropy density.
The following lemma should be seen as a generalization of
the contractivity of the relative entropy (density) between two
measures (see Lemma~3.3
in \cite{EnFeSo93})
under strictly local transforms to transforms which are not strictly
but ``sufficiently'' local.

\begin{lem}\label{VariantionalPrinciple} Let $\mu'_1,\mu'_2 \in
\mathcal{G}(\gamma
')$ for some specification
for which\break  $\log\frac{d\gamma _{\Lambda }'(\cdot|_{\Lambda }
|\omega _1')}{d\gamma _{\Lambda }'(\cdot|_{\Lambda
} |\omega _2')}$ is\vspace*{1pt} of the order $o(|\Lambda |)$ for cubes. Take a
kernel $\mu_G[\omega
'](d\omega )$ where $\log\frac{d\mu_G[\omega '](\cdot|_{\Lambda
})}{d \mu_G[\bar\omega
'](\cdot|_{\Lambda })}$ is also of the order $o(|\Lambda |)$
uniformly in all
configurations $\omega '$~and~$\bar\omega '$ which coincide on
$\Lambda $.
Then the relative entropy density
between the mapped measures equals zero, that is,
%
\begin{equation}
\label{powerpower18} \lim_{\Lambda \uparrow{\mathbb Z}^d}\frac{1}{|\Lambda |}H \biggl(\int
\mu'_1 \bigl(d\tilde\omega ' \bigr) \mu
_G \bigl[\tilde\omega ' \bigr]|_{\Lambda } \Big| \int
\mu'_2 \bigl(d\tilde\omega ' \bigr)
\mu_G \bigl[\tilde\omega ' \bigr]|_{\Lambda
}
\biggr) =0
\end{equation}
along cubes.
\end{lem}

\begin{pf}
We need to estimate the relative entropy $H$ in a volume
$\Lambda $ where $\Lambda \subset{\mathbb Z}^d$ is a finite cube
appearing in the formula above, which is
%
\begin{equation}
\label{powerpowerEntropy} \int\mu'_1 \bigl(d\tilde\omega
' \bigr) \mu_G \bigl[\tilde\omega '
\bigr]|_{\Lambda
} \biggl(\log\frac{d\int\mu
'_1(d\tilde\omega ') \mu_G[\tilde\omega ']|_{\Lambda }}{d\int\mu
'_2(d\tilde\omega ') \mu
_G[\tilde\omega ']|_{\Lambda }} \biggr).
\end{equation}
Using the DLR equation for the integrand as well as the conditions on
the Radon--Nikodym derivatives, we find
%
\begin{eqnarray}
\label{powerpowerEntropy2}
&& \log\frac{d\int\mu'_1(d\tilde\omega ')\mu_G[\tilde\omega
']|_{\Lambda }}{d\int\mu
'_2(d\tilde\omega ')\mu_G[\tilde\omega ']|_{\Lambda }}\nonumber
\\
&&\qquad =\log\frac{\int\mu'_1(d\tilde\omega '_1) ((d\mu_G[\tilde\omega '_1]|_\Lambda )/(d\lambda ^\Lambda ))
}{\int\mu'_2(d\tilde\omega '_2) ((d\mu_G[\tilde\omega
'_2]|_\Lambda )/(d\lambda ^\Lambda ))}
\nonumber\\[-8pt]\\[-8pt]
&&\qquad \leq \sup
_{\omega '_1,\omega '_2}\log\frac{\int(\gamma
'_{\Lambda })|_{\Lambda }(d\tilde\omega '_1|\omega
'_1) ((d\mu_G[(\tilde\omega '_1)_\Lambda (\omega
'_1)_{\Lambda ^c}]|_\Lambda )/(d\lambda ^\Lambda ))
}{\int(\gamma '_{\Lambda })|_{\Lambda }(d\tilde\omega '_2|\omega
'_2)((d\mu_G[(\tilde\omega '_2)_\Lambda
(\omega '_2)_{\Lambda ^c}]|_\Lambda )/(d\lambda ^\Lambda ))}\nonumber
\\
&&\qquad =o\bigl(|\Lambda |\bigr),\nonumber
\end{eqnarray}
where the estimate in the last line uses the two assumptions in the hypothesis.
Hence the relative entropy density as the limit of the relative entropy
devided by the volumes of a cofinal sequence of cubes is equal to zero.
\end{pf}


Applying the lemma and using now the Gibbs variational principle in the
form of
Theorem 15.37 of \cite{Ge11}, our desired result, stating that every
discrete Gibbs measure has a continuous preimage, follows:

\begin{prop} Let $\mu' \in\mathcal{G}_{\theta}(\gamma ')$, then
$\mu(d\omega ):=\int\mu
'(d\omega ') \mu_G[\omega '](d\omega )\in\mathcal{G}_{\theta
}(\gamma ^\Phi)$.
\end{prop}

\begin{pf}
Let $\mu_0\in\mathcal{G}_{\theta}(\gamma ^\Phi)$ be a Gibbs
measure for the
original system and $\mu_0':=T\mu_0$ its coarse-grained image. We want
to use the preceding lemma, that is, to justify the conditions and
therefore conclude that the relative entropy density 
between the two translation-invariant measures is zero. Hence, by the
variational principle applied to the original system, also $\mu\in
\mathcal{G}
_{\theta}(\gamma ^\Phi)$.

%
%
%
%
Indeed, (\ref{boundaryforgammaprime}) asserts the condition of Lemma
\ref{VariantionalPrinciple} for the coarse-grained specification
$\gamma
'$. Also we have for $\omega ',\bar\omega '$ coinciding on $\Lambda $
\begin{eqnarray*}
\log\frac{d\mu_G[\omega ']|_{\Lambda }}{
d\mu_G[\bar\omega ']|_{\Lambda }} &=&\log\frac{\int\mu_G[\omega '](d\tilde\omega _1)
((d(\gamma ^{\omega '}_{\Lambda })|_{\Lambda })/(d
\lambda ^{\Lambda }))(\cdot|\tilde\omega _1)}{
\int\mu_G[\bar\omega '](d\tilde\omega _2)((d(\gamma
^{\omega '}_{\Lambda })|_{\Lambda })/(d \lambda ^{\Lambda
}))(\cdot|\tilde\omega _2)}
\\
&\leq&\sup
_{\tilde\omega _1,\tilde\omega _2}\log\frac{d\gamma
^{\omega '}_{\Lambda }(\cdot
|\tilde\omega _1)}{d \gamma ^{\omega '}_{\Lambda } (\cdot|\tilde\omega _2)} \leq4\sum
_{A\cap\Lambda '\neq\varnothing,A\cap\Lambda '^c\neq
\varnothing}\Vert\Phi _A\Vert=o\bigl(|\Lambda |\bigr).
\end{eqnarray*}\upqed
\end{pf}

Together with the injectivity of $T$ this means that the map from the
translation-invariant Gibbs measures of the original system $\mathcal
{G}_{\theta
}(\gamma ^\Phi)$ to the translation-invariant measures for the
coarse-grained configuration $\mathcal{G}_{\theta}(\gamma ')$ is one-to-one.

\begin{rem}
This one-to-one correspondence also holds for the extremals:
if $\mu$ is tail-trivial, then so is $T\mu$ since the tail-$\sigma $-algebra
of discrete
events is contained in the tail-$\sigma $-algebra of all events,
$\mathcal{T}'
\subset\mathcal{T}$.
In particular $\operatorname{ex}\mathcal{G}(\gamma ')\supset
T(\operatorname{ex}\mathcal{G}
(\gamma ^\Phi))$. To
see that also $\mu\in\operatorname{ex}\mathcal{G}(\gamma ^\Phi)$
for $T\mu\in
\operatorname{ex}\mathcal{G}
(\gamma ')$ one can use the fact that the mapping $T$ is affine: let us
assume $T\mu\in\operatorname{ex}\mathcal{G}(\gamma ')$ and $\mu
=s\mu_1+(1-s)\mu_2$ for
$s\in[0,1]$ and $\mu_1,\mu_2\in\mathcal{G}(\gamma ^\Phi)$. Then
we have $T\mu=sT\mu
_1+(1-s)T\mu_2$ and hence $T\mu=T\mu_1=T\mu_2$
since $T\mu$ is extremal. But that means $\mu=\mu_1=\mu_2$ and thus
$\mu
\in\operatorname{ex}\mathcal{G}(\gamma ^\Phi)$.
\end{rem}

It is interesting to note that the proof of the preceding remark also
follows from the fact that tail-triviality is preserved under the kernel
(even not assuming initial Gibbs measures). This property
explains the ``essentially local'' nature of the transformation $T$ from
the perspective of the tail events.

\begin{prop} Assume that $\mu'$ is any probability measure (not
necessarily Gibbs) on ${\Omega}'$ which is
trivial on $\mathcal{T}'$. Then
$\mu(d\omega ):=\int\mu'(d\omega ') \mu_G[\omega '](d\omega )$
is trivial on the tail-$\sigma
$-algebra $\mathcal{T}$.
\end{prop}

\begin{pf}
 We assume that also $\sup_j\sum_i \widebar C_{ij}<1$ which
is guaranteed in
the fine-discretization regime ensured by our criteria.

If $A\in\mathcal{T}$ then $\mu_G[\omega '](A)$ is $\mathcal
{T}'$-measurable.
To see this, suppose that $W$ is a finite subset of $G$, that $V$
contains $W$ and
that $A$ is in $\mathcal{T}_{V}$, the $\sigma $-algebra of events not
depending on
spins inside $V$.
Assuming that $A$ is a cylinder, at first we have
%
\begin{eqnarray}
\label{powerova21}
\sup_{\omega ',\bar\omega '\dvtx \omega '_{W^c}=\bar\omega '_{W^c}} \bigl(\mu_G \bigl[
\omega ' \bigr](A)-\mu_G \bigl[\bar\omega
' \bigr](A) \bigr)&\leq&\sum_{i\in
\operatorname{supp}
(A), j\in
W}\widebar
D_{ij}
\nonumber\\[-8pt]\\[-8pt]
&\leq&\sum_{i\in V^c, j\in W}\widebar
D_{ij}.\nonumber
\end{eqnarray}
Next we note that this inequality also holds
by approximation of probabilities of general events by cylinders (by a
semiring-approximation argument)
for all \mbox{$A\in\mathcal{T}_{V}$}. Since $A\in\mathcal{T}$ is in any
$\mathcal{T}_{V}$ we may let
$V\nearrow G$
and obtain that
%
\begin{equation}
\label{powerova22} \sup_{\omega ',\bar\omega '\dvtx \omega '_{W^c}=\bar\omega '_{W^c}} \bigl(\mu_G \bigl[
\omega ' \bigr](A)-\mu_G \bigl[\bar\omega
' \bigr](A) \bigr)\leq0.
\end{equation}
Since $W$ was arbitrary, this is the tail-measurability.

Further we note that $\mu_G[\omega '](A)\in\{0,1\}$ for each fixed
$\omega '$
and $A\in\mathcal{T}$
since the original measure constrained to coarse-grained configurations
is in the Dobrushin uniqueness regime, hence tail-trivial.
So $\mu_G[\omega '](A)=1_{A'}(\omega ')$ for some $A'\in\mathcal
{T}'$ and this
implies $\mu(A)=\int\mu'(d\omega ')\mu_G[\omega '](A)=\mu
'(A')\in\{0,1\}$
by tail-triviality of $\mu'$.
\end{pf}

\section{Continuous rotations for discrete-spin models}\label{sec3}\label{3}

After the preparations of the last section we turn now to the
discussion of the rotation dynamics.
Let us specialize to a translation-invariant $S^1$-model and look at
the Markov process given in (\ref{rotationgenerator}) with rates given
in (\ref{rotationjumprates}).
%
%
%
%

Intuitively the choice of the rates
can be understood as follows: consider the single-site discrete
observable $f(\sigma ')=1_a(\sigma '_i)$ with
$a\in\{1,\ldots, q\}$ fixed, and let $\mu_{\varphi }$ be an extremal
translation-invariant Gibbs measure of the $d\geq3$ $XY$-model labeled by
the angle $\varphi $. Then we have
%
\begin{eqnarray}
\label{Intuitive}
&& \frac{d}{dt}_{|_{t=0}}(T\mu_{\varphi +t})
\bigl(1_a \bigl(\sigma _i' \bigr) \bigr)\nonumber
\\
&&\qquad =
\frac{d}{dt}_{|_{t=0}}\mu_{\varphi } \bigl(\sigma _i
\in \bigl(a|^l-t, a|^r -t \bigr) \bigr)\hspace*{-10pt}\nonumber
\\
&&\qquad =\frac{d}{dt}_{|_{t=0}}\int(T\mu_{\varphi }) \bigl(d \omega
' \bigr)\mu _{G} \bigl[\omega ' \bigr]
\bigl( \sigma _i \in \bigl(a|^l-t, a|^r -t
\bigr) \bigr)\hspace*{-10pt}\nonumber
\\
&&\qquad =\frac{d}{dt}_{|_{t=0}}\int(T\mu_{\varphi }) \bigl(d \omega
' \bigr)\int_{a|^l -t }^{a|^r
- t}\frac{d\mu_{G}[\omega ']{ |_i}}{d\lambda }(s)\,ds
\nonumber\\[-8pt]\\[-8pt]
&&\qquad =\frac{d}{dt}_{|_{t=0}}\int(T\mu_{\varphi }) \bigl(d \omega
' \bigr)\int_{a|^l -t }^{a|^r
- t}\frac{\mu_{G\setminus i}[\omega '_{G\setminus
i}](e^{-H_{i}(s,\cdot
_{i^c})}) }{\mu_{G\setminus i}[\omega '_{G\setminus i}](\lambda ^i
( e^{-H_{i}}
1_{a}) )}\,ds\hspace*{-10pt}\nonumber
\\
&&\qquad =\int(T\mu_{\varphi }) \bigl(d \omega ' \bigr) \biggl(
\frac{\mu_{G\setminus i}[\omega '_{G\setminus
i}](e^{-H_{i}(a|^l,\cdot
_{i^c})}) }{\mu_{G\setminus i}[\omega '_{G\setminus i}](\lambda ^i
( e^{-H_{i}}
1_{a}) )}-\frac{\mu_{G\setminus i}[\omega '_{G\setminus
i}](e^{-H_{i}(a|^r,\cdot_{i^c})}) }{\mu_{G\setminus i}[\omega
'_{G\setminus i}](\lambda ^i ( e^{-H_{i}}
1_{a}) )} \biggr)\hspace*{-10pt}\nonumber
\\
&&\qquad =\int(T\mu_{\varphi }) \bigl(d \omega
' \bigr) c_L \bigl(\omega ', \bigl(
\omega ' \bigr)^i \bigr) \bigl(1_a \bigl(
\bigl( \omega ' \bigr)^i \bigr)- 1_a
\bigl(\omega ' \bigr) \bigr)\nonumber\hspace*{-10pt}
\\
&&\qquad =(T\mu_\varphi ) (Lf),\nonumber\hspace*{-10pt}
\end{eqnarray}
where in the second line we wrote $a|^l:=2\pi(a-1)/q$ (resp.,
$a|^r:=2\pi
a/q$) to indicate the left (resp., right) endpoint of $a$. In the third
line we used Theorem~\ref{Bijection}. In the fourth line we rewrote the
constrained Gibbs measure as a marginal density (at site $i$) w.r.t. the
Lebesgue measure $\lambda $ on the sphere, which as indicated in the fifth
line can again be re-expressed by seperating the part of the potential
interacting with the site $i$.

\subsection{Well-definedness of the rotation generator}\label{sec3.1}

In this subsection we prove Theorem~\ref{TheRotationDynamics}, part (1).
We use methodology of Liggett \cite{Li85} via the Hille--Yosida theorem
to prove well-definedness. Let us start with an overview on function
spaces we need in the investigation of the dynamics.
%
%
\begin{defn}\label{Spaces}
Let us fix the following notation. We write:
\begin{longlist}[(3)]
\item[(1)] $\mathcal{L}':=\{f\dvtx  {\Omega}'\to{\mathbb R}\dvtx  f\mbox{ is
local}\}$ for the
local functions.\vadjust{\goodbreak}
\item[(2)] $C({\Omega}')
=\mathcal{L}'_{\|\cdot\|}$ equivalently for the space of continuous
functions on the compact configuration-space ${\Omega}'$ which, since
$q$ is
finite, coincides with the space of bounded quasilocal functions which
is just the $\|\cdot\|$-completion of the local functions. Here
$\|\cdot\|$ denotes the uniform norm.
\item[(3)] $D({\Omega}'):=\{f\in C({\Omega}')\dvtx \|\!| f\|\!|:=\sum_{i\in G}\delta
_i(f)<\infty\}$ for the core functions.
\item[(4)] $\mathcal{L}'_{\|\!|\cdot\|\!|}$ for the triple-norm
completion of the local functions.
\item[(5)] $D_{p(\varrho )}({\Omega}'):=\{g\in C({\Omega}')\dvtx  \|\!|
g\|\!|_{p(\varrho )}:=\sum_{i\in G}p(\varrho (i,0))\delta_{i}(g)<\infty\}$ for the space of weighted
triple-normed functions, where $\varrho $ is an increasing,
translation-invariant semi-metric on the site space and $p\dvtx  {\mathbb
R}_0^+\to{\mathbb R}
_0^+$ any weight-function.
\end{longlist}
\end{defn}
%
Let us clarify the relations between those spaces and specialize to $p$
being either an exponential function with some factor $\varepsilon >0$
or a
monomial function with power $m\in{\mathbb N}$. Let the semi-metric just
be the
Euclidean metric $|\cdot|$ on an ordering of $G$. We have
%
\begin{equation}
\label{Space-relations} \mathcal{L}'\subset D_{e^{\varepsilon |\cdot|}} \bigl({
\Omega}' \bigr)\subset D_{|\cdot|^{m}} \bigl({\Omega}
' \bigr)\subset D_{|\cdot|^{1}} \bigl({\Omega}'
\bigr) \subset\mathcal{L}'_{\|\!|\cdot
\|\!|
}\subset D \bigl({
\Omega}' \bigr)\subset C \bigl({\Omega}' \bigr).
\end{equation}
Notice that all of these spaces are dense in $C({\Omega}')$ with
respect to
the $\|\cdot\|$-norm.
All inclusions should be clear except $D_{|\cdot|^{1}}({\Omega
}')\subset
\mathcal{L}'_{\|\!|\cdot\|\!|}$.

\begin{prop}\label{triplenormcriteria}
$D_{|\cdot|^{1}}({\Omega}')\subset\mathcal{L}'_{\|\!|\cdot\|\!|}$.
\end{prop}
\begin{pf}Let $f\in D_{|\cdot|^{1}}({\Omega}')$ for an
ordering $o\dvtx G\to
{\mathbb N}$. Define $\Lambda _i:=\break \{j\in G\dvtx  o(j)\leq o(i)\}$ an
exhausting sequence of
finite volumes, then\break  $\sum_{i\in G}|\Lambda _i|\delta_{i}(f)=\sum_{i\geq
0}i\delta_{o^{-1}(i)}(f)<\infty$ and $\sum_{i\geq n}i\delta
_{o^{-1}(i)}(f)\to0$ for $n\to\infty$. Let $\eta\in{\Omega}'$ be
fixed, and
define a sequence of local functions $f_n(\omega ):=f(\omega
_{\Lambda _n}\eta_{\Lambda
_n^c})$, and then we have
%
\begin{eqnarray}
\label{triplenormcriterium2} \|\!| f-f_n\|\!|&=&\sum_{i\in\Lambda _n}
\delta_i(f-f_n)+\sum_{i\in
\Lambda
_n^c}
\delta_{i}(f)\leq2n\|f-f_n\|+\sum
_{i\in\Lambda
_n^c}\delta_{i}(f)\nonumber
\\
&\leq&2n\sum_{i\in\Lambda _n^c}\delta_i(f) +\sum
_{i\in\Lambda
_n^c}\delta_i(f)\leq 2\sum
_{i>n}i\delta_{o^{-1}(i)}(f)\to0
\\
\eqntext{\mbox{for }n\to
\infty.}
\end{eqnarray}
Hence $f\in\mathcal{L}'_{\|\!|\cdot\|\!|}$ and $D_{|\cdot
|^{1}}({\Omega}
')\subset\mathcal{L}'_{\|\!|\cdot\|\!|}$. \end{pf}

In the sequel, we will drop the notation $o^{-1}(i)$ and just write
$\sum_{i\geq0}i\delta_{i}(f)$.

Let us check the criteria for well-definedness proposed in \cite{Li85}.
Note the jump rates are uniformly bounded
since we assumed the potential to be absolutely summable and
translation-invariant and the coarse-graining to be finite. Further the
rates have to be of bounded variation, that is:
%
\begin{lem}\label{boundedvariation} $\sup_{i\in G}\sum_{j\neq
i}\delta
_j(c_L(\cdot,\cdot^i))<\infty$ if $\sup_{i\in G}\sum_{A\ni i}\|\!
|\Phi
_A\|\!|<\infty$.
\end{lem}
\begin{pf}This follows from the Dobrushin comparison theorem (see
\cite{Ge11}, Theorem 8.20), indeed,
\begin{eqnarray*}
\label{welldefined2}
&& \delta_j \bigl(c_L \bigl(\cdot,
\cdot^i \bigr) \bigr) 
\\
&&\qquad \leq Ce^{\|H_{i}\|}\mathop{\sup_{\omega '=\tilde\omega '/}}_{
{\mathrm{off}\ j}}\bigl|\mu _{G\setminus i} \bigl[\omega
'_{G\setminus i} \bigr] \bigl(e^{-H_{i}(\omega
'_i|^r,\cdot
_{i^c})} \bigr) -
\mu_{G\setminus i} \bigl[\tilde\omega '_{G\setminus
i} \bigr]
\bigl(e^{-H_{i}(\tilde\omega '_i|^r,\cdot_{i^c})} \bigr) \bigr|
\\
&&\quad\qquad{} +Ce^{3\|H_{i}\|}\mathop{\sup_{\omega '=\tilde\omega ' }}_{\mathrm
{off\ }j}\bigl|\mu _{G\setminus i} \bigl[\omega
'_{G\setminus i} \bigr] \bigl(\lambda ^i \bigl(
e^{-H_{i}} 1_{\omega '_{i}} \bigr) \bigr) -\mu_{G\setminus i} \bigl[\tilde\omega '_{G\setminus i} \bigr] \bigl(\lambda ^i \bigl(
e^{-H_{i}} 1_{\tilde\omega '_{i}} \bigr) \bigr)\bigr|
\end{eqnarray*}
and therefore it suffices to look at the Gibbs measures $\mu
_{G\setminus i}[\omega '_{G\setminus i}]$ and\break  $\mu_{G\setminus
i}[\tilde\omega
'_{G\setminus i}]$ on
$(S^1)^{G\setminus i}$
applied to the quasilocal functions $\psi_1^{\omega '_i}(\cdot
):=e^{-H_{i}(\omega
'_i|^r,\cdot_{i^c})}$ and $\psi_2^{\omega '_i}(\cdot):=\lambda ^i (
e^{-H_{i}(\cdot,\cdot_{i^c})}1_{\omega '_{i}})$. For any fixed first-layer
boundary condition $\omega \in{\Omega}$, the measure $\mu
_{G\setminus i}[\omega
'_{G\setminus i}](\cdot)$ is uniquely specified by the specification
%
\begin{equation}
\label{welldefined3} \gamma ^{\omega '_{G\setminus i}}:= \bigl( \bigl(\gamma ^{\omega
'_{G\setminus i}}|_{G\setminus
i}
\bigr)_{\Lambda \setminus i}(\cdot|\omega _{\Lambda ^c\setminus
i}) \bigr)_{\Lambda \subset
G\setminus i },
\end{equation}
$\Lambda $ being finite subsets of $G\setminus i$. We have for
$\omega '_{G\setminus
j}=\tilde\omega '_{G\setminus j}$
%
\begin{equation}
\label{welldefined4} \big\| \bigl(\gamma ^{\omega '_{G\setminus i}}|_{G\setminus
i}
\bigr)_{l\setminus i}(\cdot|\omega _{l^c\setminus i})- \bigl(\gamma
^{\tilde\omega '_{G\setminus
i}}|_{G\setminus
i} \bigr)_{l\setminus i}(\cdot|\omega
_{l^c\setminus i})\big\|_l\leq1_{l=j}.
\end{equation}
Hence for $\omega '_{G\setminus j}=\tilde\omega '_{G\setminus j}$
and $\psi^{\omega
'_i}\in\{\psi_1^{\omega '_i},\psi_2^{\omega '_i}\}$ the comparison
theorem gives us
%
\begin{eqnarray*}
\label{welldefined5} \bigl|\mu_{G\setminus i} \bigl[\omega '_{G\setminus i}
\bigr] \bigl(\psi^{\omega
'_i} \bigr)-\mu _{G\setminus i} \bigl[\tilde\omega
'_{G\setminus i} \bigr] \bigl(\psi^{\omega
'_i} \bigr)\bigr|&\leq&\sum
_{k\neq i}\delta_k \bigl(
\psi^{\omega '_i} \bigr)D_{kj} \bigl(\gamma ^{\omega
'_{G\setminus i}} \bigr)
\\
&\leq&\sum_{k\neq i}\delta_k \bigl(
\psi^{\omega '_i} \bigr) \widebar D_{kj},
\end{eqnarray*}
where\vspace*{1pt} we used the fact that the specifications $\gamma ^{\omega
'_{G\setminus
i}}$ are in the Dobrushin region uniformly in the
constraint $\omega '$.
%
%
Since $\bar c:=\sup_i \sum_j \widebar C_{ij}<1$ we have \mbox{$\sum_{j\in
G}\widebar
D_{kj}<\infty$} for all $k\in G$ and can therefore conclude
\begin{eqnarray*}
\label{welldefined7} \sup_{i\in G}\sup_{\omega '_i\in S'}\sum
_{j\neq i}\sum_{k\neq
i}
\delta_k \bigl(\psi ^{\omega '_i} \bigr)\widebar D_{kj}&\leq& C\sum_{k\in G}\sup_{i\in G}\sup
_{\omega '_i\in
S'}\delta_k \bigl(\psi^{\omega '_i} \bigr)
\leq C\sup_{i\in G}\sum_{k\in
G}
\delta_k \bigl(\psi ^i \bigr)
\end{eqnarray*}
with $\psi^{i}\in\{\psi_1^{i},\psi_2^{i}\}$ and $\psi_1^{i}(\cdot
):=e^{-H_{i}(\cdot,\cdot_{i^c})}$ and $\psi_2^{i}(\cdot):=\lambda
^i (
e^{-H_{i}(\cdot,\cdot_{i^c})})$. In case $\psi^i$ is a local function,
uniformly bounded in $i$ (e.g., in the $XY$-model), the sum is
finite and thus less than infinity. In the general case were the $\psi
^i$ are coming from an uniformly bounded Hamiltonian which is only
quasilocal the summability is not guaranteed. But if we stipulate $\sup_{i\in G}\sum_{A\ni i}\|\!|\Phi_A\|\!|<\infty$, we have for $\psi
^i_1$ and~$\psi^i_2$
\begin{eqnarray*}
\label{welldefined8} \sum_{k\in G}\delta_k
\bigl(\psi^i_2 \bigr) 
&\leq& C
\sum_{k\in G}\delta_k \bigl(
\psi^i_1 \bigr) 
\leq Ce^{\Vert H_0\Vert} \sum
_{k\in G}\sum_{A\ni{ik}}
\delta_k(\Phi_A)= Ce^{\Vert H_0\Vert} \sum
_{A\ni{i}}\|\!|\Phi_A\|\!|<\infty,
\end{eqnarray*}
where we used $|e^x-e^y|\leq|x-y|e^{\max(|x|,|y|)}$.
\end{pf}

Note that in particular $c_L(\cdot,\cdot^i)\in D({\Omega}')\subset
\mathcal
{L}'_{\Vert\cdot\Vert}$ for all $i\in G$ and thus the rates are quasilocal.

Later we will need even stronger regularity of the rates in the
following sense.
%
\begin{lem}\label{ExponentialRegularityRotationRates}
Suppose $\sup_{i\in G}\sum_{A\ni i}\sum_{k\in G}e^{\varrho
(i,k)}\delta
_k(\Phi_A)<\infty$ and
%
\begin{equation}
\label{expdecay} \bar c_\varrho:=\sup_{i\in G}\sum
_{j\neq i}e^{\varrho (i,j)}\widebar C_{ij}<1
\end{equation}
then $\sup_{i\in G}\sum_{j\neq i}e^{\varrho (i,j)}\delta
_j(c_L(\cdot,\cdot^i))<\infty$.
\end{lem}
Notice, the first condition given in the above lemma is independent of
the hidden temperature parameter $\beta $ and with $\varrho
(i,k):=\varepsilon |i-k|$
corresponds to condition~(\ref{ExpDecayCondition}) in Theorem~\ref
{MainTheorem}.
Condition (\ref{expdecay}) is the requirement on the fineness of
discretization $q\geq q_0(\Phi)$ formulated in Theorem~\ref{MainTheorem}.

\begin{pf*}{Proof of Lemma \ref{ExponentialRegularityRotationRates}}
As a consequence of the exponential decay condition on
the Dobrushin matrix (\ref{expdecay})
(for\vspace*{1pt} a translation-invariant semi-metric $\varrho $ on $G$), we have
$\sup_{i\in G}\sum_{j\in G}e^{\varrho (i,j)}\widebar D_{ij}\leq\frac
{1}{1-\bar c_\varrho
}$ and by the triangle inequality
%
\begin{eqnarray}
\label{welldefinedexpdecay} \sup_{i\in G}\sum_{j\neq i}e^{\varrho (i,j)}
\delta_j \bigl(c_L \bigl(\cdot,\cdot^i
\bigr) \bigr)&\leq& \sup_{i\in G}\sum_{j\in G\setminus i}
\sum_{k\in G\setminus
i}e^{\varrho
(i,j)}\delta_k
\bigl(\psi^{i} \bigr)\widebar D_{kj}
\nonumber\\[-8pt]\\[-8pt]
&\leq&\frac{1}{1-\bar c_\varrho }\sup
_{i\in G}\sum_{k\in
G\setminus i}e^{\varrho
(i,k)}
\delta_k \bigl(\psi^{i} \bigr).\nonumber
\end{eqnarray}
%
But for $\psi^i_1$ and $\psi^i_2$ using the same arguments as in the
proof of Lemma~\ref{boundedvariation},
%
\begin{eqnarray}
\label{welldefinedexpdecay1}
\sum_{k\in G}e^{\varrho (i,k)}
\delta_k \bigl(\psi^i_2 \bigr)&\leq& C \sum
_{k\in
G}e^{\varrho
(i,k)}\delta_k \bigl(
\psi^i_1 \bigr)
\nonumber\\[-8pt]\\[-8pt]
&\leq& Ce^{2K} \sum
_{A\ni{i}}\sum_{k\in
G}e^{\varrho
(i,k)}\delta_k(\Phi_A)<\infty.\nonumber
\end{eqnarray}\upqed
\end{pf*}

Instead of imposing an exponential decay property of the Dobrushin
matrix, one could just consider polynomial weights $p(\varrho (i,j))$ which
would admit Hamiltonians with polynomial dependence. In fact, for our
purposes, that would be sufficient.

After these preparations we are in the position to use Theorem 3.9 of
\cite{Li85} and assert:
(1) The closure $\widebar{L}$ of $L$ is a Markov generator of a Markov
semiproup $(S^L_t)_{t\geq0}$ connected\vspace*{1pt} to the generator via the
Hille--Yosida theorem. $D({\Omega}')$ is a core for $\widebar{L}$.
(2) For observables $f\in D({\Omega}')$ we can control the oscillation of
$S_tf$ at any site $i\in G$ via
\[
\delta_i(S_tf)\leq \bigl[e^{t\Gamma }
\delta_\cdot(f) \bigr](i),
\]
where $\Gamma\dvtx  l_1\to l_1$, $[\Gamma \delta_\cdot(f)](i):=\sum_{j\neq i}\delta
_i(c_L(\cdot,\cdot^j))\delta_j(f)$ is a bounded operator with $\|\Gamma \|=:M$. In particular for $f\in D({\Omega}')$ we have $\|\!| S_tf\|\!|\leq
e^{tM}\|\!| f\|\!|$ and thus $S_tf\in D({\Omega}')$.

\subsection{Rotation property of the generator}\label{sec3.2}\label{Rotationpropertyofthegenerator}

The goal of this subsection is to verify Theorem~\ref
{TheRotationDynamics}, part (2). We use the following strategy:
\begin{longlist}[(3)]
\item[(1)] We verify the rotation property for infinitesimal times by
comparing the generator to the derivative on the level of the
probability density. We do this directly on local observables.
\item[(2)] In order to get from infinitesimal to finite time, we
consider the associated semigroup $(S^L_t)_{t\geq0}$ and use Taylor's
expansion. To match the first-order terms it is necessary to verify the
infinitesimal rotation for local functions propagated by $S^L_t$. Those
functions are no longer strictly local but lie in a larger space,
namely $\mathcal{L}'_{\|\!|\cdot\|\!|}$.
Since later we need (and will verify) the stronger result $S^L_tf\in
D_{p(\varrho )}({\Omega}')$ for\vspace*{1pt} local $f$ and weight-function
$p(x)=x^2$, at this
point we just assume $S^L_tf\in\mathcal{L}'_{\|\!|\cdot\|\!|}$.\vspace*{1pt}
\item[(3)] The two second-order error terms need to be estimated. As
for the first one we can use the contraction property of the semigroup.
For the other one we compute the second derivative of the measure again
on the level of the probability density and local observables. It turns
out the desired upper bound exists as long as the observable lies in a
space of weighted triple-normed functions.
\item[(4)] By assuming exponential decay of the Dobrushin matrix [see
(\ref{expdecay})] the rates of the generator are elements of this
space, even for arbitrary polynomial weights. One can think of these
spaces as containing functions with a certain degree of locality. The
amount of nonlocality the semigroup injects into a local function is
controlled by the degree of locality of the rates. This can simply be
captured by looking at the operator $\Gamma $ mentioned above. We can show
under these assumptions that local observables propagated by the
semigroup stay in the space of weighted triple-normed functions.
\end{longlist}

Let us start with an infinitesimal rotation and show $\mu
'_{(t+s)_{\operatorname{mod} 2\pi}}(f)=\mu'_{t}(S^L_sf)$ for all $t\in
[0,2\pi
)$, $s>0$, $\mu
'_\varphi =T\mu_\varphi \in\operatorname{ex} \mathcal{G}(\gamma
')$ and local observables
$f$ on
${\Omega}'$. Since the coarse-graining is finite it suffices to use
$f=1_{a_{\Lambda }}$ for finite $\Lambda $ and $a_\Lambda \in\{
1,\ldots,q\}^\Lambda $. Write $\rho
_{\Lambda }=d\gamma _\Lambda ^\Phi/d\lambda ^\Lambda $ for the
Lebesgue density of the local
specification in $\Lambda $. Then we can proceed similary to the intuitive
calculations done in (\ref{Intuitive}) and write
\begin{eqnarray*}
\label{differential1}
&& \frac{d}{d\varepsilon }_{|_{\varepsilon =0}}\mu'_{t+\varepsilon
}(1_{a_{\Lambda }})
\\
&&\qquad =
\int\mu_t (d\omega )\frac
{d}{d\varepsilon }_{|_{\varepsilon =0}} \biggl(\prod
_{i\in\Lambda }\int_{a_i|^l-\varepsilon }^{a_i|^r-\varepsilon }
\biggr) \,d\varphi _{\Lambda } \rho _{\Lambda }(\varphi _{\Lambda },
\omega )
\\
&&\qquad =\sum_{j\in\Lambda } \int\mu_t (d\omega )
\biggl(\prod_{i\in\Lambda \setminus j}\int_{a_i|^l}^{a_i|^r}
\biggr) \,d\varphi _{\Lambda \setminus j} \bigl(\rho _\Lambda
\bigl(a_j|^l,\varphi _{\Lambda \setminus
j},\omega \bigr)
\\
&&\hspace*{189pt} {}-\rho
_{\Lambda } \bigl(a_j|^r,\varphi _{\Lambda \setminus
j},
\omega \bigr) \bigr)
\end{eqnarray*}
since $\mu_t$ admits $\gamma _{\Lambda }$ for all $t\in[0,2\pi)$.
On the other hand,
\begin{eqnarray*}
\label{differential2}
&& \frac{d}{d\varepsilon }_{|_{\varepsilon =0}}\mu '_{t}S^L_{\varepsilon }(1_{a_{\Lambda }})
\\
&&\qquad =
\mu'_{t}(L1_{a_{\Lambda
}})
\\
&&\qquad =\sum_{j\in\Lambda } \biggl(\sum
_{\omega '\dvtx  \omega '^j_{\Lambda
}=a_{\Lambda }}c_L \bigl(\omega ', \bigl(
\omega ' \bigr)^j \bigr)\mu '_{t}
\bigl(\omega ' \bigr)-\sum_{\omega '\dvtx  \omega '_{\Lambda }=a_{\Lambda
}}c_L
\bigl(\omega ', \bigl(\omega ' \bigr)^j
\bigr)\mu'_{t} \bigl(\omega ' \bigr)
\biggr).
\end{eqnarray*}
Looking at the individual summands we find
\begin{eqnarray*}
\label{differential3}
&& \sum_{\omega '\dvtx  \omega '_{\Lambda }=a_{\Lambda }}c_L \bigl(
\omega ', \bigl(\omega ' \bigr)^j \bigr)
\mu'_{t} \bigl(\omega ' \bigr)
\\[-1pt]
&&\qquad =\int\mu
'_{t} \bigl(d\omega '
\bigr)1_{a_{\Lambda }} \bigl(\omega ' \bigr)\frac{\mu_{G\setminus
j}[\omega '_{G\setminus
j}](e^{-H_{j}(\omega '_j|^r,\cdot_{j^c})})
}{\mu_{G\setminus j}[\omega '_{G\setminus j}](\lambda ^j ( e^{-H_{j}}
1_{\omega '_{j}}))}
\\[-1pt]
&&\qquad =\int\mu'_{t}
\bigl(d\omega ' \bigr)1_{a_{\Lambda }} \bigl(\omega
' \bigr)\frac{\mu
_{G\setminus\Lambda }[\omega
'_{G\setminus\Lambda }](\lambda ^{\Lambda \setminus
j}(e^{-H_{\Lambda }(a_j|^r,\cdot_{\Lambda
\setminus j},\cdot_{\Lambda ^c})}1_{a_{\Lambda \setminus j}}))
}{\mu_{G\setminus\Lambda }[\omega '_{G\setminus\Lambda
}](\lambda ^\Lambda ( e^{-H_{\Lambda }}1_{a_{\Lambda
}}))}
\\[-1pt]
&&\qquad =\int\mu'_{t} \bigl(d\omega '
\bigr)1_{a_{\Lambda }} \bigl(\omega ' \bigr)\mu _{G}
\bigl[\omega '_{G} \bigr] \biggl(\frac{\lambda ^{\Lambda
\setminus j}(e^{-H_{\Lambda }(a_j|^r,\cdot_{\Lambda \setminus
j},\cdot_{\Lambda
^c})}1_{a_{\Lambda \setminus j}})
}{\lambda ^\Lambda ( e^{-H_{\Lambda }}1_{a_{\Lambda }})}
\biggr)
\\[-1pt]
&&\qquad =\int\mu'_{t}
\bigl(d\omega ' \bigr)1_{a_{\Lambda }} \bigl(\omega
' \bigr)\mu _{G} \bigl[\omega '_{G}
\bigr]
\\[-1pt]
&&\hspace*{41pt}{}\times  \biggl(\frac{1
}{\gamma _{\Lambda }(1_{a_{\Lambda }}|\cdot)} \biggl(\prod_{i\in\Lambda
\setminus j}
\int_{a_i|^l}^{a_i|^r} \biggr) \,d\varphi _{\Lambda \setminus j}
\bigl( \rho _{\Lambda } \bigl(a_j|^r,\varphi
_{\Lambda \setminus
j},\cdot \bigr) \bigr) \biggr)
\\[-1pt]
&&\qquad =\int\mu'_{t} \bigl(d\omega ' \bigr)
\mu_{G} \bigl[\omega '_{G} \bigr] \biggl(
\frac
{1_{a_{\Lambda }}
}{\gamma _{\Lambda }(1_{a_{\Lambda }}|\cdot)} \biggl(\prod_{i\in\Lambda
\setminus j}\int
_{a_i|^l}^{a_i|^r} \biggr) \,d\varphi _{\Lambda \setminus j} \bigl(
\rho _{\Lambda } \bigl(a_j|^r,\varphi
_{\Lambda \setminus
j},\cdot \bigr) \bigr) \biggr)
\\[-3pt]
&&\qquad =\int\mu_{t}(d\omega ) \biggl(\prod_{i\in\Lambda \setminus j}
\int_{a_i|^l}^{a_i|^r} \biggr) \,d\varphi _{\Lambda
\setminus j}
\bigl(\rho _{\Lambda } \bigl(a_j|^r,\varphi
_{\Lambda \setminus
j},\omega \bigr) \bigr),
\end{eqnarray*}
where we used the DLR equation in the second last line and the fact that
%
\begin{equation}
\label{differential4} \frac{\mu_{G\setminus\Lambda }[\omega '_{G\setminus\Lambda
}](\varphi (\cdot_{\Lambda ^c})\lambda ^\Lambda
( e^{-H_{\Lambda }}
1_{\omega '_{\Lambda }}))
}{\mu_{G\setminus\Lambda }[\omega '_{G\setminus\Lambda
}](\lambda ^\Lambda ( e^{-H_{\Lambda }}
1_{\omega '_{\Lambda }}))} =\mu_{G} \bigl[\omega
'_{G} \bigr] \bigl(\varphi (\cdot _{\Lambda ^c})
\bigr).
\end{equation}

We proceed similarly for the other summand.
%
%
Thus we have $\frac{d}{d\varepsilon }_{|_{\varepsilon =0}}\mu
'_{t+\varepsilon }(f)=\mu'_t(Lf)$ for\vspace*{1pt}
all local observables. Later we want to apply $S^L_tf$ and will show
$S^L_tf\in\mathcal{L}'_{\|\!|\cdot\|\!|}$ if $f$ is local. So let us
prove the following proposition.
%
\begin{prop}\label{differentialontriplecompletion}
If $f\in\mathcal{L}'_{\|\!|\cdot\|\!|}$, then $\frac
{d}{d\varepsilon }_{|_{\varepsilon
=0}}\mu'_{t+\varepsilon }(f)=\mu'_t(Lf)$.
\end{prop}

\begin{pf} Assume $(f_n)_{n\in{\mathbb N}}$ to be a sequence of local
functions such that $\|\!| f-f_n\|\!|\to0$ for $n\to\infty$. Then we
have, according to Proposition 3.2 of~\cite{Li85},
%
\begin{equation}
\label{differentialontriplecompletion1} \bigl|\mu'_t(Lf)-\mu'_t(Lf_n)\bigr|
\leq\|Lf-Lf_n\|\leq C\|\!| f-f_n\|\!|
\xrightarrow{n \to\infty} 0.
\end{equation}
On the other hand, with $g:=f-f_n$ and $o\dvtx G\to{\mathbb N}$ an ordering of
$G$, we have
\begin{eqnarray*}
\label{differentialontriplecompletion2}
&& \mu'_{t+\varepsilon }(g) -\mu'_{t}(g)
\\[-3pt]
&&\qquad =\int\mu_{t}(d\tilde\omega ) \bigl(g \bigl(T(\tilde\omega - \varepsilon 1_{G}) \bigr)-g \bigl(T(\tilde\omega ) \bigr)
\bigr)
\\[-3pt]
&&\qquad =\int\mu_{t}(d\tilde\omega )\sum_{j\in G}
\bigl(g \bigl(T(\tilde\omega - \varepsilon 1_{\{0,\ldots,o(j)\}}) \bigr)-g \bigl(T(
\tilde\omega - \varepsilon 1_{\{0,\ldots,(o(j)-1)\}
}) \bigr) \bigr)
\\[-3pt]
&&\qquad \leq\sum
_{j\in G}\delta_{j}(g)\mu_{t} \bigl(
\bigl\{\tilde\omega\dvtx  T(\tilde\omega _{j}-\varepsilon )=T(\tilde\omega _{j})-1 \bigr\} \bigr),
\end{eqnarray*}
where we use a telescopic sum in the second line. Further we have with
$A_j:=\{\tilde\omega\dvtx  T(\tilde\omega _{j}-\varepsilon )=T(\tilde\omega _{j})-1\}=\{\tilde\omega\dvtx
\tilde\omega _{j}\in[a|^l,a|^l+\varepsilon ]\mbox{ for some }a\in
\{1,\ldots,q\}\}$,
%
\begin{eqnarray}
\label{differentialontriplecompletion3} \mu_{t}(A_j)&\leq&\sup
_{\omega \in{\Omega}}\gamma _j(A_j|\omega )\leq
\frac{\varepsilon q e^{\|H_j\|}}{2\pi e^{-\|H_j\|}}\leq\widebar K\varepsilon
\end{eqnarray}
uniformly in $t$ and $j$, hence $\frac{1}{\varepsilon }|\mu
'_{t+\varepsilon }(f-f_n)-\mu
'_{t}(f-f_n)|\leq\widebar K \|\!| f-f_n\|\!|\xrightarrow{n\to\infty} 0$,
and we can conclude
\begin{eqnarray*}
\label{differentialontriplecompletion4}
&&\biggl|\frac{d}{d\varepsilon }_{|_{\varepsilon =0}}
\mu'_{t+\varepsilon
}(f)-\mu'_t(Lf)\biggr|
\\
&&\qquad \leq\biggl|
\frac{d}{d\varepsilon
}_{|_{\varepsilon =0}}\mu'_{t+\varepsilon }(f-f_n)\biggr|+\bigl|
\mu '_t \bigl(L(f-f_n) \bigr)\bigr|
\xrightarrow{n \to \infty} 0.
\end{eqnarray*}\upqed
\end{pf}\eject

Assume\vspace*{2pt} for the moment $S^L_tf\in\mathcal{L}'_{\|\!|\cdot\|\!|}$ for
local $f$. In order to verify the rotation property for finite times,
we use the following iteration procedure. Let $f$ be local, $k\in{\mathbb N}$,
$t\in[0,2\pi)$, $s>0$ and $\varepsilon:=s/k$. On the one hand,
%
\begin{eqnarray}
\label{Iteration} \mu'_t \bigl(S^L_sf
\bigr)&=&\mu'_t \bigl(S^L_\varepsilon
S^L_{s-\varepsilon }f \bigr)=\mu '_t
\bigl( \bigl(1+\varepsilon L +S^L_\varepsilon - (1+\varepsilon L)
\bigr)S^L_{s-\varepsilon }f \bigr)
\nonumber\\[-8pt]\\[-8pt]
&=&\mu'_t(g)+\varepsilon \mu'_t(Lg)+
\mu'_t \bigl( \bigl(S^L_\varepsilon
- (1+\varepsilon L) \bigr)g \bigr),\nonumber
\end{eqnarray}
where we set $g:=S^L_{s-\varepsilon }f$. On the other hand we can use Taylor's
expansion in Lagrange form and write
%
\begin{eqnarray}
\label{Iteration2} \mu'_{t+\varepsilon }(g)&=&\mu'_t(g)+
\varepsilon \frac
{d}{d\varepsilon }_{|_{\varepsilon =0}}\mu'_{t+\varepsilon
}(g)+
\frac
{\varepsilon ^2}{2}\frac{d^2}{d\hat\varepsilon ^2}_{|_{\tilde \varepsilon \in[0,\varepsilon ]}}\mu'_{t+\hat\varepsilon
}(g)
\nonumber\\[-8pt]\\[-8pt]
&=&\mu'_t(g)+\varepsilon \mu'_t(Lg)+
\frac{\varepsilon ^2}{2}\frac
{d^2}{d\hat\varepsilon ^2}_{|_{\tilde \varepsilon \in[0,\varepsilon ]}}\mu'_{t+\hat\varepsilon }(g).\nonumber
\end{eqnarray}

By iteration
%
\begin{eqnarray}
\label{Iteration3} 
\mu'_t
\bigl(S^L_sf \bigr)-\mu'_{t+s}(f)&=&
\sum_{l=0}^{k-1} \mu'_{t+l\varepsilon
}
\bigl( \bigl(S^L_{\varepsilon } - (1+\varepsilon L)
\bigr)S^L_{s-(l+1)\varepsilon }f \bigr)
\nonumber\\[-8pt]\\[-8pt]
&&{} -\frac{\varepsilon ^2}{2}\sum_{l=0}^{k-1}
\frac{d^2}{d\hat\varepsilon ^2}_{|_{\tilde\varepsilon \in
[0,\varepsilon ]}}\mu'_{t+l\varepsilon +\hat\varepsilon
}
\bigl(S^L_{s-(l+1)\varepsilon }f \bigr),\nonumber
\end{eqnarray}
where the error terms should go to zero as $k$ tends to infinity.
Let us look at the first error term on the RHS of (\ref{Iteration3})
and use the uniform continuity of the Markov semigroup, we have
\begin{eqnarray*}
\label{SecondError2}
&& \sum_{l=0}^{k-1}
\mu'_{t+l\varepsilon } \bigl( \bigl(S^L_{\varepsilon } -
(1+\varepsilon L) \bigr)S^L_{s-(l+1)\varepsilon }f \bigr) 
\\
&&\qquad \leq \varepsilon \sum_{l=0}^{k-1}\biggl\|
\frac{S^L_{\varepsilon
}S^L_{s-(l+1)\varepsilon }f-S^L_{s-(l+1)\varepsilon
}f}{\varepsilon }-LS^L_{s-(l+1)\varepsilon }f\biggr\|
\\
&&\qquad \leq\varepsilon \sum
_{l=0}^{k-1}\biggl\|\frac{S^L_{\varepsilon
}f-f}{\varepsilon }-Lf\biggr\|=s\biggl\|
\frac{S^L_{\varepsilon
}f-f}{\varepsilon }-Lf\biggr\|,
\end{eqnarray*}
where the RHS goes to zero as $\varepsilon $ goes to zero since the
semigroup is
generated by $L$ and $f$ in the domain of $L$. In particular this is
true for core observables of $L$.

Let us check the second error term on the RHS of (\ref{Iteration3}).
Let $t'(l)\in[t+l\varepsilon,t+(l+1)\varepsilon ]$. Then it
suffices to find a constant
$C(s,f)$ such that
%
\begin{equation}
\label{SecondErrorBound} \frac{d^2}{d\hat\varepsilon ^2}_{|_{\hat\varepsilon =0}}\mu '_{t'(l)+\hat\varepsilon }
\bigl(S^L_{s-(l+1)\varepsilon
}f \bigr)\leq C(s,f)
\end{equation}
for all $l$, since then we have
\[
\varepsilon ^2\sum_{l=0}^{k-1}
\frac{d^2}{d\hat\varepsilon
^2}_{|_{\hat\varepsilon =0}}\mu'_{t'(l)+\hat\varepsilon }
\bigl(S^L_{s-(l+1)\varepsilon }f \bigr)\leq\frac
{s^2}{k}C(s,f)
\xrightarrow{k\to\infty}0.
\]

Consider the second derivative when we apply the extremal Gibbs measure
at first to a local indicator function $1_{a_{\Lambda }}$. Then we
have
\begin{eqnarray*}
\label{SecondDerivativ}
\hspace*{-5pt}&& \frac
{d^2}{d\varepsilon ^2}_{|_{\varepsilon =0}}\mu'_{t+\varepsilon
}(1_{a_{\Lambda }})
\\[-2pt]
\hspace*{-5pt}&&\qquad =\sum_{j\in\Lambda } \int\mu_t
(d\omega )\frac{d}{d\varepsilon
} \biggl(\prod_{i\in\Lambda \setminus j}\int
_{a_i|^l-\varepsilon }^{a_i|^r-\varepsilon } \biggr) \,d\varphi _{\Lambda
\setminus j} \bigl(
\rho _\Lambda \bigl(a_j|^l-\varepsilon,\varphi
_{\Lambda
\setminus j},\omega \bigr)
\\[-2pt]
\hspace*{-5pt}&&\hspace*{214pt}{} -\rho _{\Lambda } \bigl(a_j|^r-
\varepsilon,\varphi _{\Lambda \setminus j},\omega \bigr) \bigr)
\\[-2pt]
\hspace*{-5pt}&&\qquad =\sum_{j\in\Lambda } \int\mu_t (d\omega )
\\[-2pt]
\hspace*{-5pt}&&\hspace*{59pt}{}\times
\biggl[ \biggl(\prod_{i\in
\Lambda \setminus j}\int_{a_i|^l}^{a_i|^r}
\biggr) \,d\varphi _{\Lambda \setminus j}
\\[-2pt]
\hspace*{-5pt}&&\hspace*{76pt}{}\times  \biggl(\rho _{\Lambda }
\bigl(a_j|^r, \varphi _{\Lambda \setminus
j},\omega \bigr)
\frac{d}{d\varepsilon _j}_{|_{\varepsilon
_j=a_j|^r}}H_\Lambda (\varepsilon _j,
\varphi _{\Lambda \setminus
j},\omega )
\\[-2pt]
\hspace*{-5pt}&&\hspace*{60pt}\quad\qquad{} -\rho _\Lambda \bigl(a_j|^l,
\varphi _{\Lambda \setminus
j},\omega \bigr)\frac{d}{d\varepsilon _j}_{|_{\varepsilon
_j=a_j|^l}}H_\Lambda
(\varepsilon _j,\varphi _{\Lambda \setminus
j},\omega ) \biggr)
\\[-2pt]
\hspace*{-5pt}&&\hspace*{76pt}{} +\sum_{k\in\Lambda \setminus j} \biggl(\prod
_{i\in\Lambda \setminus \{
j,k\}}\int_{a_i|^l}^{a_i|^r} \biggr) \,d
\varphi _{\Lambda \setminus j} \bigl(\rho _\Lambda \bigl(a_j|^l,a_k|^l,
\varphi _{\Lambda \setminus \{
j,k\}},\omega \bigr)
\\[-2pt]
\hspace*{-5pt}&&\hspace*{220pt}{} -\rho _\Lambda
\bigl(a_j|^l,a_k|^r, \varphi
_{\Lambda \setminus \{j,k\}},\omega \bigr)
\\[-2pt]
\hspace*{-5pt}&&\hspace*{220pt}{} -\rho _{\Lambda } \bigl(a_j|^r,a_k|^l,
\varphi _{\Lambda
\setminus \{j,k\}},\omega \bigr)
\\[-2pt]
\hspace*{-5pt}&&\hspace*{220pt}{}+\rho _{\Lambda
} \bigl(a_j|^r,a_k|^r,
\varphi _{\Lambda \setminus \{j,k\}},\omega \bigr) \bigr) \biggr]
\\[-2pt]
\hspace*{-5pt}&&\qquad =:\sum_{j\in\Lambda } \int\mu_t (d\omega )
\biggl[A(j,a_\Lambda,\omega )+\sum_{k\in\Lambda \setminus
j}B(j,k,a_\Lambda,\omega ) \biggr],
\end{eqnarray*}
where, as we see from the formular, the Hamiltonian of the first-layer
system needs to be differentiable as a function on $S^1$.
Let us assume these partial derivatives are also uniformly bounded\vadjust{\eject} with
$K':=\sup_{i\in G}\sup_{\omega \in{\Omega}}\|\frac
{d}{d\varepsilon }H_i(\varepsilon,\break \omega _{i^c})\|<\infty$. Then we have
{\fontsize{11pt}{11pt}\selectfont{\begin{eqnarray*}
\label{SecondDerivativ1}
\hspace*{-4pt}&& A(j,a_\Lambda,\omega )
\\
\hspace*{-4pt}&&\!\qquad = \biggl(\prod
_{i\in\Lambda \setminus j}\int_{a_i|^l}^{a_i|^r} \biggr) \,d
\varphi _{\Lambda \setminus
j} \biggl(\frac{e^{-H_\Lambda (a_j|^r,\varphi _{\Lambda \setminus
j},\omega _{\Lambda ^c})}}{\int d\varphi _{\Lambda
}e^{-H_\Lambda (\varphi _\Lambda,\omega _{\Lambda ^c})}}
\biggl(\frac{d}{d\varepsilon
_j}_{|_{\varepsilon _j=a_j|^r}}H_\Lambda (
\varepsilon _j,\varphi _{\Lambda \setminus j},\omega )
\\
\hspace*{-4pt}&&\!\hspace*{225pt}{} -\frac{d}{d\varepsilon
_j}_{|_{\varepsilon _j=a_j|^l}}H_\Lambda (\varepsilon _j,
\varphi _{\Lambda \setminus
j}\omega ) \biggr)
\\
\hspace*{-4pt}&&\!\hspace*{130pt}{} +\frac{d}{d\varepsilon _j}_{|_{\varepsilon _j=a_j|^l}}H_\Lambda (\varepsilon
_j,\varphi _{\Lambda \setminus j},\omega )
\\
\hspace*{-4pt}&&\!\hspace*{159pt}{} \times \frac
{e^{-H_\Lambda (a_j|^r,\varphi _{\Lambda \setminus j},\omega
_{\Lambda ^c})}-e^{-H_\Lambda (a_j|^l,\varphi _{\Lambda \setminus
j},\omega _{\Lambda ^c})}}{\int d\varphi _{\Lambda }e^{-H_\Lambda
(\varphi _\Lambda,\omega _{\Lambda ^c})}} \biggr),
\end{eqnarray*}}}%
where $\frac{d}{d\varepsilon _j}_{|_{\varepsilon
_j=a_j|^r}}H_\Lambda (\varepsilon _j,\varphi _{\Lambda \setminus
j},\omega
)-\frac{d}{d\varepsilon _j}_{|_{\varepsilon _j=a_j|^l}}H_\Lambda
(\varepsilon _j,\varphi _{\Lambda \setminus j}\omega )
\leq\delta_j(\frac{d}{d\varepsilon _j}H_j)\leq2K'$ and
$e^{-H_\Lambda (a_j|^r,\varphi _{\Lambda \setminus j},\omega
_{\Lambda ^c})}-e^{-H_\Lambda (a_j|^l,\varphi _{\Lambda \setminus
j},\omega _{\Lambda ^c})}\leq2e^{K}e^{-\sum_{A\cap\Lambda \neq
\varnothing, j\notin
A}\Phi_A(\varphi _{\Lambda \setminus j},\omega _{\Lambda ^c})}$. Thus
%
%
%
\begin{eqnarray*}
\label{SecondDerivativ1c}
&& A(j,a_\Lambda,\omega )
\\
&&\qquad \leq 2K'e^{2K}2\pi\frac{(\prod_{i\in
\Lambda \setminus j}\int_{a_i|^l}^{a_i|^r}) \,d\varphi _{\Lambda \setminus j}e^{-\sum_{A\cap
\Lambda \neq\varnothing,
j\notin A}\Phi_A(\varphi _{\Lambda \setminus j},\omega _{\Lambda
^c})}}{\int d\varphi _{\Lambda \setminus
j}e^{-\sum_{A\cap\Lambda \neq\varnothing,  j\notin A}\Phi
_A(\varphi _{\Lambda \setminus j},\omega
_{\Lambda ^c})}}
\\
&&\quad\qquad{} +2e^{2K} \bigl(K'+\bigl(|\Lambda |-1\bigr)K \bigr)
\\
&&\qquad\qquad{}\times 2\pi
\frac{(\prod_{i\in\Lambda
\setminus j}\int_{a_i|^l}^{a_i|^r}) \,d\varphi _{\Lambda \setminus j}e^{-\sum_{A\cap
\Lambda \neq\varnothing,
j\notin A}\Phi_A(\varphi _{\Lambda \setminus j},\omega _{\Lambda
^c})}}{\int d\varphi _{\Lambda \setminus
j}e^{-\sum_{A\cap\Lambda \neq\varnothing, j\notin A}\Phi
_A(\varphi _{\Lambda \setminus j},\omega
_{\Lambda ^c})}}.
\end{eqnarray*}
Let us check the second term. We can write
\begin{eqnarray*}
\label{SecondDerivativ2}
&& B(j,k,a_\Lambda,\omega )
\\
&&\qquad \leq 4e^{4K}4
\pi^2\frac{(\prod_{i\in
\Lambda \setminus \{j,k\}}\int_{a_i|^l}^{a_i|^r}) \,d\varphi _{\Lambda \setminus \{j,k\}}e^{-\sum
_{A\cap\Lambda \neq\varnothing,\{j,k\}\not\subset A}\Phi_A(\varphi _{\Lambda \setminus \{
j,k\}},\omega _{\Lambda
^c})}}{\int d\varphi _{\Lambda \setminus \{j,k\}}e^{-\sum_{A\cap
\Lambda \neq\varnothing, \{
j,k\}\not\subset A}\Phi_A(\varphi _{\Lambda \setminus \{j,k\}
},\omega _{\Lambda ^c})}}.
\end{eqnarray*}
For convenience set $\check K:=\max\{K,K'\}$ and $\widebar K:=\max\{4\pi
\check K e^{2\check K},8\pi^2 e^{4\check K}\}$. Also we want to adopt a
notation we introduced earlier, namely,
\begin{eqnarray*}
&& \gamma _{\Lambda \setminus j}|_{j^c}(1_{a_{\Lambda \setminus
j}}|\omega
_{\Lambda ^c})
\\
&&\qquad =\frac{(\prod_{i\in\Lambda \setminus
j}\int_{a_i|^l}^{a_i|^r}) \,d\varphi _{\Lambda \setminus j}e^{-\sum
_{A\cap\Lambda \neq\varnothing, j\notin A}\Phi_A(\varphi _{\Lambda \setminus j},\omega
_{\Lambda ^c})}}{\int d\varphi _{\Lambda \setminus
j}e^{-\sum_{A\cap\Lambda \neq\varnothing, j\notin A}\Phi
_A(\varphi _{\Lambda \setminus j},\omega
_{\Lambda ^c})}}.
\end{eqnarray*}
Before we combine these estimates, let us apply the measure to a
general local function $h$ on the coarse-grained space with support
$\Lambda
$. $h$ can be written as $h(\omega ')=\sum_{a_{\Lambda }\in\{1,\ldots, q\}^{\Lambda
}}\kappa_{a_{\Lambda }}1_{a_{\Lambda }}(\omega ')$ with $\|h\|=\sup_{a_\Lambda }|\kappa
_{a_{\Lambda }}|$. Hence
\begin{eqnarray*}
\label{SecondDerivativ3}
&& \frac{d^2}{d\varepsilon ^2}_{|_{\varepsilon =0}}\mu '_{t+\varepsilon }(h)
\\
&&\qquad \leq\|h\|\sum_{j\in\Lambda }
\int\mu_t (d\omega ) \biggl[q\widebar K\bigl(|\Lambda |+1\bigr)\sum
_{a_{\Lambda \setminus j}\in\{1, \ldots, q\}^{\Lambda \setminus
j}}\gamma _{\Lambda \setminus j}|_{j^c}(1_{a_{\Lambda
\setminus j}}|
\omega _{\Lambda ^c})
\\
&&\hspace*{117pt}{}+\sum_{k\in\Lambda \setminus j}q^2\widebar K
\\
&&\hspace*{127pt}{}\times \sum
_{a_{\Lambda \setminus \{j,k\}}\in\{
1, \ldots, q\}^{\Lambda \setminus \{j,k\}}}\gamma _{\Lambda
\setminus \{j,k\}}|_{\{j,k\}^c}(1_{a_{\Lambda \setminus \{
j,k\}}}|
\omega _{\Lambda ^c}) \biggr]
\\
&&\qquad \leq\|h\| |\Lambda | \bigl(q\widebar K \bigl(|\Lambda |+1\bigr)+q^2\widebar K
\bigl(| \Lambda |-1\bigr) \bigr)
\\
&&\qquad \leq \widehat K|\Lambda |^2\| h\|.
\end{eqnarray*}
For a general quasilocal function $f$, one can write again a telescopic
sum using an ordering of $G$ and a generic configuration $\eta'$
%
\begin{eqnarray}
\label{SecondDerivativ5} f \bigl(\omega ' \bigr)&=&f \bigl(\omega
_1',\eta'_{\{1\}^c} \bigr)+ \bigl(f
\bigl(\omega _1',\omega _2',
\eta'_{\{1,2\}^c} \bigr)-f \bigl(\omega _1',
\eta'_{\{1\}^c} \bigr) \bigr)
\nonumber\\[-8pt]\\[-8pt]
&&{} +\sum_{n\geq3} \bigl(f
\bigl( \omega '_{\{1,\ldots,n\}},\eta'_{\{1,\ldots,n\}
^c}
\bigr)-f \bigl(\omega '_{\{
1,\ldots,n-1\}},\eta'_{\{1,\ldots,n-1\}^c}
\bigr) \bigr).\nonumber
\end{eqnarray}
Let us define $g_n(\omega '):= (f(\omega '_{\{1,\ldots,n\}},\eta
'_{\{1,\ldots,n\}
^c})-f(\omega '_{\{1,\ldots,n-1\}},\eta'_{\{1,\ldots,n-1\}^c}) )\in
\mathcal
{F}_{\{1,\ldots,n\}}$.
In particular
$\| g_n\| \leq\delta_n(f)$.
Hence we can write
%
\begin{eqnarray}
\label{SecondDerivativ6} \frac{d^2}{d\varepsilon ^2}_{|_{\varepsilon =0}}\mu '_{t+\varepsilon }(f)
&\leq&\| f\| \widehat{K}+\sum
_{n\geq2}\| g_n\| \widehat{K}n^2\leq 2
\widehat K\sum_{n\geq1}n^2\delta_n(f).
\end{eqnarray}
Thus in order to have (\ref{SecondErrorBound}) it suffices to show
$S^L_tf\in D_{|\cdot|^2}({\Omega}')$ for local $f$.
To do that, let us use the exponential decay property of the Dobrushin
matrix introduced in (\ref{expdecay}) and the exponentially decaying
Hamiltonian, that is, by Lemma~\ref{ExponentialRegularityRotationRates}
assume the model to satisfy $\sup_{i\in G}\sum_{j\neq i}e^{\varrho
(i,j)}\delta_j(c_L(\cdot,\cdot^i))=:\check M_\varrho <\infty$ for some
translation-invariant increasing semi-metric $\varrho $ in $G$. With
this we
can prove the following proposition.
%
\begin{prop}\label{semigroupdiffusion}
Let $f$ be a local observable on ${\Omega}'$ and $(S^L_t)_{t\geq0}$
associated to the rotation generator $L$. For all polynomials $p$ on
${\mathbb R}
_0^+$ we have $S^L_tf\in D_{p(\varrho )}({\Omega}')$.
\end{prop}

\begin{pf}Let us consider only the monomials $x^n$. It suffices
to look at $n=2^m$ for some $m\in{\mathbb N}$.
We know from Theorem 3.9 in \cite{Li85},
\[
\sum_{i\geq0}\varrho (i,0)^m
\delta_i \bigl(S^L_tf \bigr)\leq\sum
_{i\geq
0}\varrho (i,0)^m
\bigl[e^{t\Gamma } \delta_\cdot(f) \bigr](i),
\]
where $[\Gamma \delta_\cdot(f)](i):=\sum_{j\neq i}\delta
_i(c_L(\cdot,\cdot
^j))\delta_j(f)$. There exists a constant $K_{m,\varrho }$ such that for
fixed $j,m\in{\mathbb N}$, we have $\varrho (i,j)^m\leq K_{m,\varrho
}e^{\varrho (i,j)}$. Of
course local \mbox{$f\in D_{\varrho ^m}({\Omega}')$} for all $m\in{\mathbb N}$
and also for
exponential weight. Under the above condition on the jump rates, the
operator $\Gamma $ is\vspace*{1pt} bounded as well in the exponential weighted
triple-norm with norm $\check M$, indeed,
%
\begin{eqnarray}
\label{operatortriple-semi-norm} \|\!|\Gamma \|\!|_{e^\varrho }&=&\sup_{\|\!| v\|\!|_{e^\varrho }\leq
1}
\frac{\|\!|\Gamma
v\|\!|_{e^\varrho }}{\|\!| v\|\!|_{e^\varrho }}=\sup_{\|\!| v\|\!
|_{e^\varrho }\leq
1}\frac{\sum_{i\geq0}\sum_{j\neq i}e^{\varrho (i,0)}\delta
_i(c_L(\cdot,\cdot
^j))v_j}{\sum_{j\geq0}e^{\varrho (j,0)}v_j}
\nonumber\\[-8pt]\\[-8pt]
&\leq&\sup_{\|\!| v\|\!|_{e^\varrho }\leq1}
\frac{\check M \sum_{j\geq
0}e^{\varrho (j,0)}v_j}{\sum_{j\geq0}e^{\varrho (j,0)}v_j}=\check M.\nonumber
\end{eqnarray}
Then we can write
\begin{eqnarray*}
\label{TripelNorm12} \sum_{i\geq0}\varrho
(i,0)^m \bigl[e^{t\Gamma }\delta_\cdot(f) \bigr](i)&\leq& K_{m,\varrho }\sum_{i\geq0}e^{\varrho (i,0)}
\bigl[e^{t\Gamma }\delta_\cdot (f) \bigr](i)=K_{m,\varrho }\big\|\hspace*{-2pt}\big|e^{t\Gamma}\delta_\cdot(f)\big\|\hspace*{-2pt}\big|_{e^\varrho }
\\
&\leq& K_{m,\varrho }\big\|\hspace*{-2pt}\big| e^{t\Gamma }\big\|\hspace*{-2pt}\big|_{e^\varrho }\big\|\hspace*{-2pt}\big|\delta
_\cdot(f)\big\|\hspace*{-2pt}\big| _{e^\varrho }
\\
&\leq& K_{m,\varrho }e^{t\|\!|\Gamma \|\!|_{e^\varrho
}}\big\|\hspace*{-2pt}\big|\delta_\cdot (f)\big\|\hspace*{-2pt}\big|_{e^\varrho }.
\end{eqnarray*}\upqed
\end{pf}

In particular for local $f$, we have $S^L_{s-\varepsilon }f\in
D_{p(|\cdot|)}({\Omega}
')\subset\mathcal{L}'_{\|\!|\cdot\|\!|}$ for all polynomial and even
exponential weights $p$. In other words, we can control the diffusion
of the semi-group applied to a local function by looking at the decay
property of the conditional Dobrushin matrix as well as of the
first-layer Hamiltonian. In particular if those are well behaved (which
is the case for the $XY$-model with some slightly refined
coarse-graining) the second order terms in the Taylor expansion are controlled.
We can conclude
$\mu'_{t+\varepsilon }=S^L_{\varepsilon }(\mu'_{t})$ for all
extremal Gibbs measures
labeled by $t\in S^1$ and $\varepsilon >0$.

\section{Reversible dynamics for discrete-spin models}\label{sec4}\label{4}

The infinite-volume dynamics $K$
%
%
%
%
%
%
given in (\ref{GlauberGenerator}) with rates satisfying (\ref{GlauberRates})
%
%
%
%
%
is reversible. By expressing the RHS of (\ref{GlauberRates}) in terms
of the specification $\gamma '$ it is clear that $K$ has detailed balance
with respect to $\gamma '$.


%
%
These rates are bounded (by boundedness of $H_j$),
translation-invariant (by the translation-covariance of the $\mu
_{G\setminus j}[\omega '_{G\setminus j}]$ in the conditional
Dobrushin regime) and of exponentially decaying influence (however, not
strictly local).
The rates are even uniformly bounded and bounded in the triple-norm by
the same arguments as used for the rotation dynamics, so
Proposition~\ref{TheGlauberDynamics}, part (1) is true.

In the next subsection we adapt a line of arguments presented for $q=2$
in \cite{Li85} for general finite $q$.

\subsection{Translation-invariant invariant measures are Gibbs measures}\label{sec4.1}

Let us put ourselves in dimension $d\geq3$. In the right temperature
region there are multiple Gibbs measures for the $XY$-model,
ferromagnetically ordered on $S^1$.

Since in the following subsections we will only deal with second-layer
configurations, it is convenient to suppress the primes and write
$c_K(\omega,\omega ^i)$ for the up-flip at site $i\in G$ and
$c_K(\omega,\omega ^{i-})$
for the down-flip. Assume the rates to be defined as in (\ref
{GlauberRates}), in particular for the corresponding process
second-layer Gibbs measures are invariant w.r.t. $K$.

We now show,
invariant measures w.r.t. $K$ that are also translation-invariant are
second-layer Gibbs measures. This is precisely part (2) of Proposition
\ref{TheGlauberDynamics}. We use Holleys's argument \cite{Ho71}.
Recall the definition of the second-layer specification and define the
local relative entropy
%
\begin{equation}
\label{RelativeEntropy1} H_\Lambda \bigl(\nu|\gamma '_\Lambda
(\cdot|\zeta) \bigr):=\sum_{\omega
\in\{1,\ldots,q\}^\Lambda }
\nu(1_\omega )\log \frac{\nu(1_\omega )}{\gamma '_\Lambda (1_\omega |\zeta)},
\end{equation}
%
%
%
where $\Lambda \subset{\mathbb Z}^d$ is finite, $\nu\in\mathcal
{P}({\Omega}')$ and $\zeta\in
{\Omega}'$ an arbitrary but fixed boundary condition. Let
$(S^K_t)_{t\geq0}$
be the semigroup for the generator $K$, and define $\nu_t:=S^K_t(\nu)$.
Let us compute $\frac{d}{dt}_{|_{t=0}}H_\Lambda (\nu_t|\gamma
'_\Lambda (\cdot|\zeta))$
in two steps,
\begin{eqnarray*}
\label{RelativeEntropy3}
&& \frac{d}{dt}_{|_{t=0}}\sum
_{\omega \in\{1,\ldots,q\}^\Lambda
}\nu_t(1_\omega )\log\nu
_t(1_\omega ) 
\\
&&\qquad =\sum
_{\omega } \bigl[1+\log\nu(1_\omega ) \bigr]\int
K1_\omega \,d\nu
\\
&&\qquad =\sum_{\omega, i\in\Lambda }\log\nu(1_\omega )\int\nu(d
\eta ) \bigl[c_K \bigl(\eta,\eta^i \bigr)
\bigl(1_\omega \bigl(\eta^i \bigr)-1_\omega (\eta)
\bigr)
\\
&&\hspace*{108pt}\quad\qquad {}+c_K \bigl(\eta,\eta^{i-} \bigr)
\bigl(1_\omega \bigl(\eta^{i-} \bigr)-1_\omega (\eta)
\bigr) \bigr]
\\
&&\qquad =\sum_{\omega, i\in\Lambda } \biggl[
\Gamma \bigl(\omega,i^+ \bigr)\log \frac{\nu(1_{\omega ^i})}{\nu(1_{\omega
})}+\Gamma \bigl(\omega,i^-
\bigr)\log\frac{\nu(1_{\omega ^{i-}})}{\nu
(1_\omega )} \biggr],
\end{eqnarray*}
where we wrote $\Gamma (\omega,i^\pm):=\int\nu(d\eta)c_K(\eta,\eta^{i\pm})1_{\omega
}(\eta)$ for the outflows of $1_\omega $ in the direction $i^\pm$.

Note that if $\nu$ is invariant w.r.t. $K$, then $\nu(1_{\omega
_\Lambda })>0$ for
all $\omega _\Lambda \in\{1,\ldots,q\}^\Lambda $, indeed
%
\begin{eqnarray}
\label{PositivityInvariantMeasures} 0&=&\int K1_{\omega _\Lambda }\,d\nu\nonumber
\\
&=&\sum_{i\in\Lambda }
\int d\nu(d\eta) \bigl[c_K \bigl(\eta,\eta ^{i} \bigr)
\bigl(1_{\omega _\Lambda ^{i-}}(\eta )-1_{\omega _\Lambda }(\eta) \bigr)
\\
&&\hspace*{62pt}{} +c_K
\bigl(\eta,\eta^{i-} \bigr) \bigl(1_{\omega
_\Lambda ^i}(\eta)-1_{\omega _\Lambda }( \eta ) \bigr) \bigr].\nonumber
\end{eqnarray}
Since all flip-rates are positive, $\nu(1_{\omega _\Lambda })=0$
would imply $\nu
(1_{\omega _\Lambda ^i})=0=\nu(1_{\omega _\Lambda ^{i-}})$ for all
$i\in\Lambda $ and thus by
iteration $\nu(1_{\eta_\Lambda })=0$ for all $\eta_\Lambda \in\{
1,\ldots,q\}^\Lambda $,
which is a contradiction to $\nu$ being a probability measure.

Let us look at the second summand of $\frac{d}{dt}_{|_{t=0}}H_\Lambda
(\nu
_t|\gamma '_\Lambda (\cdot|\zeta))$. Since the normalizing constant
in the
specification is independent of $\omega _\Lambda $ and\break  $\sum_{\omega _\Lambda }\frac
{d}{dt}_{|_{t=0}}\int\nu_t(d\omega )=0$, we can directly compute
%
\begin{eqnarray}
\label{RelativeEntropySecondPart1}
&& \frac{d}{dt}_{|_{t=0}}\int\nu_t(d\omega )\log\mu_{G\setminus
\Lambda }[\zeta _{G\setminus\Lambda }] \bigl(
\lambda^\Lambda \bigl(e^{-H_\Lambda }1_{\omega
_\Lambda } \bigr) \bigr)\nonumber
\\
&&\qquad =\sum_{\omega _\Lambda \in\{1,\ldots,q\}^\Lambda }\int\nu (d
\eta_{\Lambda ^c})\nu(\omega _\Lambda |\eta_{\Lambda
^c})\nonumber
\\
&&\hspace*{92pt}{}\times K\log
\mu_{G\setminus\Lambda }[\zeta_{G\setminus\Lambda
}] \bigl(\lambda^\Lambda
\bigl(e^{-H_\Lambda
}1_{\cdot\Lambda } \bigr) \bigr) (\omega _\Lambda
\eta_{\Lambda ^c})\nonumber
\\
&&\qquad =\sum_{\omega _\Lambda, i\in\Lambda }
\biggl[\log\frac{\mu
_{G\setminus\Lambda }[\zeta
_{G\setminus\Lambda }](\lambda^\Lambda (e^{-H_\Lambda }1_{\omega
^i_\Lambda }))}{\mu_{G\setminus\Lambda
}[\zeta_{G\setminus\Lambda }](\lambda^\Lambda (e^{-H_\Lambda
}1_{\omega _\Lambda }))}\int\nu(d\eta )1_{\omega _{\Lambda }}(
\eta)c_K \bigl(\eta,\eta^i \bigr)
\nonumber\\[-8pt]\\[-8pt]
&&\hspace*{32pt}\quad\qquad{}+\log\frac{\mu_{G\setminus\Lambda }[\zeta
_{G\setminus\Lambda
}](\lambda^\Lambda (e^{-H_\Lambda }1_{\omega ^{i-}_\Lambda
}))}{\mu_{G\setminus\Lambda }[\zeta
_{G\setminus\Lambda }](\lambda^\Lambda (e^{-H_\Lambda }1_{\omega
_\Lambda }))}\int\nu(d\eta)1_{\omega
_{\Lambda }}(
\eta)c_K \bigl(\eta,\eta^{i-} \bigr) \biggr]\nonumber
\\
&&\qquad =\sum_{\omega _\Lambda, i\in\Lambda } \bigl[V \bigl(\omega,i^+ \bigr)\Gamma
\bigl(\omega,i^+ \bigr)+V \bigl(\omega,i^. \bigr)\Gamma \bigl(\omega,i^-
\bigr) \bigr],\nonumber
\end{eqnarray}
where we defined $V(\omega,i^\pm):=\log\frac{\mu_{G\setminus
\Lambda }[\zeta
_{G\setminus\Lambda }](\lambda^\Lambda (e^{-H_\Lambda }1_{\omega
^{i\pm}_\Lambda }))}{\mu
_{G\setminus\Lambda }[\zeta_{G\setminus\Lambda }](\lambda
^\Lambda (e^{-H_\Lambda }1_{\omega _\Lambda
}))}$. Notice we have
\begin{eqnarray*}
\frac{\mu_{G\setminus\Lambda }[\zeta_{G\setminus\Lambda
}](\lambda^\Lambda (e^{-H_\Lambda
}1_{\omega ^{i}_\Lambda }))}{\mu_{G\setminus\Lambda }[\zeta
_{G\setminus\Lambda }](\lambda^\Lambda
(e^{-H_\Lambda }1_{\omega _\Lambda }))}&=&\frac{\mu_{G\setminus
i}[\zeta_{G\setminus\Lambda
}\omega _{\Lambda \setminus i}](\lambda(e^{-H_i}1_{\omega
^i_i}))}{\mu_{G\setminus
i}[\zeta_{G\setminus\Lambda }\omega _{\Lambda \setminus
i}](\lambda(e^{-H_i}1_{\omega
_i}))}
\\
&=&\frac{c_K(\omega _{\Lambda }\zeta_{\Lambda ^c},\omega
^{i}_{\Lambda }\zeta_{\Lambda ^c})}{c_K(\omega ^{i}_{\Lambda
}\zeta_{\Lambda ^c},\omega _{\Lambda }\zeta_{\Lambda ^c})}.
\end{eqnarray*}

Combining the two summands we have
%
\begin{eqnarray}
\label{RelativeEntropyTogether}
&& \frac{d}{dt}H_\Lambda \bigl(\nu|\gamma'_\Lambda (\cdot|\zeta ) \bigr)_{|_{t=0}}\nonumber
\\
&&\qquad =
\sum_{\omega, i\in\Lambda } \biggl[\Gamma \bigl(\omega,i^+ \bigr)
\biggl(\log\frac{\nu(1_{\omega ^i})}{\nu
(1_{\omega })}-V \bigl(\omega,i^+ \bigr) \biggr)
\\
&&\hspace*{60pt}{} +\Gamma \bigl(
\omega,i^- \bigr) \biggl(\log\frac{\nu(1_{\omega ^{i-}})}{\nu
(1_\omega )}-V \bigl(\omega,i^- \bigr)
\biggr) \biggr].\nonumber
\end{eqnarray}
Since $\log\frac{\nu(1_{\omega ^i})}{\nu(1_{\omega })}-V(\omega,i^+)=-(\log\frac{\nu
(1_{\omega })}{\nu(1_{\omega ^i})}-V(\omega ^i,i^-))$, we can write
%
\begin{eqnarray}
\label{RelativeEntropyTogether2}
&& 2\frac{d}{dt}_{|_{t=0}}H_\Lambda \bigl(
\nu_t|\gamma '_\Lambda (\cdot |\zeta) \bigr)\nonumber
\\
&&\qquad =
\sum_{\omega, i\in
\Lambda } \biggl[ \bigl(\Gamma \bigl(\omega,i^+
\bigr)-\Gamma \bigl(\omega ^i,i^- \bigr) \bigr) \biggl(\log
\frac{\nu(1_{\omega ^i})}{\nu
(1_{\omega })}-V \bigl(\omega,i^+ \bigr) \biggr)
\\
&&\hspace*{28pt}\quad\qquad{} + \bigl(\Gamma \bigl(\omega,i^- \bigr)-\Gamma \bigl(\omega ^{i-},i^+
\bigr) \bigr) \biggl(\log\frac
{\nu(1_{\omega ^{i-}})}{\nu(1_\omega
)}-V \bigl(\omega,i^- \bigr) \biggr)
\biggr].\nonumber
\end{eqnarray}
Adding zeros we have
%
\begin{eqnarray*}
\label{RelativeEntropyTogether3}
&& 2\frac{d}{dt}H_\Lambda \bigl(\nu|\gamma
'_\Lambda (\cdot|\zeta ) \bigr)_{|_{t=0}}\nonumber
\\
&&\qquad = -\sum_{\omega, i\in\Lambda } \biggl[ \bigl(\Gamma \bigl(\omega,i^+
\bigr)-\Gamma \bigl(\omega ^i,i^- \bigr) \bigr)\log\frac{\Gamma (\omega,i^+)}{\Gamma (\omega ^i,i^-)}
\\
&&\hspace*{37pt}\quad\qquad{}+ \bigl(\Gamma \bigl(\omega,i^- \bigr)-\Gamma \bigl(\omega
^{i-},i^+ \bigr) \bigr)\log \frac{\Gamma (\omega,i^-)}{\Gamma (\omega
^{i-},i^+)} \biggr]
\\
&&\quad\qquad{} +\sum_{\omega, i\in\Lambda } \biggl[ \bigl(\Gamma \bigl(\omega,i^+
\bigr)-\Gamma \bigl(\omega ^i,i^- \bigr) \bigr)
\\
&&\hspace*{75pt}{} \times \biggl[\log
\frac
{\Gamma (\omega,i^+)}{\nu(1_{\omega })}-\log\frac{\Gamma
(\omega ^i,i^-)}{\nu(1_{\omega ^i})}-V \bigl(\omega,i^+ \bigr) \biggr]
\\
&&\hspace*{75pt}{}+ \bigl(\Gamma \bigl(\omega,i^- \bigr)-\Gamma \bigl(\omega
^{i-},i^+ \bigr) \bigr)
\\
&&\hspace*{85pt}{}\times \biggl[\log\frac{\Gamma (\omega,i^-)}{\nu(1_\omega )}-\log
\frac{\Gamma (\omega ^{i-},i^+)}{\nu
(1_{\omega ^{i-}})}-V \bigl(\omega,i^- \bigr) \biggr] \biggr].
\end{eqnarray*}
If $\nu$ is invariant w.r.t. $K$, it follows
\begin{eqnarray*}
\label{RelativeEntropyTogether4} &&\sum_{\omega, i\in\Lambda } \biggl[ \bigl(\Gamma
\bigl(\omega,i^+ \bigr)-\Gamma \bigl(\omega ^i,i^- \bigr) \bigr)\log
\frac{\Gamma (\omega,i^+)}{\Gamma (\omega ^i,i^-)}
\\
&&\hspace*{28pt}{}+ \bigl(\Gamma \bigl(\omega,i^- \bigr)-\Gamma \bigl(\omega
^{i-},i^+ \bigr) \bigr)\log \frac{\Gamma (\omega,i^-)}{\Gamma (\omega
^{i-},i^+)} \biggr]
\\
&&\qquad =\sum_{\omega, i\in\Lambda } \biggl[ \bigl(\Gamma \bigl(\omega,i^+
\bigr)-\Gamma \bigl(\omega ^i,i^- \bigr) \bigr)
\\
&&\hspace*{60pt}{}\times  \biggl[\log
\frac{\Gamma (\omega,i^+)}{\nu(1_{\omega })}-\log\frac{\Gamma (\omega ^i,i^-)}{\nu
(1_{\omega ^i})}-V \bigl(\omega,i^+ \bigr) \biggr]
\\
&&\hspace*{60pt}{}+ \bigl(\Gamma \bigl(\omega,i^- \bigr)-\Gamma \bigl(\omega
^{i-},i^+ \bigr) \bigr)
\\
&&\hspace*{70pt}{}\times  \biggl[\log\frac{\Gamma (\omega,i^-)}{\nu(1_\omega )}-\log
\frac{\Gamma (\omega ^{i-},i^+)}{\nu
(1_{\omega ^{i-}})}-V \bigl(\omega,i^- \bigr) \biggr] \biggr],
\end{eqnarray*}
where the left-hand side is nonnegative.
We want to exploit properties of the \mbox{$d$-}dimensional lattice in order
to show the RHS of the last equation
goes to zero for $\Lambda \nearrow G$. Let us define
%
\begin{eqnarray}
\label{DefinitionsHolleysArgument} \kappa_\Lambda \bigl(i^\pm \bigr)&:=&\sum
_{\omega } \bigl(\Gamma \bigl(\omega,i^\pm
\bigr)-\Gamma \bigl(\omega ^{i\pm},i^\mp \bigr) \bigr)\log
\frac
{\Gamma (\omega,i^\pm)}{\Gamma (\omega ^{i\pm},i^\mp)},\nonumber
\\[-2pt]
\beta_\Lambda \bigl(i^\pm
\bigr)&:=& \sum_{\omega }\bigl|\Gamma \bigl(\omega,i^\pm \bigr)-\Gamma \bigl(\omega ^{i\pm},i^\mp
\bigr)\bigr|,
\\[-2pt]
\vartheta _\Lambda (i)&:=&\sum_{j\notin\Lambda }
\sup_{\eta
_{j^c}=\tilde\eta_{j^c}}\frac
{|c_K(\eta,\eta^i)-c_K(\tilde\eta,\tilde\eta^i)|}{c_K(\eta,\eta
^i)}\nonumber
\\[-2pt]
&&{}+\sum_{j\notin\Lambda }\sup_{\eta_{j^c}=\tilde\eta
_{j^c}}
\frac
{|c_K(\eta,\eta^{i-})-c_K(\tilde\eta,\tilde\eta^{i-})|}{c_K(\eta,\eta
^{i-})}.\nonumber
\end{eqnarray}
We estimate
\begin{eqnarray*}
\label{HolleyEstimate}
&& -V(\omega,i^+)+\log\frac{\Gamma (\omega,i^+)}{\nu
(1_{\omega })}-\log\frac{\Gamma
(\omega ^i,i^-)}{\nu(1_{\omega ^i})}
\\[-3pt]
&&\qquad= \log\frac{\int\nu(d\eta)1_{\omega _{\Lambda }}(\eta)
((c_K(\omega _\Lambda \eta_{\Lambda ^c},\omega
_\Lambda ^i\eta_{\Lambda ^c}))/(c_K(\omega _{\Lambda }\zeta
_{\Lambda ^c},\omega ^{i}_{\Lambda }\zeta_{\Lambda ^c}))
}{\nu(1_{\omega _\Lambda })}
\\[-3pt]
&&\quad\qquad{} -\log\frac{\int\nu(d\eta)1_{\omega
^i_{\Lambda }}(\eta)
((c_K(\omega _\Lambda ^i\eta_{\Lambda ^c},\omega _\Lambda \eta
_{\Lambda ^c}))/(c_K(\omega ^i_{\Lambda }\zeta_{\Lambda ^c},\omega
_{\Lambda }\zeta_{\Lambda ^c}))}{\nu(1_{\omega ^i_\Lambda })}
\\[-3pt]
&&\qquad \leq \sup \biggl\{\log\frac{c_K(\eta_1,\eta^i_1)}{c_K(\eta_2,\eta
^i_2)}\dvtx  \eta _1=\eta_2\mbox{ on }\Lambda \biggr\}
\\[-3pt]
&&\quad\qquad{} +\sup \biggl\{\log\frac{c_K(\eta_1,\eta
^{i-}_1)}{c_K(\eta_2,\eta^{i-}_2)}\dvtx
\eta_1=\eta_2\mbox{ on }\Lambda \biggr\}.
\end{eqnarray*}
Using $\log a\leq a-1$ and expressing the oscillation on $\Lambda ^c$ via
single-point oscillations, we arrive at $\vartheta _\Lambda (i)$.
Similarly we get
for the second summand
\begin{eqnarray*}
\label{HolleyEstimate2}
&& -V \bigl(\omega,i^- \bigr)+\log\frac{\Gamma (\omega,i^-)}{\nu(1_\omega
)}-\log
\frac{\Gamma (\omega
^{i-},i^+)}{\nu(1_{\omega ^{i-}})}
\\[-3pt]
&&\qquad \leq \sup \biggl\{\log\frac{c_K(\eta_1,\eta^{i-}_1)}{c_K(\eta_2,\eta
^{i-}_2)}\dvtx  \eta_1=\eta_2
\mbox{ on }\Lambda \biggr\}
\\[-3pt]
&&\quad\qquad{} +\sup \biggl\{\log\frac
{c_K(\eta_1,\eta
^i_1)}{c_K(\eta_2,\eta^i_2)}\dvtx
\eta_1=\eta_2\mbox{ on }\Lambda \biggr\}
\\[-3pt]
&&\qquad\leq\vartheta _\Lambda (i).
\end{eqnarray*}
Hence
\[
\sum_{i\in\Lambda } \bigl[\kappa_\Lambda \bigl(i^+
\bigr)+\kappa_\Lambda \bigl(i^- \bigr) \bigr]\leq\sum
_{i\in\Lambda } \bigl[ \bigl(\beta _\Lambda \bigl(i^+ \bigr)+
\beta_\Lambda \bigl(i^- \bigr) \bigr)\vartheta _\Lambda (i)
\bigr].
\]
Notice that $\vartheta _\Lambda (i)\to0$ for all $i\in G$ as
$\Lambda \nearrow G$ since
our flip-rates are quasilocal and summable, indeed by the
well-definedness we have for all $i\in G$
\begin{eqnarray*}
\label{HolleyEstimate3}
\hspace*{-2pt}&& \sum_{j\notin\Lambda }\sup_{\eta_{j^c}=\tilde\eta_{j^c}}
\frac
{|c_K(\eta,\eta
^i)-c_K(\tilde\eta,\tilde\eta^i)|}{c_K(\eta,\eta^i)}
\\[-2pt]
\hspace*{-2pt}&&\qquad \leq \frac{2\pi e^{\|H_0\|}}{\min_{l\in\{1,\ldots,q\}
}\lambda(l)}\sup_{i\in G}\sum
_{j\in G}\sup_{\eta_{j^c}=\tilde\eta_{j^c}}\bigl|c_K \bigl(
\eta,\eta ^i \bigr)-c_K \bigl(\tilde\eta,\tilde\eta^i \bigr)\bigr|
=: A\times B^+<\infty.
\end{eqnarray*}
Notice also that $\kappa_{\Lambda _1}(i^\pm)\leq\kappa_{\Lambda
_2}(i^\pm)$ if
$i\in\Lambda _1\subset\Lambda _2$. Indeed if we look at the
subadditive function
$\varphi (x,y)=(x-y)\log\frac{x}{y}$ for $x,y>0$ and use
\[
\Gamma _{\Lambda _1} \bigl(\omega _1,i^\pm \bigr)=
\sum_{\omega _2=\omega
_1\ \mathrm{on}\ \Lambda _1}\Gamma _{\Lambda _2} \bigl(\omega
_2,i^\pm \bigr),
\]
we have
%
\begin{eqnarray}
\label{HolleyEstimate4} \kappa_{\Lambda _1} \bigl(i^\pm \bigr)&=&\sum
_{\omega _1\in\{1,\ldots,q\}
^{\Lambda _1}}\varphi \bigl[\Gamma _{\Lambda _1} \bigl(
\omega _1,i^\pm \bigr),\Gamma _{\Lambda _1} \bigl(
\omega _1^{i\pm},i^\mp \bigr) \bigr]
\nonumber\\[-8pt]\\[-8pt]
&\leq&\sum_{\omega _2\in\{1,\ldots,q\}^{\Lambda _2}}\varphi \bigl[\Gamma
_{\Lambda _2} \bigl(\omega _2,i^\pm \bigr),\Gamma
_{\Lambda
_2} \bigl(\omega _2^{i\pm},i^\mp
\bigr) \bigr]=\kappa_{\Lambda _2} \bigl(i^\pm \bigr).\nonumber
\end{eqnarray}
We are now in the position to finish the proof of Proposition~\ref
{TheGlauberDynamics}, part (2). This is a standard argument from \cite
{Li85} using translation-invariance and explicit control over boundary
terms, applied to the $q$-state model.

\begin{teo}\label{GlauberTheorem*}
Suppose that $G={\mathbb Z}^d$ and the Glauber dynamics flip-rates
%
\begin{equation}
\label{GlauberRatesSpecific} \frac{c_K(\omega ',(\omega ')^i)}{c_K((\omega ')^i,\omega
')}=\frac{\mu_{G\setminus i}[\omega
'_{G\setminus i}](\lambda ^{i} (e^{-H_{i}}
1_{(\omega '_i)^i}))}{\mu_{G\setminus i}[\omega '_{G\setminus
i}](\lambda ^{i} (e^{-H_{i}}
1_{\omega '_{i}}))} 
\end{equation}
are defined for a translation-invariant first-layer potential $H$. Then
a measure that is translation invariant and invariant w.r.t. $K$ must be
Gibbs for~$\gamma '$.
\end{teo}
\begin{pf}Let $\nu$ be invariant w.r.t. $K$ and
translation-invariant. Denote by
$\Lambda _n$ cubes in ${\mathbb Z}^d$ of side length $n$. Then we have
%
\begin{equation}
\label{HolleyEstimate5} \frac{1}{(kn)^d}\sum_{i\in\Lambda _{kn}} \bigl[
\kappa_{\Lambda
_{kn}} \bigl(i^+ \bigr)+\kappa_{\Lambda
_{kn}} \bigl(i^- \bigr)
\bigr] \geq\frac{1}{n^d}\sum_{i\in\Lambda _{n}} \bigl[\kappa
_{\Lambda
_{n}} \bigl(i^+ \bigr)+\kappa_{\Lambda _{n}} \bigl(i^- \bigr) \bigr].
\end{equation}
On the other hand $\beta_\Lambda (i^+)$ and $\beta_\Lambda (i^-)$
are uniformly
bounded and
%
\begin{eqnarray}
\label{HolleyEstimate6}
&& \frac{1}{n^d}\sum_{i\in\Lambda _{n}}
\vartheta _{\Lambda
_{n}}(i)\nonumber
\\
&&\qquad \leq\frac{A}{n^d}\sum
_{i\in\Lambda _{n}}\sum_{j\notin\Lambda _{n}} \bigl[
\delta_j \bigl(c_K \bigl(\cdot,\cdot^i
\bigr) \bigr)+\delta _j \bigl(c_K \bigl(\cdot,
\cdot^{i-} \bigr) \bigr) \bigr]
\nonumber\\[-8pt]\\[-8pt]
&&\qquad =\frac{A}{n^d}\sum_{i\in\Lambda _{n}}\sum
_{j\notin\Lambda
_{n}} \bigl[\delta _{j-i} \bigl(c_K
\bigl(\cdot,\cdot^0 \bigr) \bigr)+\delta_{j-i}
\bigl(c_K \bigl(\cdot,\cdot^{0-} \bigr) \bigr) \bigr]\nonumber
\\
&&\qquad = A\sum_{l\in{\mathbb Z}^d} \bigl[\delta_{l}
\bigl(c_K \bigl(\cdot,\cdot^0 \bigr) \bigr)+\delta
_{l} \bigl(c_K \bigl(\cdot,\cdot^{0-} \bigr)
\bigr) \bigr]\frac{\#\{i\in\Lambda _n\dvtx i+l\notin\Lambda _n\}}{n^d}.\nonumber
\end{eqnarray}
This tends to zero since the oscillations are bounded by $B^\pm$ and
the fact that an increasing strip of boundary of cubes goes to infinity
slower than the volume. Together we can write
%
\begin{eqnarray}
\label{HolleyEstimate7} \frac{1}{n^d}\sum_{i\in\Lambda _{n}} \bigl[
\kappa_{\Lambda
_{n}} \bigl(i^+ \bigr)+\kappa_{\Lambda
_{n}} \bigl(i^- \bigr)
\bigr]&\leq&\frac{1}{n^d}\sum_{i\in\Lambda _n} \bigl[
\bigl(\beta _{\Lambda _n} \bigl(i^+ \bigr)+\beta _{\Lambda _n} \bigl(i^-
\bigr) \bigr) \vartheta _{\Lambda _n}(i) \bigr]
\nonumber\\[-8pt]\\[-8pt]
&\leq& C\frac{1}{n^d}\sum_{i\in\Lambda _{n}}\vartheta
_{\Lambda
_{n}}(i)\to0\qquad\mbox{for }n\to\infty\nonumber
\end{eqnarray}
and hence by the nonnegativity of $\kappa_\Lambda (i^{\pm})$ we
have $\kappa
_{\Lambda _n}(i^{+})=\kappa_{\Lambda _n}(i^{-})=0$ for all $i\in
\Lambda _n$. By the
subadditivity argument $\kappa_\Lambda (i^{+})=\kappa_\Lambda
(i^{-})=0$ for all $\Lambda
$ and $i\in\Lambda $.
Thus for all finite $\Lambda $, $\omega \in\{1,\ldots,q\}^\Lambda
$ and $i\in\Lambda $
\begin{eqnarray*}
\label{InEquOutflow} 0&=&\Gamma \bigl(\omega,i^+ \bigr)-\Gamma \bigl(\omega
^i,i^- \bigr)
\\
&=&\int\nu(d\eta) \bigl[\nu(\omega _i|\omega _{\Lambda \setminus
i}
\eta_{\Lambda ^c})c_K \bigl(\omega _\Lambda
\eta_{\Lambda
^c},\omega ^i_\Lambda \eta_{\Lambda ^c}
\bigr)
\\
&&\hspace*{39pt}{}-\nu \bigl(\omega _i^i|\omega _{\Lambda \setminus i}
\eta_{\Lambda ^c} \bigr)c_K \bigl(\omega ^i_\Lambda
\eta_{\Lambda ^c},\omega _\Lambda \eta_{\Lambda ^c} \bigr) \bigr].
\end{eqnarray*}
So $\nu$-a.s. we have
\[
\label{InEquOutflow1} \frac{\nu(\omega _i^i|\omega _{\Lambda \setminus i}\eta
_{\Lambda ^c})}{\nu(\omega _i|\omega _{\Lambda \setminus i}\eta
_{\Lambda
^c})}=\frac{c_K(\omega _\Lambda \eta_{\Lambda ^c},\omega
^i_\Lambda \eta_{\Lambda ^c})}{c_K(\omega ^i_\Lambda \eta
_{\Lambda ^c},\omega _\Lambda \eta_{\Lambda ^c})}=\frac{\mu
_{G\setminus i}[\eta_{\Lambda ^c}\omega _{\Lambda
\setminus i}](\lambda(e^{-H_i}1_{\omega ^i_i}))}{\mu_{G\setminus
i}[\eta
_{\Lambda ^c}\omega _{\Lambda \setminus i}](\lambda
(e^{-H_i}1_{\omega _i}))}.
\]
Since we compare discrete measures on sites $i\in G$, it follows by the
remark below, $\nu(\omega _i|\omega _{i^c})=\gamma '_i(\omega
_i|\omega _{i^c})$ $\nu$-almost
everywhere and thus $\nu\in\mathcal{G}(\gamma ')$.
\end{pf}

\begin{rem}
Let $(a_1, \ldots, a_q)$ and $(b_1, \ldots, b_q)$ be probability
vectors with
$\frac{a_k}{a_{k+1}}=\frac{b_k}{b_{k+1}}$ for all $k\in\{1, \ldots,q\}$.
Then we have $\frac{a_l}{a_{k}}=\frac{b_l}{b_{k}}$ for all $k,l\in\{1,\ldots,q\}$ and thus
%
\begin{eqnarray}
\label{GlauberReversibleisGibbs7} a_l&=&\frac{a_l}{\sum_{k=1}^q a_k}=\frac{1}{1+\sum_{k\neq l}
({a_k}/{a_l})}=
\frac{1}{1+\sum_{k\neq l} ({b_k}/{b_l})}=b_l.
\end{eqnarray}
\end{rem}

%
%
%
%
%
%
%
%


\section{Joint dynamics}\label{sec5}\label{5}

Let us now consider the joint dynamics $L+\alpha K$ for $\alpha >0$.
Of course
well-definedness [Proposition~\ref{TheJointDynamics}, part (1)]
follows directly from the fact, that the individual rates of $L$ and
$K$ are well defined.
%
%

As a warning, we note that the generators $L$ and $K$ do not commute
(except in the limit $q\to\infty$). To see this we apply $LK-KL$
%
%
to the local observable $\psi:=1_{\eta_{\Lambda }}$ for a finite
$\Lambda \subset
G$. Evaluated,
%
%
for instance, at $\omega _{\Lambda }=\eta_{\Lambda }$, we find the expression
%
\begin{eqnarray}
\label{Commutator2}
&& \sum_{i\in\Lambda } \bigl(c_L
\bigl(\omega,\omega ^i \bigr)c_K \bigl(\omega
^i,\omega \bigr)-c_K \bigl(\omega,\omega
^{i-} \bigr)c_L \bigl(\omega ^{i-},\omega \bigr)
\bigr)
\nonumber\\[-8pt]\\[-8pt]
&&\qquad =\sum_{i\in\Lambda } \bigl(
\mu_{G\setminus i}[\omega _{G\setminus
i}] \bigl(e^{-H_{i}((\omega _i)^l,\cdot_{i^c})} \bigr)-
\mu_{G\setminus i}[\omega _{G\setminus i}] \bigl(e^{-H_{i}((\omega
_i^{i-})^l,\cdot_{i^c})} \bigr) \bigr).\nonumber
\end{eqnarray}
This does not vanish in general, and thus
the commutator is not zero. But if we consider the limit of the
coarse-graining, that is, letting $q$ the number of discrete states go
to infinity, we approach a commutative setting. This result reflects
the continuum situation in the Maes and Shlosman program~\cite{MaSh11}.

As a consequence, $S^{L+\alpha K}_t\neq S^L_tS^{\alpha K}_t$, and it
is not
immediate that the joint dynamics also rotates the discrete Gibbs
measures in the sense of Proposition~\ref{TheJointDynamics}, part (2).
To see that this is nevertheless true one has to follow the same
arguments as in Section~\ref{Rotationpropertyofthegenerator} and
notice $\|\!|\Gamma ^{\mathrm{joint}}\|\!|_{e^\varrho }<\infty$.

%
%
%
%
%
%

\subsection{The invariant measure for the joint dynamics}\label{sec5.1}
In this subsection we show Proposition~\ref{TheJointDynamics}, part
(3) and Corollary~\ref{TheUniqueInvariantMeasureForJoint}. First
let us verify that indeed the symmetrically mixed measure is invariant
and in the set of Gibbs measures this is the only one. Finally we prove
that measures that are invariant under the joint dynamics must be Gibbs.

The mixture of all translation-invariant extremal Gibbs measures $\mu
'_t$
\[
\mu'_*:=\frac{1}{2\pi}\int_0^{2\pi}
\mu'_t\,dt
\]
is invariant for the rotation dynamics and hence for the joint dynamics
$L+\alpha K$. Indeed, let $(S^L_t)_{t\geq0}$ be the semigroup for $L$ and
$f$ a quasilocal observable, we have
%
\begin{eqnarray}
\label{MixedMeasure} \int S^L_tf(\eta)\mu'_*(d
\eta)&=&\frac{1}{2\pi}\int_0^{2\pi}\int
S^L_tf(\eta)\mu'_s(d
\eta)\,ds\nonumber
\\
&=&\int f(\eta)\frac{1}{2\pi}\int_0^{2\pi}
\mu'_{s+t}(d\eta )\,ds
\\
&=& \int f(\eta )\mu'_*(d
\eta).\nonumber
\end{eqnarray}

\begin{prop}
There are no translation-invariant invariant Gibbs measures for the
rotation dynamics other then $\mu'_*$.
\end{prop}
\begin{pf}We know from Theorem 7.26 of \cite{Ge11} that every
Gibbs measure $\mu'\in\mathcal{G}_{\theta}(\gamma ')$ has a unique
representation
\[
\mu'=\int_{\operatorname{ex}\mathcal{G}_\theta(\gamma ')}\bar {\mu} w_{\mu
'}(d
\bar {\mu}),
\]
where $w_{\mu'}\in\mathcal{P}(\operatorname{ex}\mathcal{G}_\theta
(\gamma
'),\sigma (\operatorname{ex}\mathcal{G}_\theta(\gamma ')))$, and
$\sigma (\mathcal
{P})$ is the so-called
evaluation \mbox{$\sigma $-}algebra. Since the Gibbs measures can be labeled as
described above, there is a bijection
\[
b\dvtx  \operatorname{ex}\mathcal{G}_{\theta} \bigl(\gamma '
\bigr)\to[0,2\pi )=S^1,\qquad \mu '\mapsto\arg
\bigl(e_{\hat{m}'_0} \bigl(\mu' \bigr)/m_\beta \bigr),
\]
where $b$ is $(\sigma (\operatorname{ex}\mathcal{G}_{\theta
}(\gamma ')),\mathcal
{B}([0,2\pi)))$ measurable and $\arg$ denotes the argument of a number
in $S^1$. Indeed since $\hat{m}'_0$ is bounded and measurable, so is
$e_{\hat{m}'_0}\dvtx \mu'\mapsto\mu'(\hat{m}'_0)\in{\mathbb R}^2$, and
thus $b(\mu
')=\arg(e_{\hat{m}'_0}(\mu')/m_\beta)$ is a composition of measurable
functions. Hence we can consider image measures $v_{\mu'}$ of $w_{\mu
'}$ under $b$.

On the other hand for all local coarse-grained sets $A'\in\mathcal
{F}'$, the mapping
\[
c_{A'}\dvtx  [0,2\pi)\to[0,1],\qquad   t\mapsto\mu'_t
\bigl(A' \bigr)=\mu _t(A)=\lim_{\Lambda \nearrow{\mathbb Z}^d}
\gamma _\Lambda (A|\omega _t),
\]
where $A:=T^{-1}(A')$ and $\omega _t$ the homogeneous boundary
condition as
described in the \hyperref[sec1]{Introduction},
is Borel-measurable as a composition of measurable maps. We also used
the measurability of $t\mapsto\omega _t$. Hence this is true for all
$A'\in
\mathcal{F'}$.

By the transformation theorem for measurable maps we have for all
\mbox{$A'\in
\mathcal{F'}$}
%
\begin{eqnarray}
\label{MixedMeasureOnlyInvGibbs} \mu' \bigl(A' \bigr)&=&\int
_{\operatorname{ex}\mathcal{G}_\theta(\gamma
')}\bar {\mu} \bigl(A' \bigr)
w_{\mu
'}(d\bar {\mu})=\int_{\operatorname{ex}\mathcal{G}_\theta(\gamma
')}c_{A'}
\bigl(b(\bar{\mu }) \bigr) w_{\mu'}(d\bar {\mu})\nonumber
\\
&=&\int_0^{2\pi}c_{A'}(t)
w_{\mu'} \bigl(b^{-1}(dt) \bigr)=\int_0^{2\pi
}c_{A'}(t)v_{\mu'}(dt)
\\
&=& \int_0^{2\pi}\mu'_t
\bigl(A' \bigr)v_{\mu'}(dt).\nonumber
\end{eqnarray}

By looking at tail-measurable interval sets
\[
A_{[0,u)}:= \biggl\{\omega \in{\Omega}\dvtx \lim_{n\to\infty}
\frac
{1}{|\Lambda _n|}\sum_{j\in\Lambda
_n}\arg \biggl(
\frac{\omega _j}{m_\beta} \biggr)=[0,u) \biggr\}
\]
and $\varphi _{[0,u)}(\omega '):=\mu_G[\omega '](A_{[0,u)})$ we see
that $v_{\mu'}$
has to be a translation-invariant Borel-measure, indeed
%
\begin{eqnarray}
\label{MixedMeasureUnique} v_{\mu'} \bigl([0,u )\bigr)&=&\int_0^{2\pi}
\mu'_{t}(\varphi _{[0,u)})v_{\mu
'}(dt)=
\mu '(\varphi _{[0,u)})=\mu'S^L_s(
\varphi _{[0,u)})
\nonumber\\[-8pt]\\[-8pt]
&=&\int_0^{2\pi}\mu'_{t+s}(
\varphi _{[0,u)})v_{\mu'}(dt)=v_{\mu
'} \bigl([-s,u-s)\bigr)\nonumber
\end{eqnarray}
for all $s\in[0,2\pi)$. Since $\{[0,u)\dvtx u\in[0,2\pi)\}$ is a generator
for the Borel-$\sigma $-algebra, and $v_{\mu'}$ is a probability
measure, we
have $v_{\mu'}(dt)=\frac{1}{2\pi}\lambda (dt)$. \end{pf}

Since $S^{L+\alpha K}_s(\mu'_t)=\mu'_{t+s}=S^{L}_s(\mu'_t)$ we can conclude
$\mu'_*$ is the only translation-invariant measure that is also
invariant w.r.t. the joint dynamics.
The next proposition proves Proposition~\ref{TheJointDynamics},
part (3).

\begin{prop}\label{InvarianceIsGibbs}
Every translation-invariant measure that is invariant for the joint
dynamics $L+\alpha K$ with $\alpha >0$ is a Gibbs measure.
\end{prop}
\begin{pf}Let $\Lambda \subset{\mathbb Z}^d$ be a finite set, and
$\zeta\in{\Omega}'$
be an arbitrary but fixed boundary condition for the second-layer
specification; that is, consider the coarse-grained measure $\gamma
_\Lambda '(\omega
|\zeta)=\frac{\mu_{\Lambda ^c}[\zeta_{\Lambda ^c}](\lambda
^\Lambda (e^{-H_\Lambda }1_{\omega _\Lambda
}))}{\mu_{\Lambda ^c}[\zeta_{\Lambda ^c}](\lambda^\Lambda
(e^{-H_\Lambda }))}$ on $\{1, \ldots, q\}
^\Lambda $.
Our strategy for the proof is to again look at the derivative of the
local relative entropy $H_\Lambda (\nu|\gamma _\Lambda '(\cdot
|\zeta))$
for $\nu$ translation-invariant and invariant w.r.t. the joint dynamics.
We have seen in case of the Glauber dynamics how to verify Gibbsianness
for invariant measures by estimating certain terms in the derivative of
the local relative entropy. Those term are only of the order of the
boundary $|\partial\Lambda |$.
This allowed us to prove the DLR equality for the invariant measure. A
crucial ingredient is the translation-invariance of both, the model as
well as the invariant measure.

Essentially we follow the same line of arguments here, taking special
care of the contribution of the rotation. We look at an approximating
local open boundary rotation dynamics and show its relative entropy is
decreasing. This means that the approximating rotation only ``helps''
the Glauber dynamics argument. The error we make by using the
approximation instead of the infinite-volume rotation dynamics is only
of boundary order and thus again increases more slowly than the volume.

Since the time-derivative of the local relative entropy is additive as
a sum of the two terms corresponding to the two generators $K$ and $L$,
we can calculate separately for the Glauber and for the rotation
dynamics. We write $\nu_{t,L}$ (resp., $\nu_{t,K}$, $\nu_{t,L+\alpha K}$)
for the measure $\nu$ propagated only by the rotation (resp., by the
Glauber dynamics, by the joint dynamics).

Let us compute for the rotation $\frac{d}{dt}_{|_{t=0}}H_\Lambda (\nu
_{t,L}|\gamma
_\Lambda '(\cdot|{\zeta}))$ with $\nu=\nu_0$. Again we do this in
two steps.
Similarly to the computations done in (\ref{RelativeEntropy3}) we find
%
\begin{eqnarray}
\label{JointRelativeEntropy3} \qquad \frac{d}{dt}_{|_{t=0}}\sum
_{\omega \in\{1,\ldots,q\}^\Lambda }\nu _{t,L}(1_\omega )\log\nu
_{t,L}(1_\omega ) 
&=&\sum_{\omega, i\in\Lambda }
\Gamma _L \bigl(\omega,i^+ \bigr)\log\frac
{\nu(1_{\omega ^i})}{\nu(1_{\omega })},
\end{eqnarray}
where\vspace*{1pt} we
again wrote $\Gamma _L(\omega,i^+):=\int\nu(d\eta)c_L(\eta,\eta
^i)1_{\omega }(\eta)$
for the outflows of $1_\omega $ in the direction $i^+$.
For the other summand of $\frac{d}{dt}_{|_{t=0}}H_\Lambda (\nu
_{t,L}|\gamma _\Lambda
'(\cdot|{\zeta}))$ we have
%
\begin{eqnarray}
\label{JointRelativeEntropySecondPart1} && \frac{d}{dt}_{|_{t=0}}\int\nu_{t,L}(d
\omega )\log\mu_{\Lambda
^c}[\zeta_{\Lambda
^c}] \bigl(\lambda^\Lambda
\bigl(e^{-H_\Lambda }1_{\omega _\Lambda } \bigr) \bigr) 
\nonumber\\[-8pt]\\[-8pt]
&&\qquad =\sum
_{\omega, i\in\Lambda }V^\zeta \bigl(\omega,i^+ \bigr)\Gamma
_L \bigl(\omega,i^+ \bigr),\nonumber
\end{eqnarray}
where we again defined $V^\zeta(\omega,i^+):=\log\frac{\mu
_{\Lambda ^c}[\zeta_{\Lambda
^c}](\lambda^\Lambda (e^{-H_\Lambda }1_{\omega ^i_\Lambda }))}{\mu
_{\Lambda ^c}[\zeta_{\Lambda ^c}](\lambda
^\Lambda (e^{-H_\Lambda }1_{\omega _\Lambda }))}$.
Together we have
%
\begin{eqnarray}
\label{JointRelativeEntropyTogether}\qquad \frac{d}{dt}_{|_{t=0}}H_\Lambda \bigl(
\nu_{t,L}|\gamma _\Lambda '(\cdot |{\zeta})
\bigr)&=&\sum_{\omega,
i\in\Lambda }\Gamma _L \bigl(\omega,i^+ \bigr) \biggl(\log\frac{\nu
(1_{\omega ^i})}{\nu(1_{\omega })}-V^\zeta \bigl(\omega,i^+
\bigr) \biggr).
\end{eqnarray}
We define the approximating local generator $\widetilde L_\Lambda $ via the
following open boundary rates:
\[
c_{\widetilde L_\Lambda } \bigl(\eta,\eta^i \bigr):= %
\cases{
\displaystyle\frac{\lambda^{\Lambda \setminus i}(e^{-\widetilde
H_\Lambda (\eta
_i|^r,\cdot)}1_{\eta
_{\Lambda \setminus i}})}{\lambda^\Lambda (e^{-\widetilde H_\Lambda
}1_{\eta_\Lambda })}, &\quad if $i\in\Lambda $,
\vspace*{5pt}\cr
0, &\quad if $i\in
\Lambda ^c$,}
\]
where\vspace*{2pt} $\Lambda $ is a fixed finite volume, and $\widetilde H_\Lambda:=\sum_{A\subset\Lambda
}\Phi_A$ is the open boundary Hamiltonian for $\Lambda $ in the first-layer
model. Let $(S^{\widetilde L_{\Lambda }}_{t})_{t\geq0}$ be the
corresponding semigroup. Since we assume the underlying first-layer
potential to be rotation-invariant, the open boundary measure
$\tilde\gamma _\Lambda (\omega _\Lambda ):=\lambda^\Lambda
(e^{-\widetilde H_\Lambda }1_{\omega _\Lambda })/\break \lambda^\Lambda
(e^{-\widetilde H_\Lambda })$ on $\{1, \ldots,q\}^\Lambda $ is invariant
for $\widetilde L_\Lambda
$. Indeed for all $\omega _\Lambda \in\{1, \ldots, q\}^\Lambda $
we have
%
\begin{eqnarray}
\label{LocalInvariantMeasures} \qquad\tilde\gamma _\Lambda \bigl(\widetilde L_\Lambda
(1_{\omega _\Lambda }) \bigr) 
&=&\frac{\sum_{i\in\Lambda }[\lambda^{\Lambda \setminus
i}(e^{-\widetilde H_\Lambda (\omega
_i|^l,\cdot)}1_{\omega _{\Lambda \setminus i}})-\lambda^{\Lambda
\setminus i}(e^{-\widetilde
H_\Lambda (\omega _i|^r,\cdot)}1_{\omega _{\Lambda \setminus
i}})]}{\lambda^\Lambda (e^{-\widetilde H_\Lambda
})}
\nonumber\\[-8pt]\\[-8pt]
&=&\frac{d}{d\varepsilon }_{|_{\varepsilon =0}}\tilde\gamma _\Lambda
(1_{\omega _\Lambda +\varepsilon })=\frac{d}{d\varepsilon
}_{|_{\varepsilon =0}}\tilde\gamma
_\Lambda (1_{\omega _\Lambda })=0.\nonumber
\end{eqnarray}
We can employ a standard argument for the decrease of relative entropy
in finite volume in order to determine the sign of $\frac
{d}{dt}_{|_{t=0}}H_{\Lambda }(\nu_{t,\widetilde L_\Lambda }|\tilde\gamma _\Lambda )$. Indeed if we
use the convex function $\psi(x)=x\log x+x-1$, the relative entropy reads
%
\begin{eqnarray}
\label{RelativeEntropyApproxJensen} H_{\Lambda }(\nu_{t,\widetilde L_\Lambda }|\tilde\gamma
_\Lambda )&=&\sum_{\omega }\tilde\gamma
_\Lambda (1_\omega )\psi \biggl(\frac{\nu_{t,\widetilde L_\Lambda }(1_{\omega })}{\tilde\gamma
_\Lambda (1_\omega )} \biggr)
\nonumber\\[-8pt]\\[-8pt]
&=&\sum_{\omega }
\tilde\gamma _\Lambda (1_\omega )\psi \biggl(\frac
{1}{\tilde\gamma _\Lambda (1_\omega )}
\sum_{\eta
} S^{\widetilde L_{\Lambda }}_{t}(1_{\omega })
(\eta)\frac{\nu(\eta
)}{\tilde\gamma _\Lambda (\eta
)}\tilde\gamma _\Lambda (\eta) \biggr),\nonumber
\end{eqnarray}
where $\frac{S^{\widetilde L_{\Lambda }}_{t}(1_{\omega })(\eta)}{\tilde\gamma (1_\omega
)}\tilde\gamma (d\eta)=\frac{1_{\omega }(\eta)}{\tilde\gamma
(1_\omega )}\tilde\gamma (d\eta
)$ is a probability measure. Hence we can use Jensen's inequality and obtain
%
\begin{eqnarray}
\label{RelativeEntropyApproxJensen2} H_{\Lambda }(\nu_{t,\widetilde L_\Lambda }|\tilde\gamma
_\Lambda )&\leq&\sum_{\omega }\tilde\gamma
(1_\omega )\frac{1}{\tilde\gamma (1_\omega )}\sum_{\eta}
S^{\widetilde
L_{\Lambda }}_{t}(1_{\omega }) (\eta )\psi \biggl(
\frac{\nu(\eta)}{\tilde\gamma _\Lambda (\eta)} \biggr)\tilde\gamma _\Lambda (\eta)\nonumber
\\
&=&\sum_{\omega }\psi \biggl(\frac{\nu(\omega )}{\tilde\gamma _\Lambda
(\omega )} \biggr)
\tilde\gamma _\Lambda (\omega )
\\
&=& H_{\Lambda}(\nu|\tilde\gamma
_\Lambda )\nonumber
\end{eqnarray}
with equality if and only if $\nu_{t,\widetilde L_\Lambda }=\tilde\gamma _\Lambda $. Thus
the derivative must be nonpositive
%
\begin{eqnarray}
\label{RelativeEntropyLocalDerivative}
0&\geq&\frac{d}{dt}_{|_{t=0}}H_{\Lambda }(
\nu_{t,\widetilde L_\Lambda
}|\tilde\gamma _\Lambda )\nonumber
\\
&=&\sum
_{\omega, i\in\Lambda }\nu_{t}(1_{\omega })c_{\widetilde
L_\Lambda }
\bigl(\omega,\omega ^i \bigr) \biggl(\log\frac
{\nu(1_{\omega ^i})}{\nu(1_{\omega })}-\log
\frac{\tilde\gamma
_\Lambda (\omega ^i)}{\tilde\gamma _\Lambda
(\omega )} \biggr)
\\
&=:&\sum_{\omega, i\in\Lambda }\Gamma _{\widetilde L_\Lambda
} \bigl(\omega,i^+ \bigr) \biggl(\log\frac{\nu(1_{\omega
^i})}{\nu(1_{\omega })}-V_{\widetilde L_\Lambda } \bigl(\omega,i^+
\bigr) \biggr).\nonumber
\end{eqnarray}

We are going to show $|\frac{d}{dt}_{|_{t=0}}H_{\Lambda }(\nu
_{t,L}|\gamma _\Lambda
'(\cdot|{\zeta}))-\frac{d}{dt}_{|_{t=0}}H_{\Lambda }(\nu_{t,\widetilde
L_\Lambda
}|\tilde\gamma _\Lambda )|=o(|\Lambda |)$. Let us start with the
following estimate:
%
\begin{eqnarray}
\label{RelativeEntropyLocalError}
&& \frac{d}{dt}_{|_{t=0}}H_{\Lambda }\bigl(
\nu_{t,L}|\gamma _\Lambda '(\cdot|{\zeta})\bigr)-
\frac{d}{dt}_{|_{t=0}}H_{\Lambda }(\nu_{t,\widetilde L_\Lambda }|\tilde\gamma _\Lambda )\nonumber
\\
&&\qquad =\sum_{\omega, i\in\Lambda } \bigl[\Gamma
_L \bigl(\omega,i^+ \bigr)-\Gamma _{\widetilde L_\Lambda } \bigl(\omega,i^+ \bigr) \bigr]\log\frac{\nu
(1_{\omega ^i})}{\nu(1_{\omega })}\nonumber
\\
&&\quad\qquad{} +\sum_{\omega, i\in\Lambda } \bigl[
\bigl[V_{\widetilde
L_\Lambda } \bigl(\omega,i^+ \bigr)-V^\zeta \bigl(\omega,i^+ \bigr) \bigr]\Gamma _{L} \bigl(\omega,i^+ \bigr)
\nonumber\\[-8pt]\\[-8pt]
&&\hspace*{39pt}\quad\qquad{}- \bigl[\Gamma _L \bigl(\omega,i^+ \bigr)-\Gamma
_{\widetilde L_\Lambda
} \bigl(\omega,i^+ \bigr) \bigr]V_{\widetilde L_\Lambda } \bigl(\omega,i^+ \bigr) \bigr]\nonumber
\\
&&\qquad \leq \sum_{\omega, i\in\Lambda }A \bigl(\omega,i^+ \bigr)
\nu(1_\omega )\biggl|\log\frac{\nu(1_{\omega ^i})}{\nu
(1_{\omega })}\biggr|\nonumber
\\
&&\quad\qquad{} +\sum_{\omega, i\in\Lambda } \bigl[B \bigl(\omega,i^+
\bigr)\Gamma _{L} \bigl(\omega,i^+ \bigr)-A \bigl(\omega,i^+ \bigr)
\nu(1_\omega )\bigl|V_{\widetilde L_\Lambda } \bigl(\omega,i^+ \bigr)\bigr| \bigr],\nonumber
\end{eqnarray}
where we defined $B(\omega,i^+):=|V_{\widetilde L_\Lambda }(\omega,i^+)-V^\zeta(\omega,i^+)|$
and used the following estimate and definition:
\begin{eqnarray*}
\label{RelativeEntropyLocalErrorOrder1}
&& \Gamma _L \bigl(\omega,i^+ \bigr)-\widetilde\Gamma
_{\widetilde L_\Lambda } \bigl(\omega,i^+ \bigr)
\\
&&\qquad =\int\nu(d\eta)c_L \bigl(
\eta, \eta ^i \bigr)1_{\omega }(\eta)-\nu(1_{\omega })c_{\widetilde L_\Lambda }
\bigl(\omega,\omega ^i \bigr)
\\
&&\qquad =\int\nu(d\eta_{\Lambda ^c})\nu(1_{\omega }|\eta_{\Lambda
^c})
\bigl[c_L \bigl(\omega _\Lambda \eta_{\Lambda ^c},\omega
^i_\Lambda \eta_{\Lambda ^c} \bigr)-c_{\widetilde L_\Lambda }
\bigl( \omega,\omega ^i \bigr) \bigr]
\\
&&\qquad \leq \sup_{\eta_{\Lambda ^c}}\biggl|\frac{\mu_{G\setminus i}[\eta
_{\Lambda ^c}\omega _{\Lambda
\setminus i}](e^{-H_{i}(\omega _i^r,\cdot_{i^c})})
}{\mu_{G\setminus i}[\eta_{\Lambda ^c}\omega _{\Lambda \setminus
i}](\lambda (e^{-H_{i}(\cdot
)}1_{\omega _{i}}))}-\frac{\lambda^{\Lambda \setminus i}(e^{-\widetilde
H_\Lambda (\omega
_i^r,\cdot_{\Lambda \setminus i})}1_{\omega _{\Lambda \setminus
i}})}{\lambda^\Lambda
(e^{-\widetilde H_\Lambda (\cdot_\Lambda )}1_{\omega _\Lambda })}\biggr|\nu
(1_{\omega })
\\
&&\qquad =: A \bigl(\omega,i^+ \bigr)\nu(1_\omega ).
\end{eqnarray*}

We first verify $\sup_{\omega }\sum_{i\in\Lambda }A(\omega,i^+)=o(|\Lambda |)$ and $\sup_{\omega
}\sum_{i\in\Lambda }B(\omega,i^+)=o(|\Lambda |)$. We do this in
the two following lemmata.
%
\begin{lem}\label{EstimateFirstErrorTerm}
$\sup_{\omega }\sum_{i\in\Lambda }A(\omega,i^+)=o(|\Lambda |)$.
\end{lem}
\begin{pf}In order to see cancellations we define for a given
second-layer boundary condition inside $\Lambda $, namely $\omega
_{\Lambda }$, and open
boundary conditions outside~$\Lambda $, the conditional first-layer
probability measures on $(S^1)^{G\setminus i}$
%
\begin{equation}
\label{RelativeEntropyLocalErrorOrder2} \tilde\mu_{G\setminus i}[\omega _{\Lambda \setminus i}](\varphi ):=
\frac{\lambda
^{{i^c}}(\varphi e^{-\sum_{i\notin A\subset\Lambda } \Phi
_A}1_{\omega _{\Lambda \setminus
i}})}{\lambda^{{i^c}}(e^{-\sum_{i\notin A\subset\Lambda } \Phi
_A}1_{\omega _{\Lambda
\setminus i}})}.
\end{equation}
In particular
$\frac{\tilde\mu_{G\setminus i}[\omega _{\Lambda \setminus
i}](e^{-\widetilde H_i(\omega
_i|^r,\cdot_{\Lambda \setminus i})})}{\tilde\mu_{G\setminus
i}[\omega _{\Lambda
\setminus i}](\lambda (e^{-\widetilde H_i(\cdot_{\Lambda })}1_{\omega
_i}))}=\frac{\lambda
^{\Lambda \setminus i}(e^{-\widetilde H_\Lambda (\omega _i|^r,\cdot
_{\Lambda \setminus i})}1_{\omega
_{\Lambda \setminus i}})}{\lambda^\Lambda (e^{-\widetilde H_\Lambda
(\cdot_\Lambda )}1_{\omega _\Lambda })}$.
These fractions again give rise to a specification $\tilde\gamma $ on the
second layer when we look at subvolumes, keeping the $\Lambda $ fixed.

In essence we want to exploit the Dobrushin comparison theorem.
Since we can bound every term by some constant times $e^{\pm\|H_i\|}=e^{\pm\|H_i\|}:=e^K$ it suffices to estimate the distance of the
conditional first-layer Gibbs measures $\mu_{G\setminus i}[\eta
_{\Lambda
^c}\omega _{\Lambda \setminus i}]$ and $\tilde\mu_{G\setminus
i}[\omega _{\Lambda \setminus
i}]$ applied to the quasilocal functions
%
\begin{eqnarray}
\label{QuasilocalFunctions1} \psi_1^{\omega _i}(\cdot)&:=&e^{-H_{i}(\omega _i|^r,\cdot
_{i^c})},\qquad
 \psi _2^{\omega _i}(\cdot):=\lambda \bigl(
e^{-H_{i}(\cdot,\cdot
_{i^c})}1_{\omega _{i}} \bigr)\quad\mbox{and}
\nonumber\\[-8pt]\\[-8pt]
\tilde\psi_1^{\omega _i}(\cdot)&:=&e^{-\widetilde H_{i}(\omega
_i|^r,\cdot
_{i^c})},\qquad
\tilde\psi_2^{\omega _i}(\cdot):=\lambda \bigl(
e^{-\widetilde
H_{i}(\cdot,\cdot_{i^c})}1_{\omega _{i}} \bigr).\nonumber
\end{eqnarray}
(Notice that we have done computations of the same flavor in the
section about the well-definedness of the rotation dynamics.)
For any fixed first-layer boundary condition $w\in{\Omega}$ the
measure $\mu
_{G\setminus i}[\eta_{\Lambda ^c}\omega _{\Lambda \setminus i}]$
is uniquely admitted
by the specification
%
\begin{equation}
\label{RelativeEntropyLocalErrorOrder4} \gamma ^{\eta_{\Lambda ^c}\omega _{\Lambda \setminus
i}}|_{i^c}:= \bigl( \bigl(\gamma
^{\eta_{\Lambda ^c}\omega _{\Lambda
\setminus i}}|_{i^c} \bigr)_{\Delta }(\cdot|w_{\Delta ^c\setminus
i})
\bigr)_{\Delta \subset
i^c }
\end{equation}
and $\tilde\mu_{G\setminus i}[\omega _{\Lambda \setminus i}]$ is
admitted by
$\tilde\gamma ^{\omega _{\Lambda \setminus i}}|_{i^c}:=
((\tilde\gamma ^{\omega _{\Lambda \setminus
i}}|_{i^c})_{\Delta }(\cdot|w_{\Delta ^c\setminus i})
)_{\Delta \subset i^c }$,
$\Delta $ being finite subsets of $i^c$.
The total variational distance between the two specifications on the
site $l\neq i$ can be estimated by
\begin{eqnarray*}
\label{RelativeEntropyLocalErrorOrder5} b_l&:=&\sup_{\eta_{\Lambda ^c}\omega _{\Lambda \setminus
i},w_{l^c\setminus i}}\bigl\| \bigl(\gamma ^{\eta_{\Lambda ^c}\omega _{\Lambda \setminus i}}|_{i^c} \bigr)_{l}(\cdot
|w_{l^c\setminus
i})- \bigl(\tilde\gamma ^{\omega _{\Lambda \setminus
i}}|_{i^c}
\bigr)_{l}(\cdot|w_{l^c\setminus
i})\bigr\|_l
\\
&=&\sup_{\eta_{\Lambda ^c}\omega _{\Lambda \setminus
i},w_{l^c\setminus i}, B\in\mathcal
{S}^1}\biggl|\frac{\lambda (1_B1_{\eta_{\Lambda ^c}\omega _{\Lambda
\setminus i}}e^{-\sum_{i\notin
A\ni l}\Phi_A(\cdot,w_{l^c\setminus i})})}{\lambda (1_{\eta
_{\Lambda ^c}\omega _{\Lambda
\setminus i}}e^{-\sum_{i\notin A\ni l}\Phi_A(\cdot,w_{l^c\setminus
i})})}
\\
&&\hspace*{77pt}{}-\frac{\lambda (1_B1_{\omega _{\Lambda \setminus
i}}e^{-\sum_{i\notin A\ni
l, A\subset\Lambda }\Phi_A(\cdot,w_{l^c\setminus i})})}{\lambda
(1_{\omega _{\Lambda \setminus
i}}e^{- \sum_{i\notin A\ni l, A\subset\Lambda }\Phi_A(\cdot,w_{l^c\setminus
i})})}\biggr|
\\
&\leq& \cases{ 1, &\quad if $l\in\Lambda
^c$,
\vspace*{5pt}\cr
\displaystyle K\sum_{l\in A\not\subset\Lambda }\|\Phi_A\|, &\quad if $l\in \Lambda \setminus i$,}
\end{eqnarray*}
where $K$ some constant and
again we used $|e^x-e^y|\leq|x-y|e^{\max(|x|,|y|)}$.
Notice, for any fixed $l$ when $\Lambda $ tends to ${\mathbb Z}^d$,
because of the
absolute summability of the Hamiltonian, this goes to zero.
Further we want to estimate
%
\begin{eqnarray}
\label{RelativeEntropyLocalErrorOrder6} &&\sup_{\eta_{\Lambda ^c}\omega _\Lambda } \bigl|\mu_{G\setminus
i}[
\eta_{\Lambda ^c}\omega _{\Lambda
\setminus i}] \bigl(\psi_1^{\omega _i}
\bigr) -\tilde\mu_{G\setminus i}[\omega _{\Lambda \setminus i}] \bigl(\tilde\psi_1^{\omega _i} \bigr) \bigr|\quad\mbox{and}
\nonumber\\[-8pt]\\[-8pt]
&&\sup_{\eta_{\Lambda ^c}\omega _\Lambda } \bigl|\mu_{G\setminus
i}[\eta_{\Lambda ^c}\omega
_{\Lambda
\setminus i}] \bigl(\psi_2^{\omega _i} \bigr) -\tilde\mu_{G\setminus i}[\omega _{\Lambda \setminus i}] \bigl(\tilde\psi_2^{\omega _i}
\bigr) \bigr|.\nonumber
\end{eqnarray}
We do this for both terms simultaneously by just writing $\psi$ instead
of $\psi_1,\psi_2$.
%
\begin{eqnarray}
\label{RelativeEntropyLocalErrorOrder7} &&\sup_{\eta_{\Lambda ^c}\omega _\Lambda } \bigl|\mu_{G\setminus
i}[
\eta_{\Lambda ^c}\omega _{\Lambda
\setminus i}] \bigl(\psi^{\omega _i} \bigr) -
\tilde\mu_{G\setminus i}[\omega _{\Lambda \setminus i}] \bigl(\tilde\psi^{\omega _i} \bigr) \bigr|\nonumber
\\
&&\qquad \leq \sup_{\eta_{\Lambda ^c}\omega _\Lambda } \bigl|\mu _{G\setminus i}[\eta_{\Lambda ^c}
\omega _{\Lambda
\setminus i}] \bigl(\tilde\psi^{\omega _i} \bigr) -\tilde\mu_{G\setminus i}[\omega _{\Lambda \setminus i}] \bigl(\tilde\psi^{\omega _i}
\bigr) \bigr|
\\
&&\quad\qquad{} +\sup_{\eta_{\Lambda ^c}\omega _\Lambda } \bigl|\mu _{G\setminus i}[
\eta_{\Lambda
^c}\omega _{\Lambda \setminus i}] \bigl(\psi^{\omega _i} \bigr) -
\mu_{G\setminus i}[\eta_{\Lambda ^c}\omega _{\Lambda \setminus
i}] \bigl(\tilde\psi^{\omega
_i} \bigr) \bigr|.\nonumber
\end{eqnarray}
For the second part in (\ref{RelativeEntropyLocalErrorOrder7})
we have
\begin{eqnarray*}
\label{RelativeEntropyLocalErrorOrder8} &&\sup_{\eta_{\Lambda ^c},\omega _\Lambda } \bigl|\mu_{G\setminus
i}[
\eta_{\Lambda ^c}\omega _{\Lambda
\setminus i}] \bigl(\psi^{\omega _i} \bigr) -
\mu_{G\setminus i}[\eta_{\Lambda ^c}\omega _{\Lambda \setminus
i}] \bigl(\tilde\psi^{\omega
_i} \bigr) \bigr|
\\
&&\qquad \leq \sup_{\omega _i}\big\Vert\psi^{\omega _i}-\tilde\psi^{\omega_i}\big\Vert
\\
&&\qquad \leq K\sum_{i\in A\not\subset\Lambda }\|\Phi_A\|,
\end{eqnarray*}
which tends to zero as $\Lambda \nearrow{\mathbb Z}^d$ by the absolute
summability of
the potential. In particular there exists a radius $r\in{\mathbb N}$ such that
$\sup_{i\in\Lambda _{n-r}}\sum_{i\in A\not\subset\Lambda _n}\|\Phi_A\|<\varepsilon $ for
all centered cubes $\Lambda _n$ such that $n-r\geq0$. Hence
%
\begin{eqnarray}
\label{RelativeEntropyLocalErrorOrder9} \frac{1}{|\Lambda _n|}\sum_{i\in\Lambda _n}\sum
_{i\in A\not
\subset\Lambda _n}\|\Phi_A\| 
<\varepsilon +\|H_0\|\!
\frac{|\Lambda _n\setminus\Lambda
_{n-r}|}{|\Lambda _n|},
\end{eqnarray}
where the RHS becomes arbitrarily small as $n\to\infty$.

Let us look at the first part of (\ref
{RelativeEntropyLocalErrorOrder7}) and use the Dobrushin
comparison theorem, which states
%
\begin{eqnarray}
\label{RelativeEntropyLocalErrorOrder10}
&& \sup_{\eta_{\Lambda ^c}\omega _\Lambda } \bigl|\mu_{G\setminus
i}[
\eta_{\Lambda ^c}\omega _{\Lambda
\setminus i}] \bigl(\tilde\psi^{\omega _i}
\bigr) -\tilde\mu_{G\setminus i}[\omega _{\Lambda \setminus i}] \bigl(\tilde\psi^{\omega _i} \bigr) \bigr|\nonumber
\\
&&\qquad \leq\sup_{\eta_{\Lambda ^c}\omega _\Lambda }\sum_{k\neq i,l\neq
i}
\delta_k \bigl(\tilde\psi ^{\omega _i} \bigr)D_{kl}
\bigl(\gamma ^{\eta_{\Lambda ^c}\omega _{\Lambda
\setminus i}} \bigr)b_l
\\
&&\qquad \leq\sup_{\omega _i}\sum_{k\in\Lambda \setminus i}\sum
_{l\in
\Lambda ^c}\delta_k \bigl(\tilde\psi^{\omega _i} \bigr)\widebar D_{kl}+\sup_{\omega _i}
\sum_{k\in
\Lambda \setminus i}\sum_{l\in\Lambda \setminus i}
\delta_k \bigl(\tilde\psi^{\omega
_i} \bigr)\widebar
D_{kl}b_l.\nonumber
\end{eqnarray}
As for the second term on the RHS of (\ref
{RelativeEntropyLocalErrorOrder10}), we have
\begin{eqnarray*}
&& \sum_{i\in\Lambda } \sum
_{k\in\Lambda \setminus i} \sum_{l\in
\Lambda \setminus i}\sup
_{\omega
_i}\delta_k \bigl(\tilde\psi^{\omega _i}
\bigr)\widebar D_{kl}b_l 
\\
&&\qquad \leq
\biggl(\sum_{l\in\Lambda \setminus i}b_l \biggr) \biggl(\sup
_{l}\sum_{k\in
\Lambda \setminus
i}\widebar
D_{kl} \biggr) \biggl(\sup_k\sum
_{i\in\Lambda }\sup_{\omega
_i}\delta_k
\bigl(\tilde\psi ^{\omega _i} \bigr) \biggr)
\\
&&\qquad \leq K\sum_{l\in\Lambda \setminus i}b_l=o\bigl(|\Lambda |\bigr).
\end{eqnarray*}
Indeed we have for all $k$,
%
\begin{equation}
\label{RelativeEntropyLocalErrorOrder11} \sum_{i\in\Lambda }\sup_{\omega _i}
\delta_k \bigl(\tilde\psi^{\omega _i} \bigr) 
\leq e^C\sum
_{i\in\Lambda }\sum_{\{i,k\}\subset A\subset\Lambda
}
\delta_k(\Phi_A) 
\leq
e^C\sum_{0\in A}\|\!|\Phi_A\|\!|<\infty
\end{equation}
and also $\sum_{k\in\Lambda \setminus i}\widebar D_{kl}\leq\sum_{k}\widebar
D_{kl}=\sum_{k}\widebar D_{0,l-k}=\sum_{j}\widebar D_{0,j}<\infty$
for all $l$. Finally,
\[
\label{RelativeEntropyLocalErrorOrder12} \frac{1}{|\Lambda |}\sum_{l\in\Lambda \setminus i}b_l
\leq\frac
{K}{|\Lambda |}\sum_{l\in\Lambda
}\sum
_{l\in A\not\subset\Lambda }\|\Phi_A\|\to0\qquad\mbox{as }
\Lambda \nearrow{\mathbb Z}^d 
\]
by the Ces\`{a}ro argument as in (\ref
{RelativeEntropyLocalErrorOrder9}). Let us consider for the first
term on the RHS of~(\ref{RelativeEntropyLocalErrorOrder10})
%
\begin{eqnarray}
\label{RelativeEntropyLocalErrorOrder13} \sum_{i\in\Lambda }\sum
_{k\in\Lambda \setminus i}\sum_{l\in
\Lambda ^c}\sup
_{\omega _i}\delta _k \bigl(\tilde\psi^{\omega _i}
\bigr)\widebar D_{kl}&=&\sum_{i\in\Lambda
}\sum
_{l\in\Lambda ^c}\sum_{k\in\Lambda \setminus i}\sup
_{\omega _i}\delta_k \bigl(\tilde\psi ^{\omega _i}
\bigr)\widebar D_{kl}.
\end{eqnarray}
Notice we assume the model to have the exponential decay property (\ref
{expdecay}) with increasing translation invariant semi-metric
$\varrho $ on
$G$ and again summability of the potential in the triple-norm. Thus for
all $i$ and $l$ by the triangle inequality
%
\begin{eqnarray}
\label{RelativeEntropyLocalErrorOrder14}
&& \sum_{k\in\Lambda \setminus i}e^{\varrho (i,l)}\sup
_{\omega
_i}\delta_k \bigl(\tilde\psi^{\omega
_i}
\bigr)\widebar D_{kl}\nonumber
\nonumber\\[-8pt]\\[-8pt]
&&\qquad \leq  \biggl(\sup
_{i\in\Lambda }\sum_{k}e^{\varrho (i,k)}
\sup_{\omega
_i}\delta_k \bigl(\tilde\psi
^{\omega _i} \bigr) \biggr) \Bigl(\sup_{l,k}e^{\varrho (k,l)}
\widebar D_{kl} \Bigr) 
\leq 
\widetilde C.\nonumber
\end{eqnarray}
%
Hence we can write
%
\begin{eqnarray}
\label{RelativeEntropyLocalErrorOrder15} \sum_{i\in\Lambda }\sum
_{l\in\Lambda ^c}\sum_{k\in\Lambda
\setminus i}\sup
_{\omega _i}\delta _k \bigl(\tilde\psi^{\omega _i}
\bigr)\widebar D_{kl} 
&=&\widetilde C \sum
_{j\in{\mathbb Z}^d}e^{-\varrho (0,j)}\#\{i\in\Lambda\dvtx i+j\notin\Lambda \},
\end{eqnarray}
which again tends to infinity slower than $|\Lambda |$.
\end{pf}
%
\begin{lem}\label{EstimateSecondErrorTerm}
$\sup_{\omega }\sum_{i\in\Lambda }B(\omega,i^+)=o(|\Lambda |)$.
\end{lem}
\begin{pf}For the next error term in (\ref
{RelativeEntropyLocalError}) we have
%
\begin{eqnarray}
\label{RelativeEntropyLocalErrorOrderSecond1} B \bigl(\omega,i^+ \bigr) 
&=&\biggl|
\log \frac{\tilde\mu_{G\setminus i}[\omega _{\Lambda
\setminus i}](\lambda
(e^{-\widetilde H_i}1_{\omega ^i_i}))}{\mu_{G\setminus i}[\zeta
_{\Lambda ^c}\omega _{\Lambda
\setminus i}](\lambda(e^{-H_i}1_{\omega ^i_i}))}-\log\frac{\tilde\mu
_{G\setminus i}[\omega _{\Lambda \setminus i}](\lambda(e^{-\widetilde
H_i}1_{\omega
_i}))}{\mu_{G\setminus i}[\zeta_{\Lambda ^c}\omega _{\Lambda
\setminus i}](\lambda
(e^{-H_i}1_{\omega _i}))}\biggr|\nonumber
\\
&\leq&\frac{|\tilde\mu_{G\setminus i}[\omega _{\Lambda \setminus
i}](\lambda
(e^{-\widetilde H_i}1_{\omega ^i_i}))
-\mu_{G\setminus i}[\zeta_{\Lambda ^c}\omega _{\Lambda \setminus
i}](\lambda
(e^{-H_i}1_{\omega ^i_i}))|}{\mu_{G\setminus i}[\zeta_{\Lambda
^c}\omega _{\Lambda
\setminus i}](\lambda(e^{-H_i}1_{\omega ^i_i}))}
\\
&&{} +\frac{|\tilde\mu_{G\setminus i}[\omega _{\Lambda
\setminus
i}](\lambda(e^{-\widetilde H_i}1_{\omega _i}))
-\mu_{G\setminus i}[\zeta_{\Lambda ^c}\omega _{\Lambda \setminus
i}](\lambda
(e^{-H_i}1_{\omega _i}))|}{\mu_{G\setminus i}[\zeta_{\Lambda
^c}\omega _{\Lambda \setminus
i}](\lambda(e^{-H_i}1_{\omega _i}))},\nonumber
\end{eqnarray}
where we assumed $\tilde\mu_{G\setminus i}[\omega _{\Lambda
\setminus i}](\lambda
(e^{-\widetilde H_i}1_{\omega ^i_i}))\geq\mu_{G\setminus i}[\zeta
_{\Lambda ^c}\omega _{\Lambda
\setminus i}](\lambda(e^{-H_i}1_{\omega ^i_i}))$ and $\tilde\mu
_{G\setminus i}[\omega _{\Lambda \setminus i}](\lambda(e^{-\widetilde
H_i}1_{\omega
_i}))\geq\mu_{G\setminus i}[\zeta_{\Lambda ^c}\omega _{\Lambda
\setminus i}](\lambda
(e^{-H_i}1_{\omega _i}))$. In this case, as well as in all other
cases, we
can follow the exact same arguments as before and get
%
\begin{equation}
\label{RelativeEntropyLocalErrorOrderSecond2} \frac{1}{|\Lambda |}\sup_\omega \sum
_{i\in\Lambda }B \bigl(\omega,i^+ \bigr)\to0\qquad\mbox{as }\Lambda
\nearrow {\mathbb Z}^d.
\end{equation}\upqed
\end{pf}

Given the fact that $\sup_i\sum_\omega \Gamma _L(\omega,i^+)<\infty$ and
\[
\sup_{i\in\Lambda }\sum_{\omega }\nu(1_{\omega })\bigl|V_{\widetilde
L_\Lambda }(\omega,i^+)\bigr|
\leq K\sup_{i\in\Lambda }\sum_{\omega }\nu(1_{\omega })=K<\infty,
\]
we have by now verified that the second summand in the last line of
(\ref{RelativeEntropyLocalError}) is indeed $o(|\Lambda |)$.
The first summand in the last line of (\ref
{RelativeEntropyLocalError}) requires some extra care. We prepare by writing
%
\begin{eqnarray}
\label{LastErrorterm1} \qquad\sum_{\omega, i\in\Lambda }A \bigl(\omega,i^+
\bigr) \nu(1_\omega )\biggl|\log \frac{\nu(1_{\omega ^i})}{\nu(1_{\omega
})}\biggr|&\leq&\sum
_{i\in\Lambda }\sup_{\omega }A \bigl(\omega,i^+ \bigr)
\sum_{\omega }\nu(1_{\omega })\biggl|\log
\frac
{\nu(1_{\omega ^i})}{\nu(1_{\omega })}\biggr|. 
\end{eqnarray}
%
%
%
%
%
The next step is then to show boundary order of the RHS of (\ref
{LastErrorterm1}), in other words to show the following lemma.
%
\begin{lem}\label{EstimateThirdErrorTerm}
$\sum_{i\in\Lambda }\sup_{\omega }
A(\omega,i^+)\sum_{\omega }\nu(1_{\omega })
|\log\frac{\nu(1_{\omega ^i})}{\nu(1_{\omega })}|
=o(|\Lambda |)$.
\end{lem}

Notice, since the rates are bounded from below away from zero and
bounded from above, that is, $e^{-2\Vert H_0\Vert}\leq c_L(\omega,\omega ^i)\leq
\tilde ke^{2\Vert H_0\Vert}$, $e^{-\Vert H_0\Vert}\leq c_K(\omega,\omega ^i)\leq
\tilde ke^{\Vert H_0\Vert}$, $e^{-\Vert H_0\Vert}\leq c_K(\omega,\omega
^{i-})\leq\tilde ke^{\Vert H_0\Vert}$ and $\nu$ is invariant, that is,
$0=\int(L+\alpha K)1_{\omega _\Lambda }\,d\nu$, we have for all
$\omega \in\{1, \ldots, q\}^{\Lambda
}$ [after separation of the terms in this equation proportional to $\nu
(1_\omega )$ from $\nu(1_{\omega ^i})$ and $\nu(1_{\omega ^{-i}})$],
%
\begin{eqnarray}
\label{InvariantMeasureEstimate1} \widetilde K&\geq&\frac{1}{|\Lambda |}\sum
_{i\in\Lambda } \biggl[\frac{\nu
(1_{\omega ^{i-}})}{\nu
(1_{\omega })}+\frac{\nu(1_{\omega ^{i}})}{\nu(1_{\omega })} \biggr]
\geq \frac{1}{|\Lambda |}\sum_{i\in\Lambda }\frac{\nu(1_{\omega })}{\nu(1_{\omega ^i})}
\nonumber\\[-8pt]\\[-8pt]
&=&
\frac
{1}{|\Lambda |}\sum_{i\in\Lambda
}\frac{\nu(\omega _i| \omega _{\Lambda \setminus i})}{\nu
(\omega ^i_i| \omega _{\Lambda \setminus i})}.\nonumber
\end{eqnarray}
To control possibly small arguments of the logarithm, we need to bound
\mbox{$\nu$-}probabi\-lities from below. For this the following lemma will be useful.
%
\begin{lem}\label{InvariantMeasureEstimate0}
Let $\nu\in\mathcal{P}({\Omega}')$ be translation-invariant and
invariant for
the joint dynamics. There exists a constant $\widehat K$ such that for all
finite sets $\Delta $ we have
\[
\int\nu(d\omega )\frac{1}{\sqrt{\nu(\omega _0^0| \omega
_{\Delta \setminus 0})}}<\widehat K.
\]
\end{lem}
\begin{pf*}{Proof of Lemma~\ref{InvariantMeasureEstimate0}}
By the
Jensen inequality, it suffices to show this for centered cubes $\Delta $.
%
%
Let us consider the $\nu$-expectation of our essential estimate~(\ref
{InvariantMeasureEstimate1}) and apply Jensen's inequality to obtain
\begin{eqnarray*}
\label{InvariantMeasureEstimate2} \widetilde K&\geq&\frac{1}{|\Lambda |}\sum
_{i\in\Lambda } \int\nu (d\omega )\frac{\nu(\omega _i| \omega
_{\Lambda \setminus i})}{\nu(\omega ^i_i| \omega _{\Lambda
\setminus i})} 
\\
&\geq& \frac{1}{|\Lambda |}\sum_{i\in\Lambda }\sum
_{k=1}^q \biggl[\int\nu (d\omega )
\frac{\nu(k| \omega
_{\Lambda \setminus i})}{\sqrt{\nu((k+1)| \omega _{\Lambda
\setminus i})}} \biggr]^2
\\
&\geq&\frac{1}{|\Lambda |}\sum_{i\in\Lambda }\min
_l\int\nu (d\omega )\frac{\nu(l| \omega _{\Lambda
\setminus i})}{\sqrt{\nu((l+1)| \omega _{\Lambda \setminus
i})}}\sum
_{k=1}^q\int\nu(d\omega )\frac
{\nu(k| \omega _{\Lambda \setminus i})}{\sqrt{\nu((k+1)| \omega
_{\Lambda \setminus i})}}
\\
&\geq&\frac{1}{|\Lambda |}\sum_{i\in\Lambda }\min
_l\int\nu (d\omega )\nu(l| \omega _{\Lambda \setminus
i})\sum
_{k=1}^q\int\nu(d\omega )
\frac{1}{\sqrt{\nu(\omega _i^i|
\omega _{\Lambda \setminus i})}}
\\
&\geq&\frac{1}{|\Lambda |}\sum_{i\in\Lambda }\min
_l\nu(l)\sum_{k=1}^q
\int\nu(d\omega )\frac{1}{\sqrt{\nu(\omega _i^i| \omega _{\Lambda \setminus
i})}}
\\
&\geq&\frac{\varepsilon _0}{|\Lambda |}\sum
_{i\in\Lambda }\int\nu(d\omega )\frac{1}{\sqrt{\nu(\omega
_i^i| \omega _{\Lambda \setminus i})}},
\end{eqnarray*}
where we used $\min_l\nu(l)\geq\varepsilon _0>0$.

%
\begin{rem}
In fact, take $\Lambda =\{0\}$ in (\ref
{InvariantMeasureEstimate1}), and
then we have $\frac{\nu(k)}{\nu(k+1)}\leq\widetilde K$ and hence
$\nu(k)=\frac{1}{\sum_{l=1}^q (\nu(l)/\nu(k))}
\geq\frac{1}{1+\widetilde K+\widetilde K^2+\cdots+\widetilde
K^{(q-1)}}=\frac{1}{\sum_{l=0}^{q-1}\widetilde K^l}=:\varepsilon _0$.
\end{rem}

Consider $\Lambda _n:=[-n,n]^d$ and $m:=n-\lfloor n^\kappa \rfloor$
with $\kappa \in
(0,1)$. Then the above inequality
can be further estimated by
\begin{eqnarray*}
\label{InvariantMeasureEstimate3} \frac{\widetilde K}{\varepsilon _0} 
&\geq&
\frac{1}{|\Lambda _n|}\sum_{i\in\Lambda _m}\int\nu(d\omega )
\frac{1}{\sqrt{\nu(\omega
_i^i| \omega _{\Lambda _n\setminus i})}}\geq\frac{1}{|\Lambda
_n|}\sum_{i\in\Lambda _m}
\int\nu(d\omega )\frac{1}{\sqrt{\nu(\omega _i^i| \omega _{(\Delta +i)\setminus i})}},
\end{eqnarray*}
where $\Delta $ is the largest centered cube such that for all $i\in
\Lambda _m$ we
have $\Delta +i\subset\Lambda $. We used the conditional Jensen
inequality in the
last line. Because of translation-invariance we have $\frac{\widetilde
K}{\varepsilon
_0}\geq\frac{|\Lambda _m|}{|\Lambda _n|}\int\nu(d\omega )\frac
{1}{\sqrt{\nu(\omega _0^0| \omega
_{\Delta \setminus 0})}}$. Since $\frac{|\Lambda _m|}{|\Lambda _n|}
\to1$ for $n\to\infty$ and $\Lambda _n\setminus \Lambda _m$
allows $\Delta $ to become
arbitrarily large, the result of Lemma~\ref{InvariantMeasureEstimate0} follows.
\end{pf*}

\begin{pf*}{Proof of Lemma~\ref{EstimateThirdErrorTerm}}
Consider centered cubes $\Lambda _n$ of side-length \mbox{$2n+1$}, and write
$\partial_{r}(\Lambda _n):=\{i\in\Lambda\dvtx d(i,\Lambda ^c)=r\}$
where $d(\cdot ,\cdot)$ is the
uniform norm. We show for $i\in\partial_{r}(\Lambda _n)$
%
\begin{equation}
\label{LastErrorterm2} \sup_{\omega }A \bigl(\omega,i^+ \bigr)\leq f(r)
\end{equation}
with $\lim_{r\to\infty}f(r)=0$. Indeed, let us look at (\ref
{RelativeEntropyLocalErrorOrder7}) again. We can estimate the
second part by
%
\begin{eqnarray}
\label{EstimatesOnA1} \qquad \sup_{i\in\partial_r(\Lambda _n)}\sup_{\omega _i}\big\Vert
\psi ^{\omega _i}-\tilde\psi^{\omega
_i}\big\Vert&\leq&\sup
_{i\in\partial_r(\Lambda _n)}K\sum_{i\in A\not
\subset\Lambda
_n}\Vert
\Phi_A\Vert 
\leq K\sum
_{0\in A\not\subset\Lambda _r}\Vert\Phi_A\Vert,
\end{eqnarray}
which goes to zero as $r$ tends to infinity.
For the first part of (\ref{RelativeEntropyLocalErrorOrder7}) we have
%
\begin{eqnarray}
\label{EstimatesOnA2} && \sup_{\eta_{\Lambda _n^c}\omega _\Lambda } \bigl|\mu _{G\setminus i}[
\eta_{\Lambda _n^c}\omega _{\Lambda
_n\setminus i}] \bigl(\tilde\psi^{\omega _i}
\bigr) -\tilde\mu_{G\setminus i}[\omega _{\Lambda _n\setminus i}] \bigl(\tilde\psi^{\omega _i} \bigr) \bigr|\nonumber
\\
&&\qquad \leq\sup_{\omega _i}\sum
_{k\neq i,l\neq i}\delta_k \bigl(\tilde\psi ^{\omega
_i}
\bigr)\widebar D_{kl}b_l
\\
&&\qquad \leq \sup_{\omega _i}\sum_{k\in\Lambda _n\setminus i}\sum
_{l\in
\Lambda _n^c}\delta _k \bigl(\tilde\psi^{\omega _i} \bigr)\widebar D_{kl}+\sup_{\omega _i}
\sum_{k\in\Lambda
_n\setminus i}\sum_{l\in\Lambda _n\setminus i}
\delta_k \bigl(\tilde\psi ^{\omega
_i} \bigr)\widebar
D_{kl}b_l.\nonumber
\end{eqnarray}
%
%
%
By looking at (\ref{RelativeEntropyLocalErrorOrder14}) we notice
for $i\in\partial_r(\Lambda _n)$
\begin{eqnarray*}
\label{EstimatesOnA4} \sup_{\omega _i}\sum_{k\in\Lambda _n\setminus i}
\sum_{l\in
\Lambda _n^c}\delta_k \bigl(\tilde\psi^{\omega _i} \bigr)\widebar D_{kl}&\leq&\widetilde C\sum
_{l\in\Lambda
_n^c}e^{-\varrho
(i,l)}\leq\widetilde C\sum
_{l\in(\Lambda _n-i)^c}e^{-\varrho
(0,l)}
\\
&\leq&\widetilde C\sum
_{l\in\Lambda _r^c}e^{-\varrho (0,l)},
\end{eqnarray*}
which goes to zero as $r$ tends to infinity. Let us define
$n_l:=d(l,\Lambda
_n^c)$ for the distance between the site $l$ and $\Lambda _n^c$. For the
other summand in (\ref{EstimatesOnA2}) we have
%
\begin{eqnarray*}
\label{EstimatesOnA5} \sup_{\omega _i}\sum_{k\in\Lambda _n\setminus i}
\sum_{l\in
\Lambda _n\setminus i}\delta _k \bigl(\tilde\psi^{\omega _i} \bigr)\widebar D_{kl}b_l
&\leq&\sum_{k}\sup_{\omega _0}
\delta_k \bigl(\tilde\psi^{\omega
_0} \bigr)\sum
_{l}\widebar D_{0,l}\sum
_{0\in A\not\subset\Lambda _{n_{i+k+l}}}\Vert\Phi _A\Vert.
\end{eqnarray*}
Notice if $n_i\to\infty$ then $n_{i+l+k}\to\infty$ for every fixed
$l,k$. In particular since $\sum_{l}\widebar D_{0,l}<\infty$ we have
for $n_i\to\infty$
\[
\label{EstimatesOnA6} \sum_{l}\widebar
D_{0,l}\sum_{0\in A\not\subset\Lambda
_{n_{l+k+l}}}\Vert\Phi
_A\Vert\to0
\]
%
%
%
by the dominated convergence theorem.
Similar to (\ref{RelativeEntropyLocalErrorOrder11}) we have
%
\begin{equation}
\label{REstimatesOnA8} \sum_{k}\delta_k
\bigl(\tilde\psi^{\omega _0} \bigr) 
\leq e^C\sum_{k}\sum
_{\{i,k\}\subset A}\delta_k(\Phi_A)
\leq e^C\sum_{0\in A}\|\!|\Phi_A\|\!|<\infty
\end{equation}
and thus for $n_i\to\infty$ we can conclude again with the dominated
convergence theorem
\[
\label{EstimatesOnA9} \sum_{k}\sup_{\omega _0}
\delta_k \bigl(\tilde\psi^{\omega _0} \bigr)\sum
_{l}\widebar D_{0,l}\sum
_{0\in A\not\subset\Lambda _{n_{l+k+l}}}\Vert\Phi _A\Vert\to0.
\]

Since $-\log(x)<\frac{1}{\sqrt{x}}<\frac{1}{x}$ on $(0,1]$, we
finally have
%
\begin{eqnarray}
\label{LastErrorterm3} && \sum_{r=1}^n\sum
_{i\in\partial_{r}(\Lambda _n)}\sup_{\omega
}A \bigl(\omega,i^+ \bigr)\int\nu(d\omega )\biggl|\log\frac{\nu(\omega ^i_i|\omega _{\Lambda _n\setminus
i})}{\nu(\omega _i|\omega _{\Lambda _n\setminus i})}\biggr|\nonumber
\\
&&\qquad \leq\sum_{r=1}^n
f(r)\sum_{i\in\partial_{r}(\Lambda _n)} \biggl[\int -\log\nu(\omega
_i|\omega _{\Lambda _n\setminus i})\nu(d\omega )\nonumber
\\
&&\hspace*{108pt}{} +\int-\log\nu \bigl(\omega
^i_i|\omega _{\Lambda _n\setminus i} \bigr)\nu(d\omega ) \biggr]
\nonumber\\[-8pt]\\[-8pt]
&&\qquad \leq\sum_{r=1}^n f(r)\sum
_{i\in\partial_{r}(\Lambda _n)} \biggl[\int \nu(d\omega )\frac
{1}{\nu(\omega _i|\omega _{\Lambda \setminus i})}+\int\nu (d
\omega )\frac{1}{\sqrt{\nu(\omega _i^i|\omega _{\Lambda
_n\setminus i})}} \biggr]\nonumber
\\
&&\qquad =\sum_{r=1}^n f(r)\sum
_{i\in\partial_{r}(\Lambda _n)} \Biggl[\sum_{k=1}^q
\int\nu(d\omega )\frac{\nu(k|\omega _{\Lambda \setminus i})}{\nu(k|\omega
_{\Lambda \setminus i})}+\int\nu(d\omega )\frac
{1}{\sqrt{\nu(\omega _i^i|\omega _{\Lambda _n\setminus i})}} \Biggr]\nonumber
\\
&&\qquad \leq\sum_{r=1}^n f(r)\sum
_{i\in\partial_{r}(\Lambda _n)}[q+\widehat K]=K\sum_{r=1}^n
f(r)\bigl|\partial_{r}(\Lambda _n)\bigr|,\nonumber
\end{eqnarray}
where we used the last lemma in the last line. Since $f(r)\to0$ for
$r\to\infty$ there exists a $R\in{\mathbb N}$ such that for all $r\geq
R$ we
have $f(r)<\varepsilon $, and hence for large $n$
\begin{eqnarray*}
\label{LastErrorterm5}
&& \frac{1}{|\Lambda _n|}\sum_{r=1}^n
f(r)\bigl|\partial_{r}(\Lambda _n)\bigr|
\\
&&\qquad =\frac{1}{|\Lambda
_n|}
\Biggl(\sum_{r=R+1}^n f(r)\bigl|
\partial_{r}(\Lambda _n)\bigr|+\sum
_{r=1}^R f(r)\bigl|\partial _{r}(\Lambda
_n)\bigr| \Biggr)
\\
&&\qquad \leq \varepsilon +K\frac{|\partial(\Lambda _n)|}{|\Lambda _n|},
\end{eqnarray*}
where the second summand goes to zero as $n$ tends to infinity.
\end{pf*}

%
%
%
%
%
%
%
%
%
%
%
%
%
%
%
%
%

Together we see the combined error caused by the finite-volume
approximation vanishes from the point of view of difference between
time derivatives of relative entropy densities, that is, $\frac
{d}{dt}_{|_{t=0}}H_{\Lambda }(\nu_{t,L}| \gamma _\Lambda '(\cdot
|{\zeta}))-\frac
{d}{dt}_{|_{t=0}}H_{\Lambda }(\nu_{t,\widetilde L_\Lambda }|\break\tilde\gamma _\Lambda )=o(|\Lambda |)$.
For a translation-invariant measure $\nu$, that is also invariant w.r.t.
the joint dynamics, we have
\begin{eqnarray*}
\label{RelativeEntropyBigEstimate1} 0&=&\frac{d}{dt}_{|_{t=0}}H_{\Lambda } \bigl(
\nu_{t,L+\alpha K}|\gamma _\Lambda '(\cdot|{\zeta}) \bigr)
\\
&=&\frac{d}{dt}_{|_{t=0}}H_{\Lambda } \bigl(\nu_{t,L}|
\gamma _\Lambda '(\cdot|{\zeta}) \bigr)+\alpha
\frac
{d}{dt}_{|_{t=0}}H_{\Lambda } \bigl(\nu_{t,K}|\gamma
_\Lambda '(\cdot |{\zeta}) \bigr).
\end{eqnarray*}
Hence with the notation given in (\ref{DefinitionsHolleysArgument})
we can write
\begin{eqnarray*}
\label{RelativeEntropyBigEstimate2}
&&\alpha \sum_{i\in\Lambda } \bigl[
\kappa_{\Lambda,K} \bigl(i^+ \bigr)+\kappa _{\Lambda,K} \bigl(i^- \bigr)
\bigr]
\\
&&\qquad \leq \alpha \sum_{i\in\Lambda } \bigl[ \bigl(
\beta_{\Lambda,K} \bigl(i^+ \bigr)+\beta _{\Lambda,K} \bigl(i^- \bigr)
\bigr) \vartheta _{\Lambda,K}(i) \bigr]+2\frac{d}{dt}_{|_{t=0}}H_{\Lambda }(
\nu_{t,L}|\gamma _\Lambda ' \bigl(\cdot|({\zeta
}) \bigr)
\\
&&\qquad \leq \alpha \widehat C\sum_{i\in\Lambda }\sum
_{j\notin\Lambda
}\delta_{j} \bigl(c_K \bigl(
\cdot, \cdot ^i \bigr) \bigr)
\\
&&\quad\qquad{} +2\widetilde C\biggl|\frac{d}{dt}_{|_{t=0}}H_{\Lambda }
\bigl(\nu_{t,L}|\gamma _\Lambda '(\cdot |{\zeta})
\bigr)-\frac{d}{dt}_{|_{t=0}}H_{\Lambda }(\nu_{t,\widetilde
L_\Lambda }|
\tilde\gamma _\Lambda )\biggr|
\\
&&\qquad \leq\alpha \widehat Co\bigl(|\Lambda |\bigr)+2\widetilde C\widetilde Ko\bigl(|\Lambda |\bigr)=o\bigl(|\Lambda |\bigr),
\end{eqnarray*}
where in the second line we dropped the contribution of the
finite-volume part since it is only negative.

However, this estimate implies the single-site DLR equation and thus
$\nu$ must be Gibbs. This completes the proof of Proposition~\ref
{InvarianceIsGibbs}.
\end{pf}

\section*{Acknowledgments}
We thank C. Maes and A. van Enter for stimulating discussions at the workshop
``Dynamical Gibbs--non-Gibbs Transitions'' at Eurandom, and an anonymous
referee for useful comments.



\printaddresses

\end{document}